\PassOptionsToPackage{normalem,normalbf}{ulem}
\documentclass[AMA,STIX1COL]{WileyNJD-v2}
\usepackage[utf8]{inputenc}

\usepackage{amsmath, amsthm, bm}
\usepackage{cancel}

\usepackage{graphicx,xcolor}

\usepackage[draft]{changes}
\definechangesauthor[color=magenta]{NB}
\definechangesauthor[color=orange]{JB}
\definechangesauthor[color=green!75!black]{SV}

\articletype{Research Article}%

\received{26 January 2019}
\revised{6 June 2019}
\accepted{6 June 2019}

\renewcommand*{\doi}[1]{doi: \href{http://dx.doi.org/\detokenize{#1}}{\ttfamily\detokenize{#1}}}

\DeclareMathAlphabet{\mathbbmss}{U}{bbmss}{m}{n}

\newcommand{\fatvec}[1]{\bm{#1}}
\newcommand{\dint}[1]{\,\mathrm{d}{#1}}

\newcommand{\pptext}[2]{{\tfrac{\partial #1}{\partial #2}}}

\newcommand*{\tran}{^{\mkern-1.5mu\mathsf{T}}}

\newcommand{\tolwet}{\mathrm{TOL}_\mathrm{wet}}
\newcommand{\Tend}{T_\mathrm{end}}

\newcommand{\Cx}{\fatvec{x}}
\newcommand{\vu}{\fatvec{u}}

\begin{document}

\title{A limiter-based well-balanced discontinuous Galerkin method for shallow-water flows with wetting and drying: Triangular grids}

\author[uhh,cen]{Stefan Vater}
\author[ucd]{Nicole Beisiegel}
\author[uhh,cen]{Jörn Behrens}

\address[uhh]{\orgdiv{Department of Mathematics}, \orgname{Universität Hamburg}, \orgaddress{Bundesstraße 55, \state{20146 Hamburg}, \country{Germany}}}
\address[cen]{\orgdiv{CEN -- Center for Earth System Research and Sustainability}, \orgname{Universität Hamburg}, \orgaddress{Grindelberg 5, \state{20144 Hamburg}, \country{Germany}}}
\address[ucd]{\orgdiv{School of Mathematics \& Statistics}, \orgname{University College Dublin}, \orgaddress{Belfield, \state{Dublin 4}, \country{Ireland}}}

\authormark{VATER \textsc{et al.}}

\corres{Stefan Vater, Institute of Mathematics, Freie Universität Berlin, Arnimallee 6, 14195 Berlin, Germany. \email{vater@math.fu-berlin.de}}


\fundingInfo{%
\fundingAgency{Volkswagen Foundation} \fundingNumber{Grant: ASCETE}
\fundingAgency{European Union, FP7-ENV2013 6.4-3} \fundingNumber{Grant: 603839}
\fundingAgency{Science Foundation Ireland} \fundingNumber{Grant: 14/US/E3111}}

\abstract[Summary]{A novel wetting and drying treatment for second-order Runge-Kutta discontinuous Galerkin (RKDG2) methods solving the non-linear shallow water equations is proposed. It is developed for general conforming two-dimensional triangular meshes and utilizes a slope limiting strategy to accurately model inundation. The method features a non-destructive limiter, which concurrently meets the requirements for linear stability and wetting and drying. It further combines existing approaches for positivity preservation and well-balancing with an innovative velocity-based limiting of the momentum. This limiting controls spurious velocities in the vicinity of the wet/dry interface. It leads to a computationally stable and robust scheme -- even on unstructured grids -- and allows for large time steps in combination with explicit time integrators. The scheme comprises only one free parameter, to which it is not sensitive in terms of stability. A number of numerical test cases, ranging from analytical tests to near-realistic laboratory benchmarks, demonstrate the performance of the method for inundation applications. In particular, super-linear convergence, mass-conservation, well-balancedness, and stability are verified.}


\jnlcitation{\cname{%
\author{Vater S.},
\author{N. Beisiegel}, and
\author{J.  Behrens}}
\ctitle{A limiter-based well-balanced discontinuous Galerkin method for shallow-water flows with wetting and drying: Triangular grids}, submitted to \cjournal{International Journal for Numerical Methods in Fluids}, 
}

\maketitle


\section{Introduction} \label{sec:Introduction}

In order to successfully compute near- and onshore propagation of ocean waves, depth-integrated equations such as the shallow water equations are commonly employed. While these are derived under the assumption that vertical velocities are negligible, they efficiently model large scale horizontal flows and wave propagation with high accuracy. Computational problems occur in the coastal area, where the shoreline, theoretically defined as the line of zero water depth, represents a moving boundary condition, which must be considered in the numerical scheme. Here, it is essential to construct a computational method, which concurrently fulfills the physical requirements for accurate coastal modeling:
\begin{itemize}
  \item conservation of mass,
  \item exact representation of the shoreline (wetting and drying), and
  \item robust computation of perturbations from the steady state at rest (well-balancedness).
\end{itemize}
Although the moving shoreline can be incorporated into the numerical scheme by adjusting the computational domain, its implementation is difficult and often only simple flow configurations have been successfully considered with this approach \citep{Bates1999}. Only recently, an Arbitrary Lagrangian Eulerian (ALE) method together with a moving $r$-adaptive mesh was applied to more complex flow situations \citep{Arpaia2018}. Most commonly, an Eulerian approach is considered, where the mesh points are fixed and the numerical scheme is applied to all cells, regardless if they are wet or dry. Cells are flooded or run dry through the interaction with other cells, i.e., by fluxes.

In recent years discontinuous Galerkin (DG) methods have gained a lot of scientific interest within the geophysical fluid dynamics (GFD) modeling community. They combine a number of desirable properties important for coastal applications such as conservation of physical quantities, geometric flexibility and accuracy \citep{Bernard2007,Giraldo2008,Dawson2013}. While Godunov-type finite volume (FV) methods are considered one of the most comprehensive tools for hydrologic modeling of coastal inundation problems \citep{Toro2007,Kesserwani2014}, several studies have investigated possible merits of the more complex DG formulation. General comparisons between the two discretization techniques point out that FV methods need wide stencils to attain high order of accuracy, whereas DG discretizations remain compact \citep{Zhou2001,Zhang2005}. This feature of DG formulations makes them particularly viable for $p$-adaptivity and parallelization. On the other hand the linear stability limit (CFL condition) is more generous for FV methods. Concerning coastal modeling, \citet{Kesserwani2014} argue that the extra complexity associated with DG discretizations pays off by providing higher-quality solution behavior on very coarse meshes. Furthermore, DG methods better approximate the near-shore velocity in certain situations \citep{Kesserwani2014,Vater2017}.

While there is theoretically no barrier to extend the DG method to higher-order accuracy, several practical aspects impede this endeavor. For the time discretization appropriate integration schemes -- such as strong stability preserving methods -- have to be applied, which maintain crucial properties of the governing equations. Such methods can be hard to derive for high-order methods. In many practical applications data sets expose a large amount of roughness leading to small-scale variations in high order methods requiring appropriate filtering with associated order reduction. With respect to coastal inundation, no high-order convergence can be expected at the wet-dry interface due to the non-differentiable transition from wet to dry. Therefore, a second-order DG model seems to be a reasonable choice for practical coastal modeling \cite{Kesserwani2014}. This is also reflected in the DG literature, where most wetting and drying treatments are build on top of a second-order accurate DG discretization \citep{Gourgue2009,Kesserwani2010a,Bunya2009,Ern2008}, and only few go beyond formal second-order accuracy (e.g., \citet{Xing2013}).

Several concepts addressing the above mentioned physical requirements for wetting and drying in a DG framework proved to be practical. To guarantee positivity of the fluid depth and conserve mass at the same time, various authors proposed a redistribution of mass within each cell \citep{Bunya2009,Ern2008,Xing2013}. This reduces the problem of positivity preservation to only requiring positivity in the mean, which can be guaranteed by a CFL time step restriction \citep{Xing2013}. Since the DG discretization may not exactly resolve the wet/dry interface, artificial gradients can occur in the surface elevation, which render the method unbalanced. To preserve a discrete steady state at rest in this case, it was proposed to neglect the gravity terms in such cells \citep{Bunya2009,Gourgue2009}. To retain a well-defined velocity computation, which is usually not a primary variable and calculated through the quotient of momentum and fluid depth, a so-called thin-layer approach was introduced, where a point is considered dry, if the fluid depth drops below a given tolerance \citep{Bunya2009}. Using this tolerance, the velocity can be set directly to zero in such a situation and the problem of possibly dividing by a zero fluid depth is circumvented.

Although there are other approaches to deal with the wetting and drying problem \citep{Gourgue2009,Kaernae2011}, the most common procedure is based on slope-modification techniques \citep{Ern2008,Bunya2009,Kesserwani2010a,Xing2013}. Based on the works of \citet{Cockburn1998}, a generalized Minmod total variation bounded (TVB) limiter is usually combined with the wetting and drying treatment to guarantee linear stability and prevent oscillations in case of discontinuities. However, several authors point out that this slope limiting and the handling of wetting and drying may activate each other, such that their concurrent use can lead to instabilities. This conflict is circumvented by applying the TVB limiter only in those cells, where the wetting and drying algorithm is not activated. This is the starting point of the current study, in which we base our wetting and drying treatment on Barth/Jespersen-type limiters \citep{Barth1989}. To the authors' experience such limiters do not severely alter a discrete state at rest and small perturbations around it, when limiting in surface elevation.

The application of the new non-destructive limiter leads to another problem, which does not seem to be as prominent for the TVB limiter. Because both, fluid depth and momentum diminish close to the shoreline, the quotient of both -- representing the velocity -- becomes numerically ill-conditioned. This may lead to spurious errors in the velocity values and result in severe time step restrictions induced by the CFL stability condition inherent in the discretization of hyperbolic problems. In order to solve this issue, our approach incorporates the velocity field into the limiting procedure for the momentum. The idea is borrowed from FV methods, where interface values are often reconstructed from other variables than the primary prognostic variables (see \citet{Leer2006} and references therein). Here, we develop an approach for DG methods to modify the primary variables based on other secondary variables. We stress that our scheme maintains common time step restrictions unlike implicit methods, such as in \citet{Meister2014}.

The aim of this work is to introduce a new approach that addresses all of the above issues by combining existing with novel strategies for a robust, efficient, and accurate inundation scheme for explicit DG computations on general triangular grids. In this course, a previously developed one-dimensional limiter \citep{Vater2015} is generalized to the two-dimensional case. Here, we rely on a nodal DG formulation along the lines of \citet{Giraldo2008}, where we work with geometrical nodes and basis function expansions defined by nodal values. To preserve positivity, we adopt the ``positive-depth operator'' from \citet{Bunya2009}, but only applied to the fluid depth. Furthermore, cancellation of gravity terms is applied in cells adjacent to the wet/dry interface. The presented method is based on limiting total fluid height $H=h+b$, which is the sum of fluid depth $h$ and bathymetry $b$, and velocity in the momentum-based flux computation. This approach is based on the original idea of hydrostatic reconstruction for FV methods \citep{Audusse2004,Zhou2001,Noelle2006}. Two Barth/Jespersen-type limiters are employed -- the original edge-based version \citep{Barth1989} and a modified vertex-based development \citep{kuzmin_vertex-based_2010}, the latter being particularly suitable for triangular grids. We find only minor differences between these options and utilize them for computational efficiency reasons. Although the limiter from \citet{kuzmin_vertex-based_2010} also works for higher than second-order accuracy, for the reasons given above we stick to piecewise linear basis functions, and leave the extension of our approach to higher-order for future research. A set of six commonly used test cases is implemented to demonstrate stability, accuracy, convergence, well-balancedness and mass-conservation of our scheme in the presence of moving wet-dry interfaces.

This manuscript is organized into four further sections. Following this introduction, we briefly introduce the equations and review the DG discretization scheme. We then detail the wetting and drying algorithm in section \ref{sec:WettingDrying}, before demonstrating rigorously the properties of the new limiting approach with numerical examples in section \ref{sec:Results}. We conclude with final remarks and an outlook for future applications.

\section{The shallow water equations and their RKDG discretization} \label{sec:RKDGMethod}

To model two-dimensional waves in shallow water and their interaction with the coast, this study employs the nonlinear shallow water equations. They are derived from the principles of conservation of mass and momentum and can be written compactly in conservative form
\begin{equation} \label{eq:SWESystem}
  \fatvec{U}_t + \nabla \cdot \fatvec{F}(\fatvec{U}) = \fatvec{S}(\fatvec{U}),
\end{equation}
where the vector of unknowns is given by $\fatvec{U} = (h, h\vu)\tran$. Here and below, we have written the partial derivative with respect to time $t$ as $\fatvec{U}_t \equiv \pptext{\fatvec{U}}{t}$ and the divergence with respect to the spatial horizontal coordinates $\Cx=(x,y)\tran$ as $\nabla \cdot\fatvec{F}$, which is applied to each component of $\fatvec{F}$. The quantity $h=h(\Cx,t)$ denotes the fluid depth of a uniform density water layer and $\vu=\vu(\Cx,t)= (u(\Cx,t),v(\Cx,t))\tran$ is the depth-averaged horizontal velocity. The flux function is defined by $\fatvec{F}(\fatvec{U}) = \begin{pmatrix} h\vu \\ h \vu \otimes \vu + \frac{g}{2} h^2\fatvec{I_{2}} \end{pmatrix} $, where $g$ is the gravitational acceleration and $\fatvec{I_{2}}$ the $2\times 2$ identity matrix. Furthermore, the bathymetry (bottom topography) $b=b(\Cx)$ is represented by the source term $\fatvec{S}(\fatvec{U}) = (0, -gh \nabla b)\tran$.

Note that realistic simulations might require further source terms such as bottom friction or Coriolis forcing, since these can significantly influence the wetting and drying as well as the propagation behavior. However, this paper focuses on novel slope limiting techniques for robust inundation modeling. Hence, we only consider the influence of bathymetry and leave the inclusion of further source terms for future investigation to avoid additional (stabilizing) diffusive effects caused, for example, by bottom friction.

For the discretization using the DG method, the governing equations are solved on a polygonal domain $\Omega\subset \mathbb R^2$, which is divided into conforming elements (triangles) $C_i$. On each element, system \eqref{eq:SWESystem} is multiplied by a test function $\varphi$ and integrated. Integration by parts of the flux term leads to the weak DG formulation
\begin{equation} \label{eq:WeakDGForm}
  \int_{C_i} \fatvec{U}_t \varphi \dint{\Cx} -
  \int_{C_i} \fatvec{F}(\fatvec{U})\cdot \nabla\varphi \dint{\Cx} +
  \int_{\partial C_i} \fatvec{F}^*(\fatvec{U}) \cdot \fatvec{n} \ \varphi \dint{\sigma} =
  \int_{C_i} \fatvec{S}(\fatvec{U}) \varphi \dint{\Cx} \ ,
\end{equation}
where $\fatvec{n}$ is the unit outward pointing normal at the edges of element $C_i$. The interface flux $\fatvec{F}^*$ is not defined in general, as the solution can have different values at the interface in the adjacent elements. This problem is circumvented in the discretization by using an (approximate) solution of the corresponding Riemann problem. For the simulations in this study we used the Rusanov solver \citep{Toro2009}.
Another integration by parts of the volume integral over the flux yields the so-called strong DG formulation \citep{Hesthaven2008}
\begin{equation} \label{eq:StrongDGForm}
  \int_{C_i} \fatvec{U}_t \varphi \dint{\Cx} +
  \int_{C_i} \nabla \cdot \fatvec{F}(\fatvec{U})\ \varphi \dint{\Cx} +
  \int_{\partial C_i} \left(\fatvec{F}^*(\fatvec{U}) - \fatvec{F}(\fatvec{U})\right) \cdot \fatvec{n} \ \varphi \dint{\sigma} =
  \int_{C_i} \fatvec{S}(\fatvec{U}) \varphi \dint{\Cx} \ ,
\end{equation}
which recovers the original differential equations within a cell, but with an additional term accounting for the jumps at the interfaces. We will deal with both formulations \eqref{eq:WeakDGForm} and \eqref{eq:StrongDGForm} in this work.

The system of equations is further discretized using a semi-discretization in space with a piecewise polynomial ansatz for the discrete solution components and test functions $\varphi_k$. We obtain formally second-order accuracy by using piecewise linear functions, which are represented by nodal Lagrange basis functions \cite{Hesthaven2008,Giraldo2008}, based on the cell vertices as nodes. The solution in one element is then given by $\fatvec{U}(\Cx,t) = \sum_j \fatvec{\tilde{U}}_j(t) \varphi_j(\Cx)$, where $(\fatvec{\tilde{U}}_j(t))_j$ is the vector of degrees of freedom. The integrals are computed using 3-point and 2-point Gau\ss-Legendre quadrature for volume and line integrals, respectively. Note, that the divergence of the the flux is computed analytically at each quadrature point based on piecewise linear distributions of $h$ and $h\vu$. For the gravitational part this leads to the identity $\nabla \cdot ( \frac{g}{2} h^2 \fatvec{I_{2}} ) = g h\nabla h$ within each cell. This discretization in space leads to a system of ordinary differential equations (ODEs)
\begin{align*}
 \frac{\partial \fatvec{\tilde{U}}_\Delta}{\partial t} = \fatvec{H}_\Delta(\fatvec{\tilde{U}}_\Delta),
\end{align*}
where $\fatvec{\tilde{U}}_\Delta$ contains the degrees of freedom for all cells. The right-hand-side $\fatvec{H}_\Delta(\fatvec{\tilde{U}}_\Delta)$ represents the discretized flux and source term. This system is then solved using a total-variation diminishing (TVD) $s$-stage Runge-Kutta scheme \cite{Shu1988,Gottlieb2001}. In each Runge-Kutta stage a limiter is applied, which is usually employed to stabilize the scheme in case of discontinuities. However, as stated above, it can also be used for dealing with the problem of wetting and drying. In this study, we employ Heun's method, which is the second-order representative of a standard Runge-Kutta TVD scheme.

For the discretization of the bottom topography we use a piecewise linear representation, which is continuous across the interfaces. It is derived from the given data and fixed throughout a simulation. In order to achieve well-balancedness in the DG formulations \eqref{eq:WeakDGForm} and \eqref{eq:StrongDGForm}, exact quadrature of the terms involving $\tfrac{g}{2} h^2$ and the source term is a basic requirement \cite{Xing2006}. This ensures that the sum of cell integrals over flux and source terms together with the line integrals representing interface fluxes vanishes in the lake at rest case. It is achieved by the given quadrature rules. Note that one could also employ a discretization of the bottom topography, which has jumps at the cell interfaces (e.g., resulting from a $L^2$ projection of the exact data). In this case one has to include higher order correction terms into the source term discretization \citep{Xing2006,Xing2010}, which is based on hydrostatic reconstruction of the interface values \citep{Audusse2004}. However, throughout this work we use discretely continuous representations of the bathymetry and exact quadrature rules.

\section{Wetting and drying algorithm} \label{sec:WettingDrying}

After having introduced the governing equations and the general DG discretization, we describe our approach for dealing with wetting and drying, which is a direct generalization of the one-dimensional limiter described in \citet{Vater2015} to the two-dimensional case of triangular meshes. It consists of a flux modification applied in cells with dry nodes and a specially designed limiter, which is non-destructive for the steady state at rest, ensures positivity of the fluid depth and leads to a stable velocity computation. The flux modification is needed to maintain the lake at rest. While the positivity preservation is mostly taken from previous works \citep{Bunya2009,Xing2013}, we introduce a new strategy for the momentum, where we essentially limit the velocities but keep piecewise linear momentum distributions.

\subsection{The wet/dry interface}

\begin{figure}
  \centering
  \includegraphics[width=0.5\textwidth]{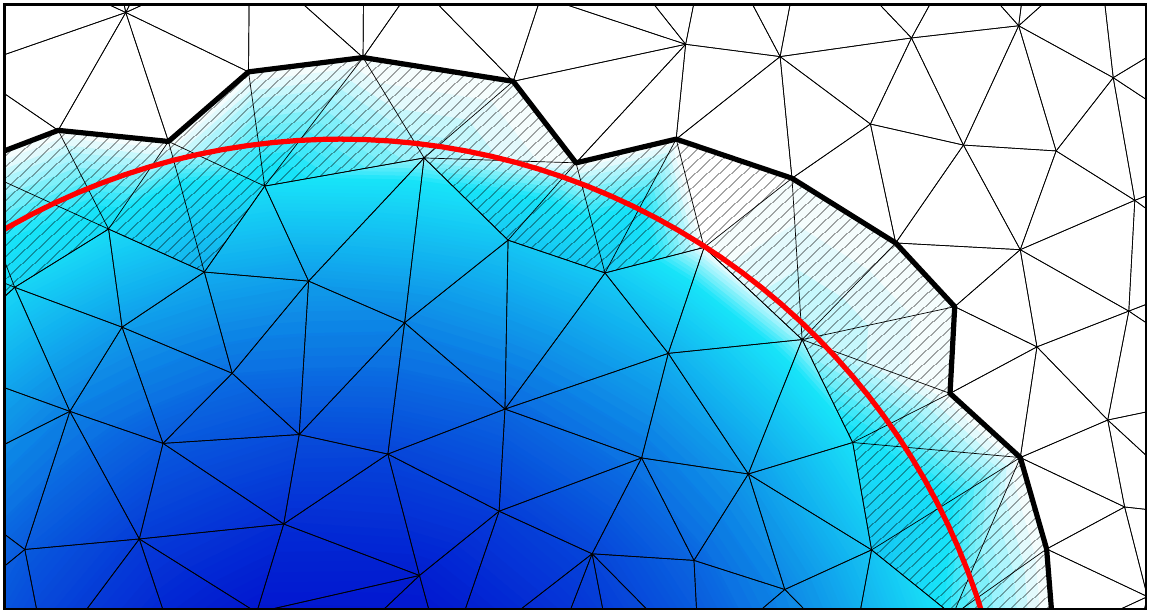}
  \caption{Discontinuous Galerkin discretization of a partly dry domain using a triangular grid. The exact shoreline is printed in red, while the discrete is in black. Semi-dry cells where at least one node has zero fluid depth are hatched.\label{fig:WetDryInterface}}
\end{figure}

\begin{figure}
  \includegraphics[width=0.3\textwidth]{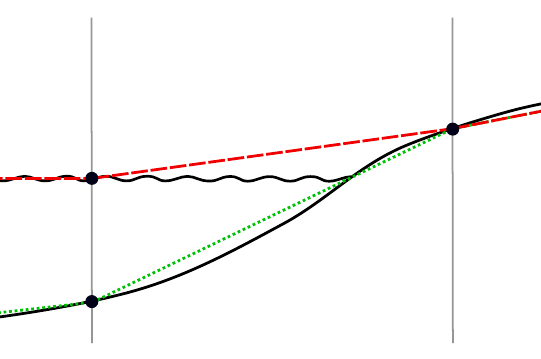}\hfill%
  \includegraphics[width=0.3\textwidth]{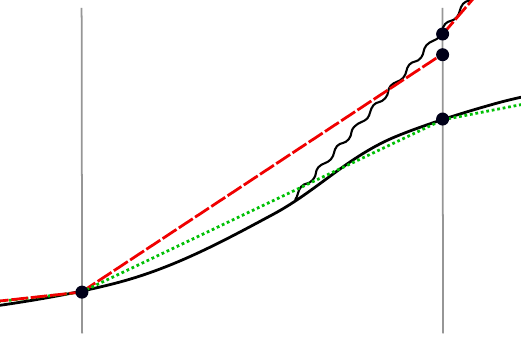}\hfill%
  \includegraphics[width=0.3\textwidth]{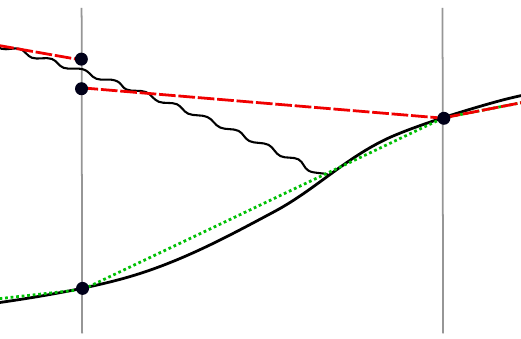}
  \caption{Different configurations of semi-dry cells using a DG discretization with piecewise linear elements in 1D (red dashed: surface elevation, green dotted: bottom topography). Black circles denote nodal values of fluid depth and bathymetry. Displayed are a configuration where \eqref{eq:wetdrycrit} is fulfilled and the gradient in surface elevation should be neglected (left) and two situations which are discretized physically correct (middle and right).\label{fig:SemidryCellTypes}}
\end{figure}

By using a piecewise polynomial DG discretization and enforcing positivity, the discrete shoreline is located at cell interfaces. This can introduce artificial gradients of the fluid depth for the lake at rest in cells which have at least one node with zero fluid depth, the latter we call ``semi-dry'' cells (cf.\ figure \ref{fig:WetDryInterface}). Here, we define the discrete lake at rest by interpolating the exact surface elevation $H=h+b$ at triangle nodes and setting the momentum to zero. This results in a continuous representation of the fluid depth. On the other hand, there may be semi-dry cells that are approximated physically correct (e.g., in a dam-break situation where the water comes from higher elevation) and which must be distinguished from the lake-at-rest case (figure \ref{fig:SemidryCellTypes}). We do this by comparing the maximum total height with the maximum bottom topography within a cell. If the maximum total height is not larger than the maximum bottom topography within cell $C_i$, i.e.,
\begin{equation} \label{eq:wetdrycrit}
  \max_{\Cx \in C_i} H(\Cx) - \max_{\Cx \in C_i} b(\Cx) < \tolwet \ ,
\end{equation}
we are possibly in a local lake-at-rest situation, and the cell must be specially treated. Here, we have introduced the parameter $\tolwet$, which is a threshold for the fluid depth under which a point is considered dry. At such points the velocity is set to zero, which is computed by division of $(hu)/h$ elsewhere. Otherwise, if \eqref{eq:wetdrycrit} is not fulfilled, the cell is treated as a completely wet cell.

To render the method well-balanced for the discrete lake at rest, we neglect the volume integrals over the terms involving $\frac{g}{2} h^2\fatvec{I_{2}}$ and the source term in semi-dry cells where \eqref{eq:wetdrycrit} is fulfilled. This is equivalent to setting $g$ to zero. For the strong DG formulation \eqref{eq:StrongDGForm} no further modifications are needed. In the case of the lake at rest where no advection is present, all volume integrals vanish for wet and semi-dry cells, and the flux difference $\fatvec{F}^*(\fatvec{U}) - \fatvec{F}(\fatvec{U})$, which is computed at the interfaces, automatically cancels due to continuity of the surface elevation and consistency of the numerical flux function. This is different for the weak DG formulation \eqref{eq:WeakDGForm}, where also the volume integrals associated with gravitational forces are neglected. The interface contributions involving $\frac{g}{2} h^2$ cannot be simply neglected, since the numerical flux at the wet interfaces couples to adjacent cells, which is needed in case of perturbations. For these wet interfaces an additional flux term is introduced which balances the numerical flux. It includes only the gravitational part and is based on the fluid depth of the semi-dry cell at the wet interface. The momentum equation in a semi-dry cell then reads
\begin{gather*}
  \int_{C_i} \varphi (h\vu)_t \dint{\Cx} -
  \int_{C_i} \nabla \varphi \cdot \Bigl( h \vu \otimes \vu + \cancel{\frac{g}{2} h^2\fatvec{I_{2}}} \Bigr) \dint{\Cx}\, + \\
  \sum_{j=1}^3 \int_{E_j^i} \Bigl( \fatvec{F}_{h\vu}^* - \frac{g}{2} (h^{-})^2 \fatvec{I_{2}} \Bigr) \cdot \fatvec{n} \ \varphi \dint{\sigma} =  - \int_{C_i} \cancel{\varphi g h \nabla b} \dint{\Cx} \ ,
\end{gather*}
where $h^{-}$ is the value of the fluid depth at the interface based on cell $C_i$, and $E_j^i$, $j\in\{1,2,3\}$ are the three edges of $C_i$. In the equation, we have canceled the above mentioned volume integrals associated with gravitational forces. If this discretization is applied to the lake at rest with $\vu \equiv \fatvec{0}$, then all the edge contributions in a semi-dry cell vanish. For dry edges this is because $h$ is zero, and for wet cells the difference computed under the integral cancels to zero. Note, however, that well-balancing is easier accomplished by using the strong form, and this form also leads to slightly better results as we will see in section \ref{sec:Results}. These described modifications can be also interpreted as a flux limiting approach.

\subsection{Limiting of the fluid depth} \label{sec:LimitDepth}

Limiting with respect to fluid depth is based on a Barth/Jespersen-type limiter \citep{Barth1989}, which fulfills the requirement to not alter a well-balanced discrete solution, when limiting in total height $H=h+b$. Positivity is attained by ensuring positivity in the mean and redistribution of fluid depth within each cell.

Barth/Jespersen-type limiters modify the solution within a cell, such that it does not exceed the maximum and minimum of the cell mean values of adjacent cells. In this work we study the original version by Barth and Jespersen, which incorporates the cells which are connected by a common edge to the cell under consideration. Additionally, we consider a generalization for triangular grids, which was introduced by \citet{kuzmin_vertex-based_2010}. This limiter incorporates the cells which are connected by a common vertex for comparison. We refer to these two versions as ``edge-based'' and ``vertex-based'' limiter, respectively. Note, that the main goal is not to compare these two limiters. We introduce the two versions to offer some flexibility in the computational setup, since some algorithms require edge based computations for efficiency, whereas others are organized by nodal representations.

Given the cell average or centroid value $H_c = \overline{H} = \overline{h+b}$ of the total height, the piecewise linear in-cell distribution can be described by $H(\Cx) = H_c + (\nabla H)_c \cdot (\Cx-\Cx_c)$. A Barth/Jespersen-type limiter modifies this to
\begin{equation*}
  \hat{H}(\Cx) = H_c + \alpha_e (\nabla H)_c \cdot (\Cx-\Cx_c), \quad 0 \le \alpha_e \le 1 ,
\end{equation*}
where $\alpha_e$ is chosen, such that the reconstructed solution is bounded by the maximum and minimum centroid values of a given cell neighborhood:
\begin{equation*}
  H_c^\mathrm{min} \le \hat{H}(\Cx) \le H_c^\mathrm{max} .
\end{equation*}
For the original Barth/Jespersen limiter this cell neighborhood is given by the considered cell and the three cells sharing an edge with this cell. In case of the limiter described by \citet{kuzmin_vertex-based_2010}, the cell neighborhood is given by the considered cell and the surrounding cells sharing a vertex with this cell. The correction factor is explicitly defined as
\begin{equation*}
  \alpha_e = \min_i
  \begin{cases}
    \min \left\{ 1, \frac{H_c^\mathrm{max} - H_c}{H_i - H_c}\right\}, &
      \text{if } H_i - H_c > 0 \\
    1, & \text{if } H_i - H_c = 0 \\
    \min \left\{ 1, \frac{H_c^\mathrm{min} - H_c}{H_i - H_c}\right\}, &
      \text{if } H_i - H_c < 0
  \end{cases}
\end{equation*}
where $H_i$ are the in-cell values of $H$ at the three vertices of the triangle. The limited fluid depth $\hat{h}$ is then recovered by $\hat{h} = \hat{H} - b$.

Positivity is enforced in a second step by the positive depth operator originally proposed by \citet{Bunya2009}. Note that this approach is closely related to the positivity preserving limiter introduced in \citet{Xing2013}, the latter being also suitable for higher order elements. Here we also follow \citet{Xing2013} by relying on a CFL time step restriction to preserve positivity in the mean. Applied to the RKDG2 method, the CFL limit is 1/3, which is less restrictive than the time step restriction to ensure linear stability. Let us express the piecewise linear function $\hat{h}$ by its nodal representation with Lagrange basis functions
\begin{equation*}
  \hat{h}(x,y) = \sum_{i=1}^3 \hat{h}_i \, \varphi_i(x,y) \quad \text{for } (x,y) \in C_k ,
\end{equation*}
where we take as nodes the vertices of the triangle $C_k$. Then positivity on the whole triangle is obtained by requiring positivity for the nodal values. Denoting the final limited values by $h^{\lim}$ and $h^{\lim}_i$ with $h^{\lim} = \sum h^{\lim}_i \, \varphi_i$, $h^{\lim}_i$ is determined by the following procedure: If $\hat{h}_i  \ge 0 \  \forall i \in \{ 1, 2, 3\}$, then
\begin{equation*}
  h^{\lim}_i = \hat{h}_i, \ i \in \{ 1, 2, 3\} .
\end{equation*}
Otherwise we determine the order of nodal indices $n_i \in \{ 1, 2, 3\}$ that satisfy $\hat{h}_{n_1} \le \hat{h}_{n_2} \le \hat{h}_{n_3}$ and compute the values in the following sequence:
\begin{align*}
  h^{\lim}_{n_1} &= 0\\
  h^{\lim}_{n_2} &= \max\{ 0, \hat{h}_{n_2} - (h^{\lim}_{n_1}-\hat{h}_{n_1}) / 2 \}\\
  h^{\lim}_{n_3} &= \hat{h}_{n_3} - (h^{\lim}_{n_1}-\hat{h}_{n_1}) - (h^{\lim}_{n_2}-\hat{h}_{n_2})
\end{align*}
As \citet{Bunya2009} show, this algorithm lowers the depths at nodes $n_2$ and $n_3$ by equal amounts, and the algorithm is mass conserving.

\subsection{Velocity-based limiting of the momentum}


In a last step the momentum distribution is modified by limiting the in-cell variation of the resulting velocity distribution. This provides a stable computation near the wet/dry interface in situations when both, $h$ and $(h\vu)$ are small. It is designed to keep a piecewise linear momentum distribution with fixed cell mean values. As noted in the introduction, this approach originates from FV methods, where the reconstruction of interface values by means of other than the primary variables has a long tradition \cite{Leer2006,Audusse2004}. For FV methods this does not pose a problem, since the reconstructed values are only used for flux computation. In DG methods the situation is more complicated, since the whole in-cell solution is limited and used throughout the computations. Therefore, one is usually bound to use the primary variables for limiting.

For the momentum limiting we first compute preliminary limited velocity components $\hat{u}_i$ (and similarly $\hat{v}_i$) at each node $i$ of the triangle
\begin{equation*}
  \hat{u}_i = \max\{\min\{ u_i, u_c^\mathrm{max}\}, u_c^\mathrm{min}\}
\end{equation*}
where $u_i = (hu)_i/h_i$ and $u_c = (hu)_c/h_c$ are the $x$-velocities computed from the nodal and centroid values of momentum and (the unlimited) fluid depth. Note that in case $h_i<\tolwet$ the velocity is set to 0. The minimum and maximum values $u_c^\mathrm{min/max}$ are determined as in subsection \ref{sec:LimitDepth} for the total height by considering the centroid values of the neighboring cells which share a common edge (edge-based limiter) or a common vertex (vertex-based limiter) with the cell.

This results in three different linear momentum distributions based on two of the three preliminary nodal velocities, by keeping the cell mean value of the momentum and the distribution of the fluid depth. For the three possibilities we obtain a velocity for the third node with
\begin{gather*}
  \hat{u}^{23}_1 = \frac{3 (hu)_c - h^{\lim}_2\cdot\hat{u}_2 - h^{\lim}_3\cdot\hat{u}_3}{h^{\lim}_1}, \quad
  \hat{u}^{13}_2 = \frac{3 (hu)_c - h^{\lim}_1\cdot\hat{u}_1 - h^{\lim}_3\cdot\hat{u}_3}{h^{\lim}_2}, \\
  \hat{u}^{12}_3 = \frac{3 (hu)_c - h^{\lim}_1\cdot\hat{u}_1 - h^{\lim}_2\cdot\hat{u}_2}{h^{\lim}_3},
\end{gather*}
where the lower index denotes the node for which the velocity is computed and the upper index defines, which nodal velocities this is based on. The final limited momentum component is then chosen to produce the smallest in-cell velocity variation. Set
\begin{equation}\label{eq:grad_velo}
\begin{aligned}
  \delta \hat{u}_1 &= \max\{ \hat{u}^{23}_1, \hat{u}_2, \hat{u}_3 \} -
                      \min\{ \hat{u}^{23}_1, \hat{u}_2, \hat{u}_3 \} , \\
  \delta \hat{u}_2 &= \max\{ \hat{u}_1, \hat{u}^{13}_2, \hat{u}_3 \} -
                      \min\{ \hat{u}_1, \hat{u}^{13}_2, \hat{u}_3 \} , \\
  \delta \hat{u}_3 &= \max\{ \hat{u}_1, \hat{u}_2, \hat{u}^{12}_3 \} -
                      \min\{ \hat{u}_1, \hat{u}_2, \hat{u}^{12}_3 \} .
\end{aligned}
\end{equation}
If $\delta \hat{u}_1 \le \delta \hat{u}_i$ for $i \in \{2, 3\}$, then
\begin{equation} \label{eq:hulim1}
  (hu)^{\lim}_1 = h^{\lim}_1 \cdot \hat{u}^{23}_1 , \quad
  (hu)^{\lim}_2 = h^{\lim}_2 \cdot \hat{u}_2 ,      \quad
  (hu)^{\lim}_3 = h^{\lim}_3 \cdot \hat{u}_3 .
\end{equation}
If $\delta \hat{u}_2 \le \delta \hat{u}_i$ for $i \in \{1, 3\}$, then
\begin{equation}\label{eq:hulim2}
  (hu)^{\lim}_1 = h^{\lim}_1 \cdot \hat{u}_1 ,      \quad
  (hu)^{\lim}_2 = h^{\lim}_2 \cdot \hat{u}^{13}_2 , \quad
  (hu)^{\lim}_3 = h^{\lim}_3 \cdot \hat{u}_3 .
\end{equation}
Otherwise
\begin{equation}\label{eq:hulim3}
  (hu)^{\lim}_1 = h^{\lim}_1 \cdot \hat{u}_1 , \quad
  (hu)^{\lim}_2 = h^{\lim}_2 \cdot \hat{u}_2 , \quad
  (hu)^{\lim}_3 = h^{\lim}_3 \cdot \hat{u}^{12}_3 .
\end{equation}
The final limited momentum is then given by $(hu)_h^{\lim} = \sum_i [(hu)^{\lim}_i \, \varphi_i(x,y)]$. The same procedure is applied to the $y$-momentum. In conclusion, the wetting and drying algorithm can be summarized as follows:
\clearpage
\begin{boxtext}
    \section*{Limiter-Based Wetting and Drying Treatment}
    \begin{enumerate}
      \item Flux modification
      \begin{enumerate}
        \item Set $g$ to 0 in volume integrals of semi-dry cells, which satisfy \eqref{eq:wetdrycrit}, and add additional interface flux if using the weak DG formulation.
      \end{enumerate}
      \item Limiting of fluid depth
      \begin{enumerate}
        \item\label{alg:limha} Apply edge-based \cite{Barth1989} OR vertex-based \cite{kuzmin_vertex-based_2010} limiter to total height $H=h+b$.
        \item Apply positive depth operator \cite{Bunya2009} to limited $\hat{h}$ obtained from $\hat{H}$ in step \ref{alg:limha}.
      \end{enumerate}
      \item Limiting of momentum
      \begin{enumerate}
        \item\label{alg:limma} Apply edge-based \cite{Barth1989} OR vertex-based \cite{kuzmin_vertex-based_2010} limiter to velocities at triangle nodes.
        \item\label{alg:limmb} Extrapolate in-cell velocity distribution from two out of three nodal values obtained in step \ref{alg:limma}.
        \item Determine discrete in-cell velocity variation from the three distributions obtained in step \ref{alg:limmb}.
        \item Compute limited momentum from velocities with smallest variation and limited fluid depth (cf.\ \eqref{eq:hulim1}--\eqref{eq:hulim3}).
      \end{enumerate}
    \end{enumerate}
\end{boxtext}

\section{Numerical results} \label{sec:Results}

In the following we demonstrate the major functionalities of the limiting procedure described in the last section. Using a hierarchy of test cases, where we start with configurations where the exact solution to the shallow water equations is known, we show the scheme's mass conservation and well-balancedness, as well as the correct representation of the shoreline. Two test cases, which originate from \citet{Thacker1981}, particularly demonstrate the scheme's ability of representing a moving shoreline. Two further test cases are derived from laboratory experiments, which, together with the runup onto a linearly sloping beach, are standard test cases for the evaluation of operational tsunami models \citep{Synolakis2007}.

For the simulations, we use both versions of the limiter, i.e., vertex-based and edge-based limitation of total height and velocity and show that although they differ slightly in computational complexity and added numerical diffusion, they both yield comparable and accurate results.  The presented limiter depends on one free parameter -- the wet/dry tolerance $\tolwet$ -- that defines the fluid depth threshold, below which a point is considered dry. This is especially important for the computation of the velocity. We comment on this tolerance for each test case and, overall, conclude that the stability of the limiting strategy is not sensitive to it. However, it can influence the discrete location of the wet/dry interface.

Apart from the first test case, where we compare the weak and strong DG formulations concerning well-balancedness, we only present results using the strong DG form in the simulations. Although the strong DG form shows somehow better results in case of the lake at rest, the other test cases produce nearly indistinguishable results for both DG formulations.

Throughout this section, we set the acceleration due to gravity to $g=9.80616\, \mathrm{m/s^2}$ and omit the units of measurement of the physical quantities, which should be thought in standard SI units with $\mathrm{m}$ (meters), $\mathrm{s}$ (seconds) etc. For the spatial discretization we mostly use regular grids, which are usually derived from one or more rectangles, each divided into two triangles. Such a grid is then repeatedly uniformly refined by bisection to obtain the desired resolution (see \citet{Behrens2005} for details on grid generation). The discrete initial conditions and the bottom topography are derived from the analytical ones by interpolation at the nodal points (vertices of triangles).

Explicit methods for the solution of hyperbolic problems are usually subject to a CFL time step restriction \citep{Courant1928}, which is of the form $\Delta t \le \mathrm{cfl} \, h_\Delta / c_\mathrm{max}$. Here, $h_\Delta$ defines a grid parameter and $c_\mathrm{max}$ is the maximum propagation speed of information. For one-dimensional problems, \citet{Cockburn1991} proved that $\mathrm{cfl_{1D}} = 1/3$ for the RKDG2 method. However, this cannot be directly transferred to two-dimensional triangular grids. We follow the work of \citet{Kubatko2008}, who propose as grid parameter the radius of the smallest inner circle of the triangles surrounding a vertex. These nodal values are further aggregated to each triangle by taking the minimum of its three vertex values. The 2D CFL number then approximately relates to its 1D counterpart by $\mathrm{cfl_{2D}} \approx 2^{-1/(p-1)} \mathrm{cfl_{1D}}$, where $p$ is the order of discretization. This results in $\mathrm{cfl_{2D}} \approx 0.233$ for our RKDG2 method. Note that this limit is more restrictive than the time step restriction of 1/3 to ensure positivity in the mean. If not stated otherwise, we choose a constant time step size $\Delta t^n = \Delta t$ which guarantees that the CFL condition is essentially satisfied.

Besides fluid depth and momentum we often show the velocity $\vu = (h\vu) / h$ with its in-cell distribution, which is derived diagnostically by the quotient of the two other quantities.

\subsection{Lake at rest}

As a first test we conduct two simple ``lake at rest'' experiments with different bathymetries to examine the well-balancedness of our scheme. In a quadratic and periodic domain $\Omega=[0,1]^2$ the first bathymetry is defined by $b(\Cx) = \max\left\{ 0, 0.25 - 5 ((x-0.5)^2 + (y-0.5)^2)\right\}$, which features a local, not fully submerged parabolic mountain centered around the mid point $(0.5,0.5)^{\top}$. The initial fluid depth and velocity are given by
\begin{equation}\label{eq:iniLakeAtRest}
\begin{aligned}
  h(\Cx, 0) &= \max\left\{0, 0.1 - b(\Cx) \right\}, \\
 \vu(\Cx,0) &= \fatvec{0}.
\end{aligned}
\end{equation}
This is a steady state solution which should be reproduced by the numerical scheme. We run simulations using the strong and weak DG formulations with both limiters and a grid resolution of $\Delta x \approx 0.022$ (leg of right angled triangle) until $\Tend = 40$. A time step of $\Delta t = 0.002$ is used, which results in 20\,000 time steps. The wet/dry tolerance is chosen as $\tolwet=10^{-6}$.

\begin{figure}
  \centering
  \includegraphics[width=0.5\textwidth]{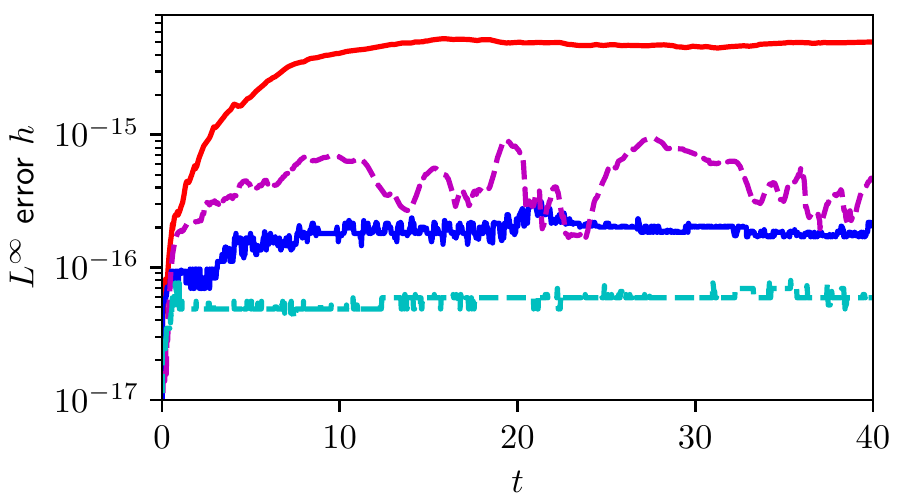}%
  \includegraphics[width=0.5\textwidth]{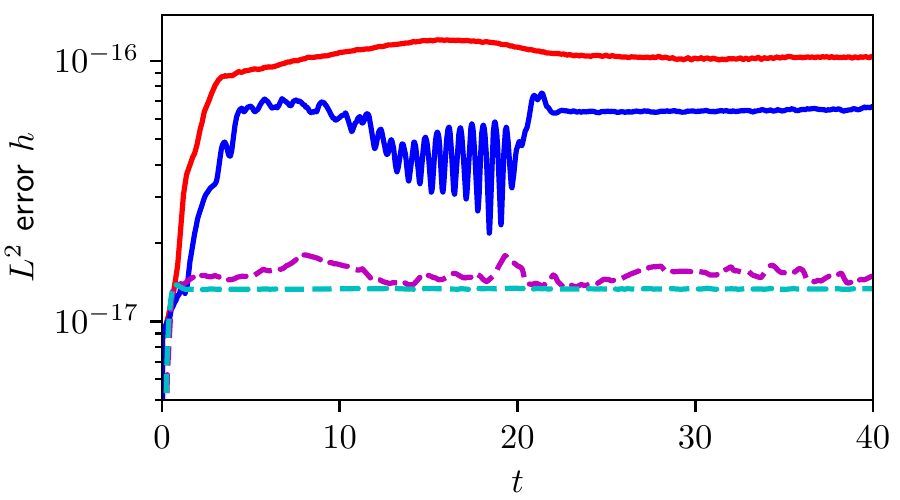}\\
  \includegraphics[width=0.5\textwidth]{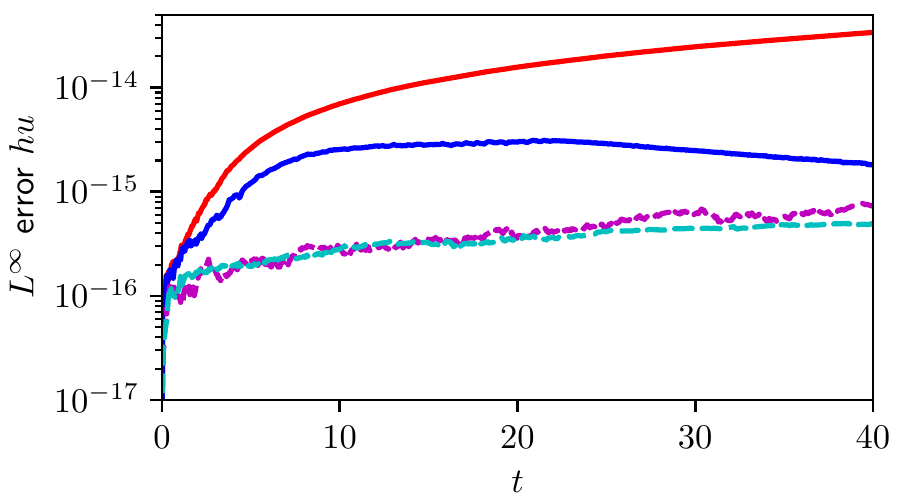}%
  \includegraphics[width=0.5\textwidth]{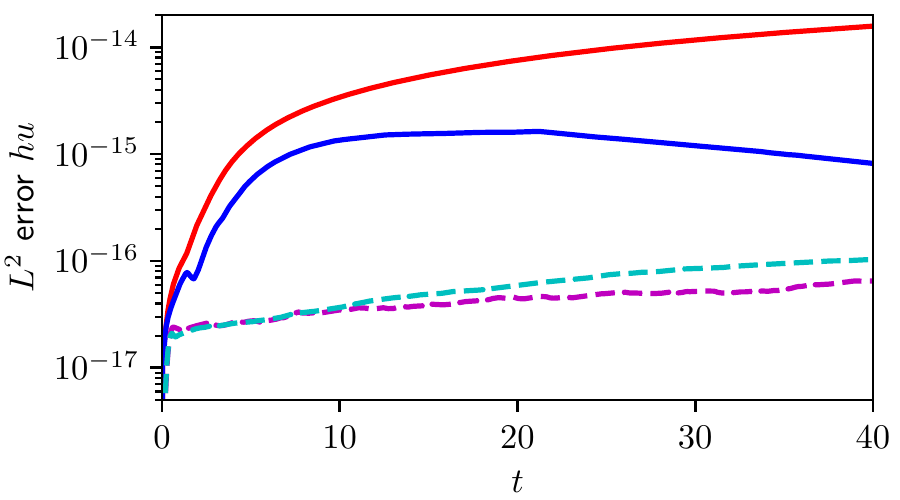}
  \caption{Lake at rest: errors over time for the  for fluid depth (top) and momentum (bottom) in the $L^\infty$ (left) and $L^2$ norms (right) for the first bathymetry setup. Vertex-based limiter with weak DG formulation (red), edge-based limiter with weak DG formulation (blue), vertex-based limiter with strong DG formulation (magenta dashed), edge-based limiter with strong DG formulation (cyan dashed).\label{fig:wellbalancing}}
\end{figure}

The results are depicted in figure \ref{fig:wellbalancing}. We show the error in the $L^2$ as well as the maximum ($L^\infty$) norm for the fluid depth (top row) and momentum (bottom row) for all four possible configurations. All combinations are well-balanced almost up to machine precision considering fluid depth. The momentum is also balanced, except for the vertex-based limiter in combination with the weak DG form, which shows a slowly growing -- yet well controlled -- momentum error. The simulations using the strong DG form show generally smaller errors.This behavior can be explained by the superior balancing property of the strong DG formulation, which has been already observed in \citet{Beisiegel2014}.

\begin{figure}
  \centering
  \includegraphics[width=0.5\textwidth]{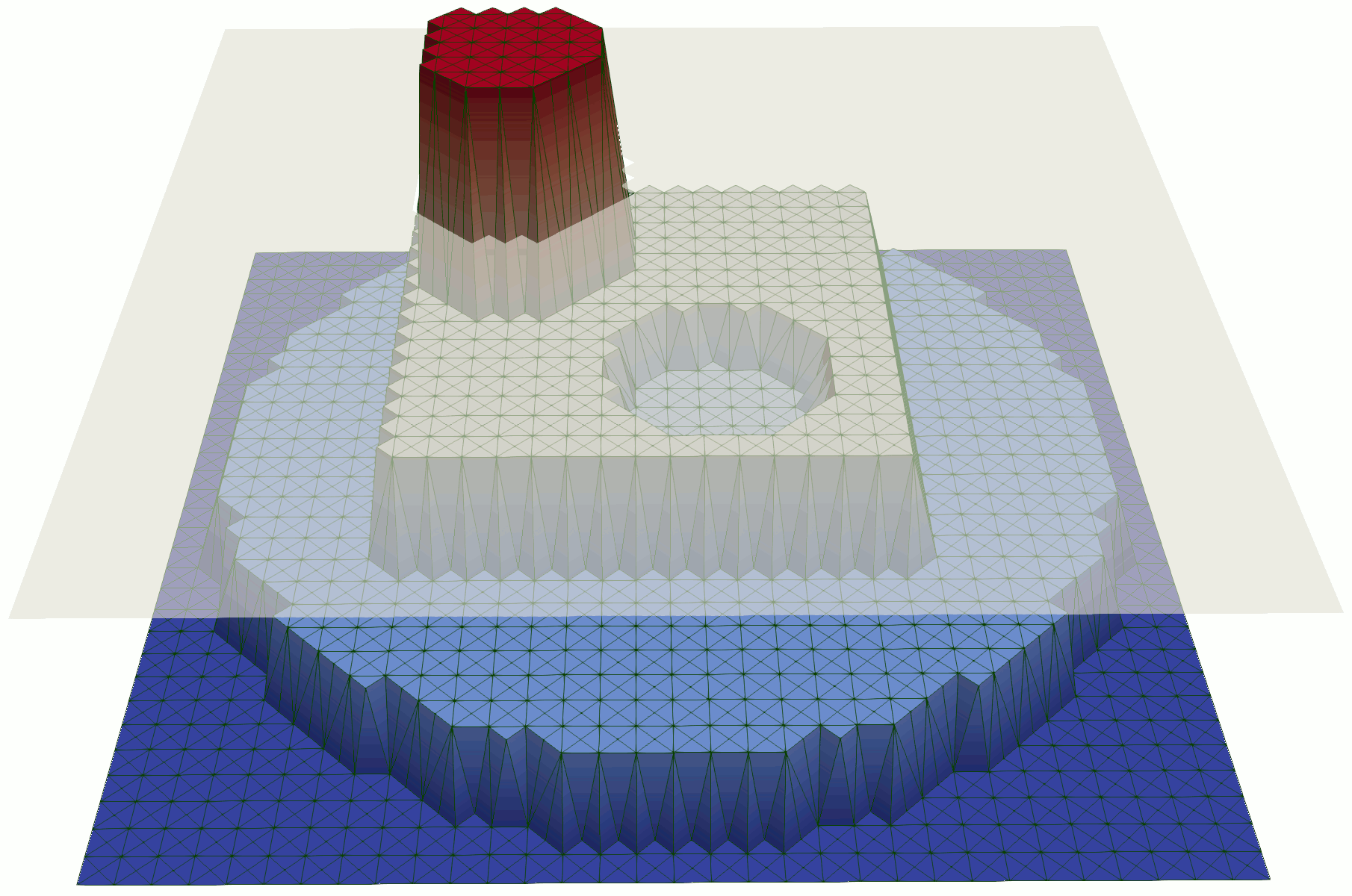}%
  \includegraphics[width=0.5\textwidth]{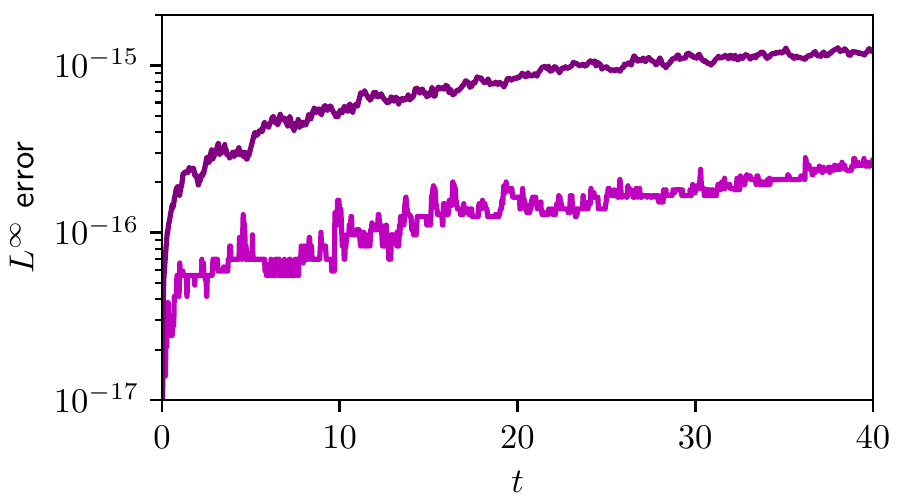}%
  \caption{Lake at rest: Initial setup for the second bathymetry configuration (left), and $L^\infty$ errors over time for fluid depth (magenta) and momentum (dark purple). Vertex-based limiter with strong DG formulation.\label{fig:wellbalstep}}
\end{figure}

For the second bathymetry setup, we define the sub-domains
\begin{align*}
  \Omega_1 &= \left\{ \Cx \in \Omega \big| \|\Cx - (0.35,0.65)^{\top}\| < 0.1 \right\},\\
  \Omega_2 &= \left\{ \Cx \in \Omega \big| \|\Cx - (0.55,0.45)^{\top}\| < 0.1 \right\},\\
  \Omega_3 &= \left\{ \Cx \in \Omega \big| | x - 0.47 | < 0.25 \wedge | y - 0.55 | < 0.25 \right\} \text{ and}\\
  \Omega_4 &= \left\{ \Cx \in \Omega \big| \|\Cx - (0.5,0.5)^{\top}\|   < 0.45 \right\}.
\end{align*}
The bathymetry is given by
\begin{equation*}
  b (\Cx) =
  \begin{cases}
    0.15 & \text{if } \Cx \in \Omega_1 \\
    0.05 & \text{if } \Cx \in \Omega_2 \\
    0.07 & \text{if } \Cx \in \Omega_3 \setminus \{\Omega_1 \cup \Omega_2\} \\
    0.03 & \text{if } \Cx \in \Omega_4 \setminus \{\Omega_3\} \\
    0 & \text{otherwise.}
  \end{cases}
\end{equation*}
The initial conditions are given as in \eqref{eq:iniLakeAtRest} (cf.\ figure \ref{fig:wellbalstep} left). Although the analytical setup has discontinuities in the bathymetry, we also interpolate bathymetry and initial condition by piecewise linear and continuous approximations at the cell vertices. Using the same grid and time step size as in the first setup, we obtain $L^\infty$ errors for fluid depth and momentum over time as displayed in figure \ref{fig:wellbalstep}, right. One can see that the scheme is also well-balanced in this case with steps in the bathymetry.

\subsection{Tsunami runup onto a linearly sloping beach}

\begin{figure}
  \centering
  \includegraphics[width=0.5\textwidth]{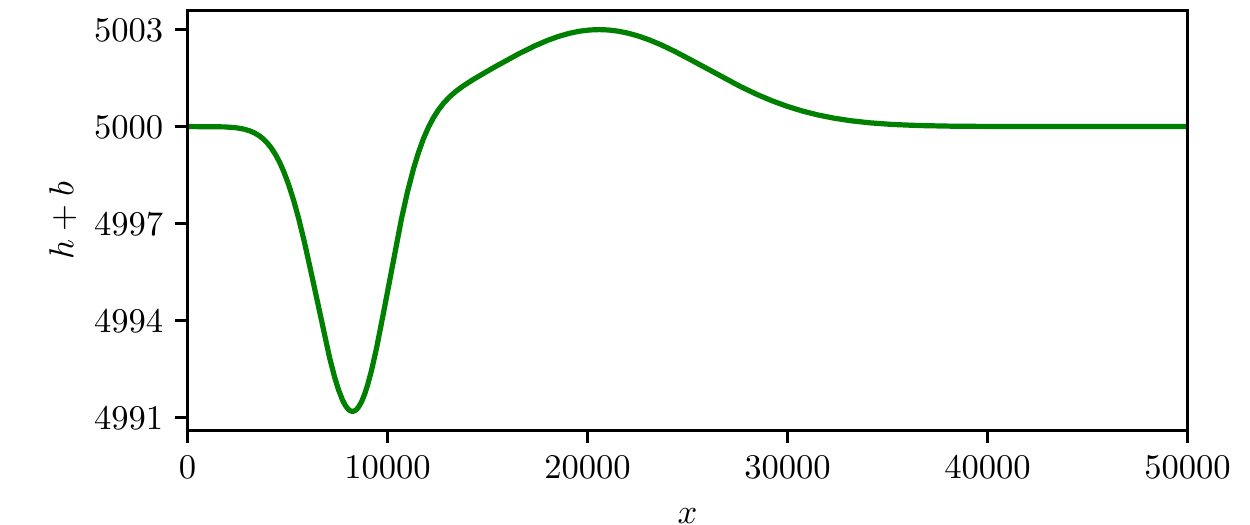}
  \caption{Tsunami runup onto a beach: initial surface elevation at $t=0$.\label{fig:TsunmaiRunupIni}}
\end{figure}

A standard benchmark problem to evaluate wetting and drying behavior of a numerical scheme is the wave runup onto a plane beach. We perform this quasi one-dimensional test case \citep{LongWaveRunupModels1_2004} to compare the results to the ones already obtained with the one-dimensional version of the scheme in \citet{Vater2015}. The test case admits an exact solution following a technique developed in \citet{Carrier2003}.
In a rectangular domain $\Omega = [-400, 50\,000]\times[0,400]$ with linearly sloping bottom topography $b(\Cx) = 5000 - \alpha x$, $\alpha = 0.1$, and initial velocity $\vu(\Cx,0)=\fatvec{0}$, an initial surface elevation is prescribed in non-dimensional variables by
\begin{equation*}
  \eta'(x') = a_1 \exp\left\{ k_1 (x'-x_1)^2 \right\} -
              a_2 \exp\left\{ k_2 (x'-x_2)^2 \right\} .
\end{equation*}
The parameters are given by $a_1 = 0.006$, $a_2 = 0.018$, $k_1 = 0.4444$, $k_2 = 4$, $x_1 = 4.1209$ and $x_2 = 1.6384$. Then, the initial surface profile is recovered by taking $x=Lx'$ and $\eta = \alpha L \eta'$ with the reference length $L=5000$ (cf.\ figure~\ref{fig:TsunmaiRunupIni}). As the solution near the left boundary of $\Omega$ is always dry and on the right boundary outgoing waves should not be reflected, we set wall and transparent boundary conditions for the left and right boundary, respectively. The boundaries in $y$-direction are set periodic, to avoid any artifacts coming from the definition of the boundary conditions.

\begin{figure}
  \centering
  \includegraphics[width=0.33\textwidth]{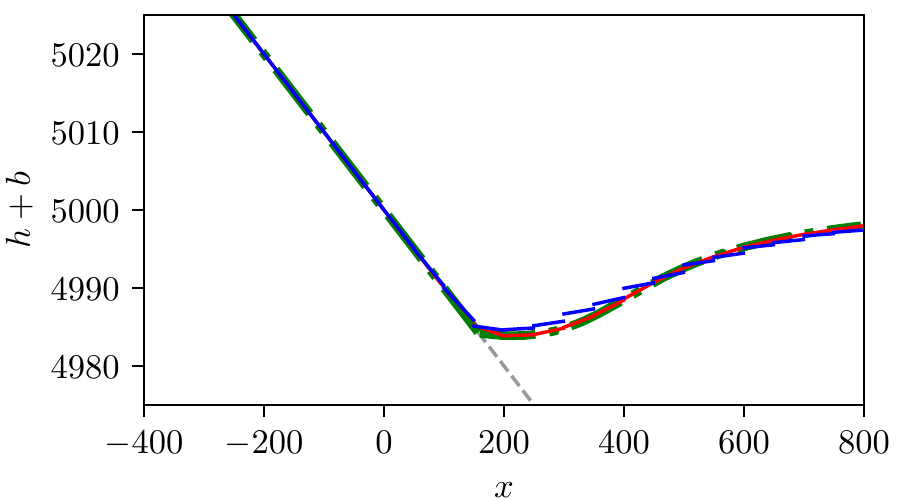}%
  \includegraphics[width=0.33\textwidth]{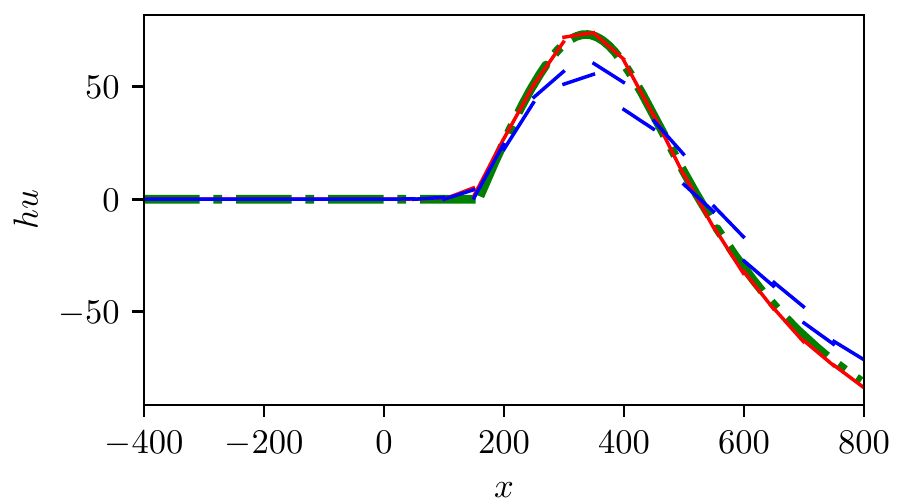}%
  \includegraphics[width=0.33\textwidth]{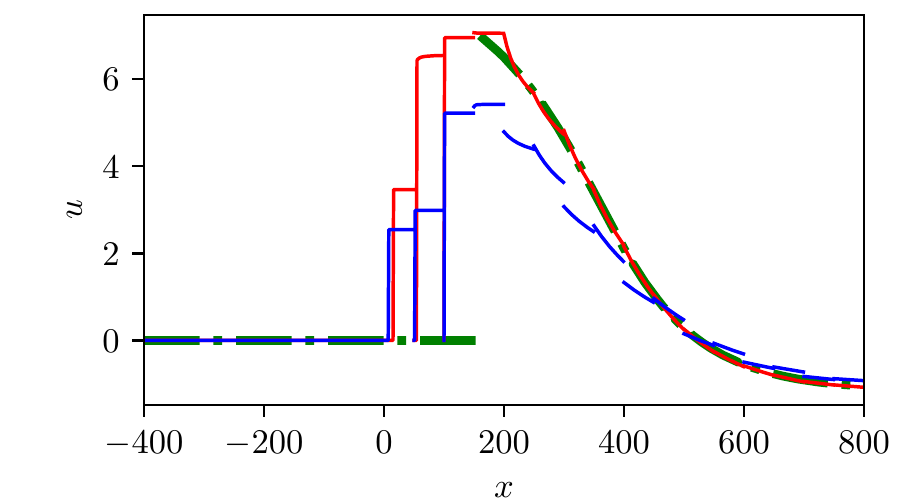}\\
  \includegraphics[width=0.33\textwidth]{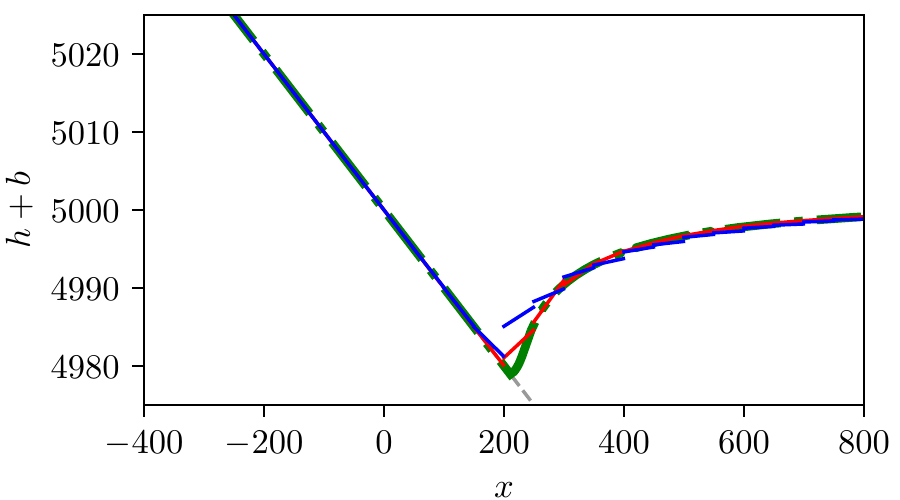}%
  \includegraphics[width=0.33\textwidth]{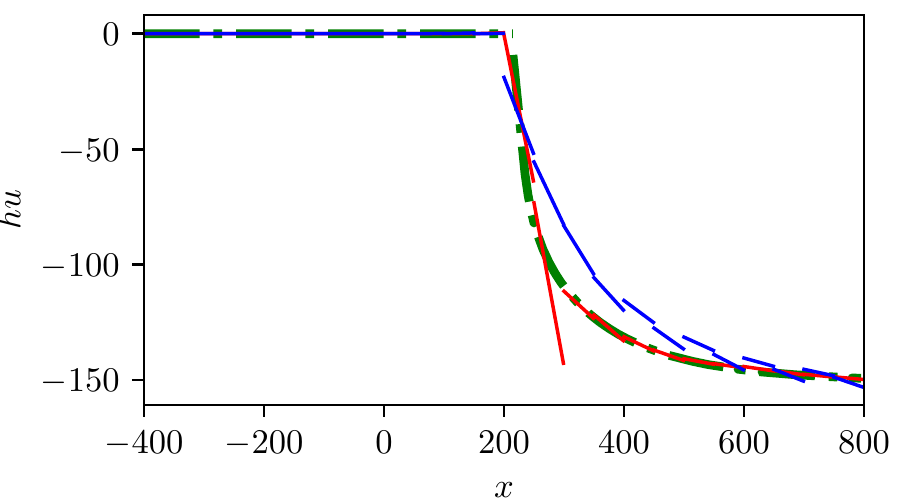}%
  \includegraphics[width=0.33\textwidth]{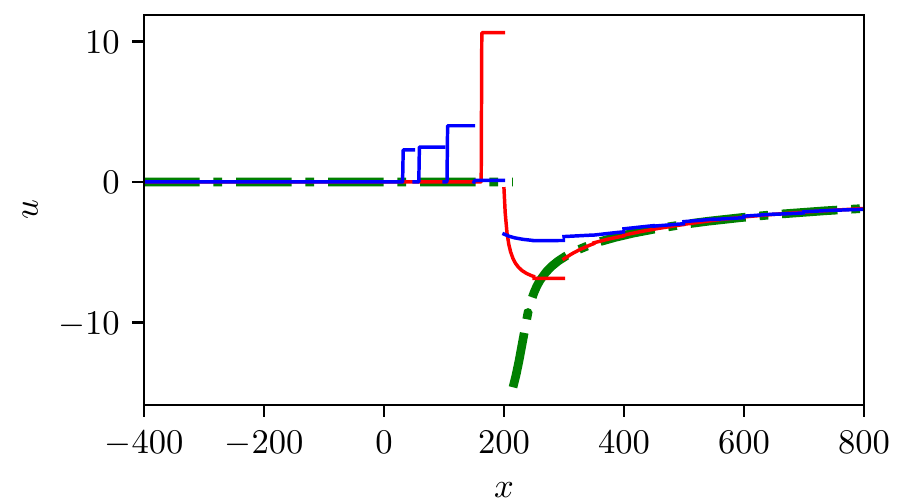}\\
  \includegraphics[width=0.33\textwidth]{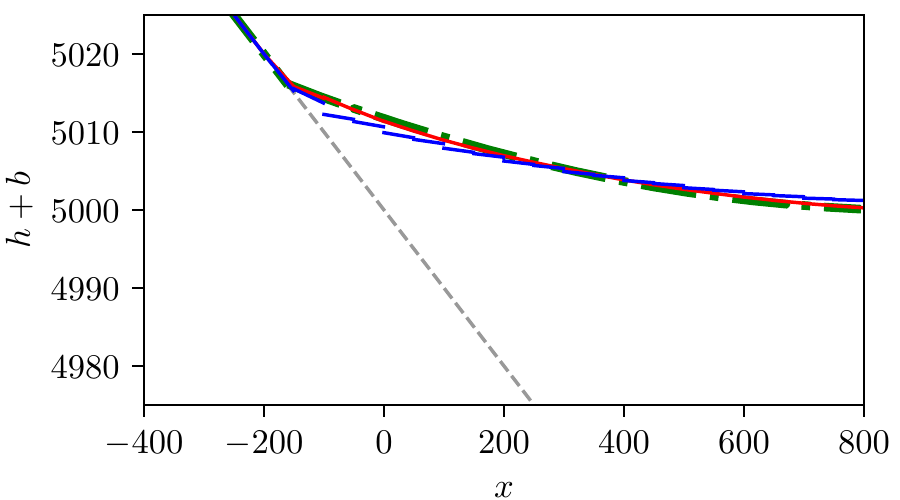}%
  \includegraphics[width=0.33\textwidth]{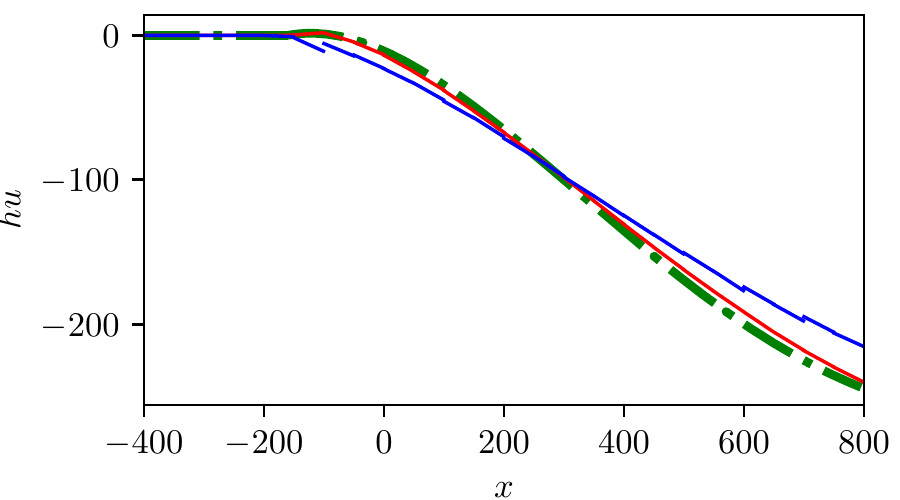}%
  \includegraphics[width=0.33\textwidth]{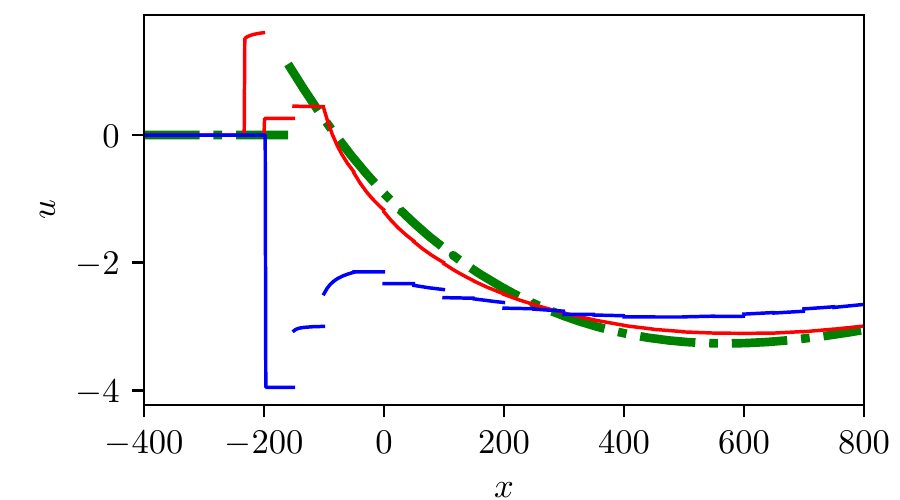}
  \caption{Tsunami runup onto a beach: surface elevation, $x$-momentum and $x$-velocity (derived by $u=(hu)/h$) along line $y=200$ at times $t=160$ (top), $t=175$ (middle) and $t=220$ (bottom). Exact solution (green dash-dotted), vertex-based limiter (red), edge-based limiter (blue).\label{fig:runup}}
\end{figure}

The simulations are run with a time step of $\Delta t = 0.04$, and a spatial resolution of $\Delta x = 50$ which corresponds to the length of the leg of a right angled triangle. The resulting Courant number is approximately $0.25$, which is attained offshore where the fluid depth is largest. The wet/dry tolerance is set to $\tolwet=10^{-2}$. We compare our numerical results with the analytical solution on the interval $x\in[-400, 800]$ at times $t = 160$, $175$ and $220$. The results are depicted in figure \ref{fig:runup}.

The left column shows the free surface elevation, where the simulations with both limiters (red and blue lines) match the exact solution (green line). However, it can be observed that due to its smaller stencil the edge-based limiter develops spurious discontinuities at cell edges. The middle and right column show the momentum and the reconstructed velocity at the respective times. As expected, the momentum is reproduced well while the velocity shows some spurious over-, and undershoots in the near-dry area, but good results elsewhere.
In general, the vertex-based limiter yields qualitatively better results for this test case. The results are comparable to those in one space-dimension given in \citet{Vater2015} and demonstrate the similarity of our two-dimensional extension of the limiter.

\subsection{Long wave resonance in a paraboloid basin}

The following two test cases particularly address the correct representation of a moving shoreline. They were originally defined in \citet{Thacker1981} and have an analytical solution. The first problem is a purely radially symmetric flow. Here, we also discuss the impact of the wet/dry tolerance on our method. Note, that we work with a scaled version of the problem as given in \citet{Lynett2002}. In a quadratic domain $\Omega = [-4000, 4000]^2$ with a parabolic bottom topography given by $b(\Cx) = \tilde{b}(r) = H_0 \tfrac{r^2}{a^2}$ where $r = |\Cx| = \sqrt{x^2 + y^2}$, the initial fluid depth and velocity are prescribed by
\begin{align*}
  h(\Cx,0) &= \max\left\{ 0,
    H_0 \left(\frac{\sqrt{1-A^2}}{1-A} - \frac{|\Cx|^2 (1-A^2)}{a^2 (1-A)^2}\right)\right\} \\
  \vu(\Cx,0) &= \fatvec{0}
\end{align*}
where
\begin{equation*}
  A = \frac{a^4 - r_0^4}{a^4 + r_0^4} ,
\end{equation*}
$H_0 = 1$, $r_0 = 2000$, $a = 2500$. The exact radially symmetric solution is then given by
\begin{align*}
  h(\Cx, t) &= \max\left\{ 0,
    H_0 \left(\frac{\sqrt{1-A^2}}{1-A \cos(\omega t)} - \frac{|\Cx|^2 (1-A^2)}{a^2 (1-A \cos(\omega t))^2}\right)\right\} \\
  (u,v)(\Cx,t) &=
    \begin{cases}
      \displaystyle\frac{\omega A \sin(\omega t)}{2 (1-A \cos(\omega t))} \, \Cx & \text{if } h(\Cx, t) > 0\\
      \fatvec{0} & \text{otherwise,}
    \end{cases}
\end{align*}
where $\omega$ is the frequency defined as $\omega = \sqrt{8 g H_0} / a$.

The simulations are run for two periods ($P$) of the oscillation, i.e.\ until $\Tend = 2P = 2\cdot(2\pi/\omega)$, with a time step of $\Delta t = P/700 \approx 2.534$, and a spatial resolution of $\Delta x = 125/\sqrt{2} \approx 88.39$ (leg of right angled triangle). The initial Courant number is approximately $0.16$. This is lower than the theoretically maximal Courant number because of possibly occurring spurious velocities in nearly dry regions affecting the Courant number at later times of the simulation. The maximum Courant number that is obtained for $\tolwet = 10^{-14}$ is $0.22$, whereas it is 0.16 for $\tolwet = 10^{-2}$.

\begin{figure}
  \centering
  \includegraphics[width=0.33\textwidth]{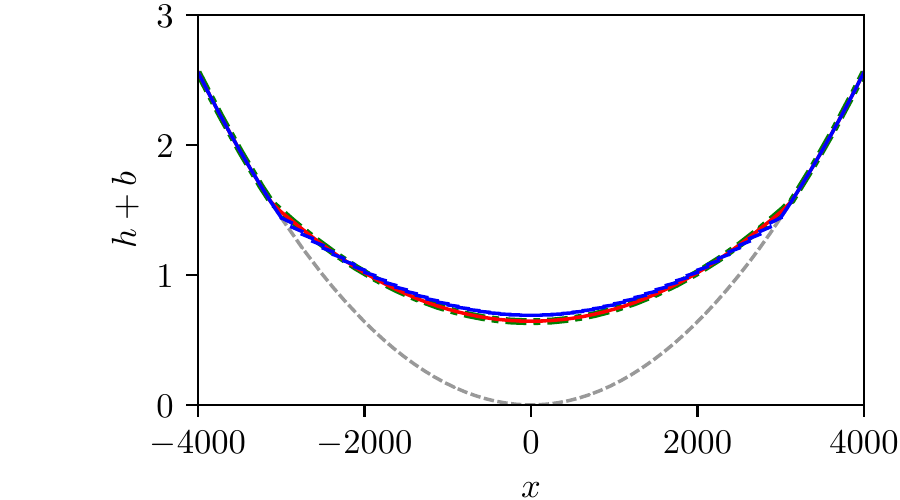}%
  \includegraphics[width=0.33\textwidth]{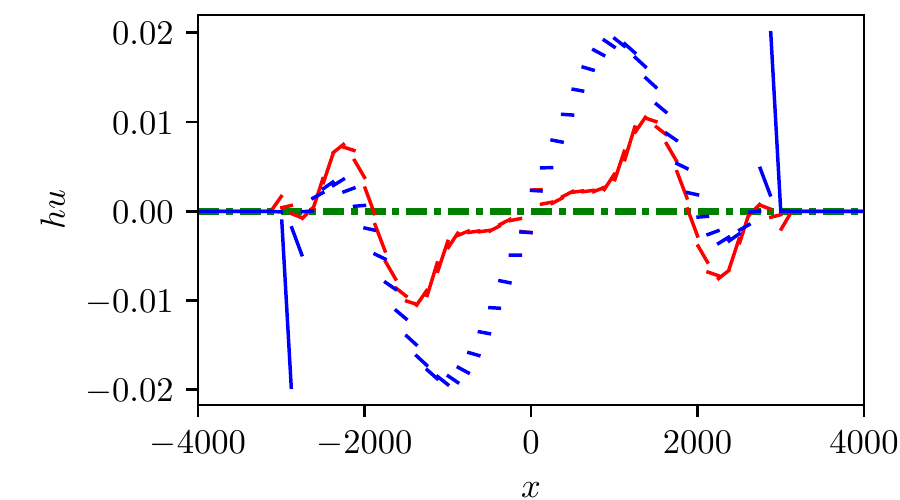}%
  \includegraphics[width=0.33\textwidth]{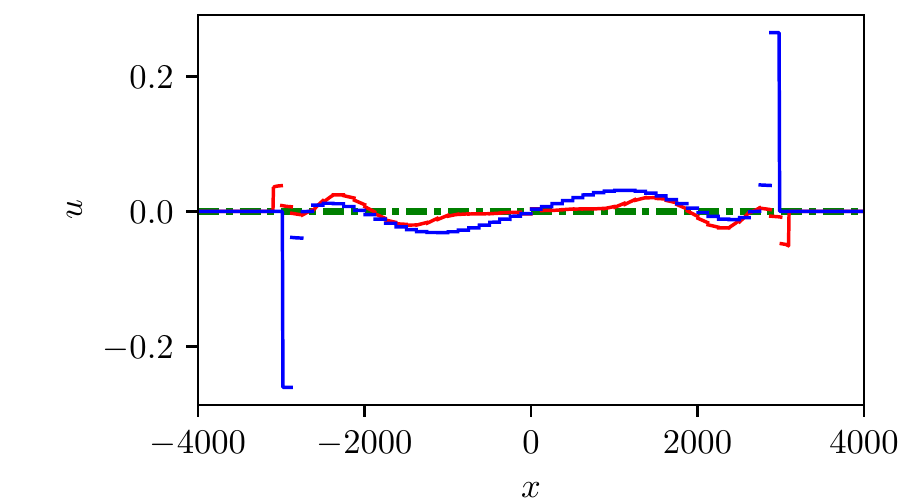}\\
  \includegraphics[width=0.33\textwidth]{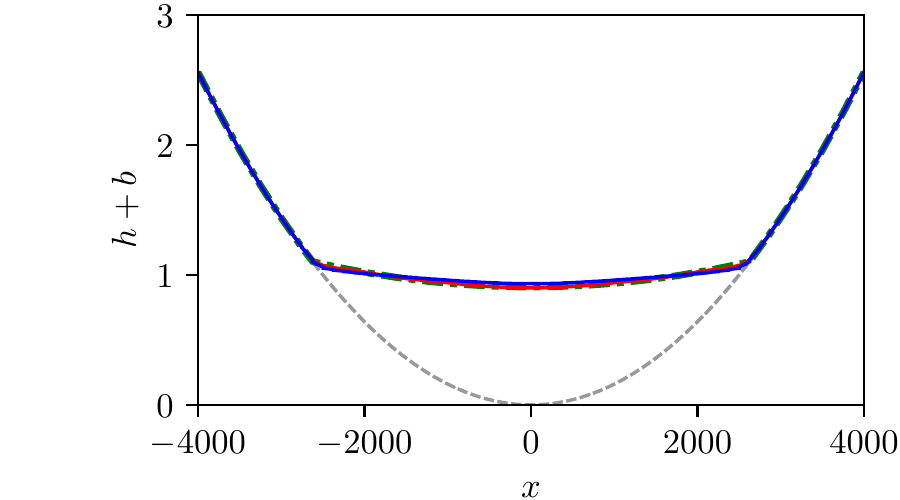}%
  \includegraphics[width=0.33\textwidth]{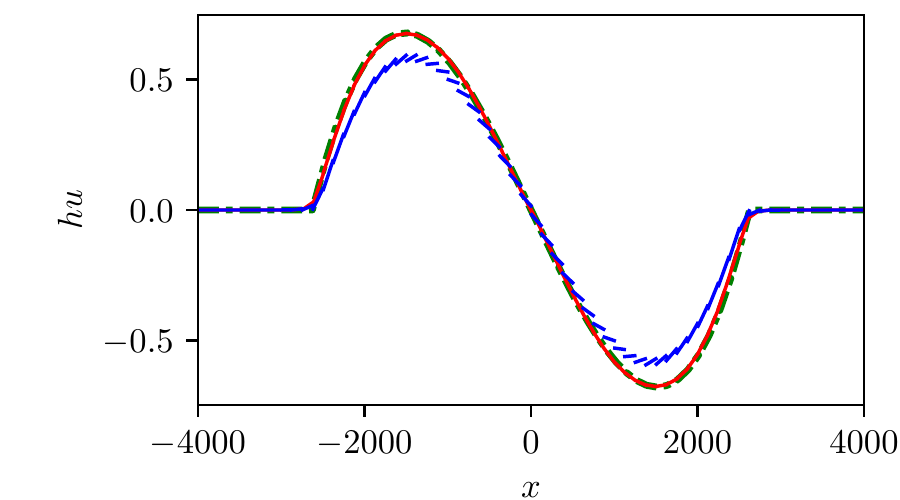}%
  \includegraphics[width=0.33\textwidth]{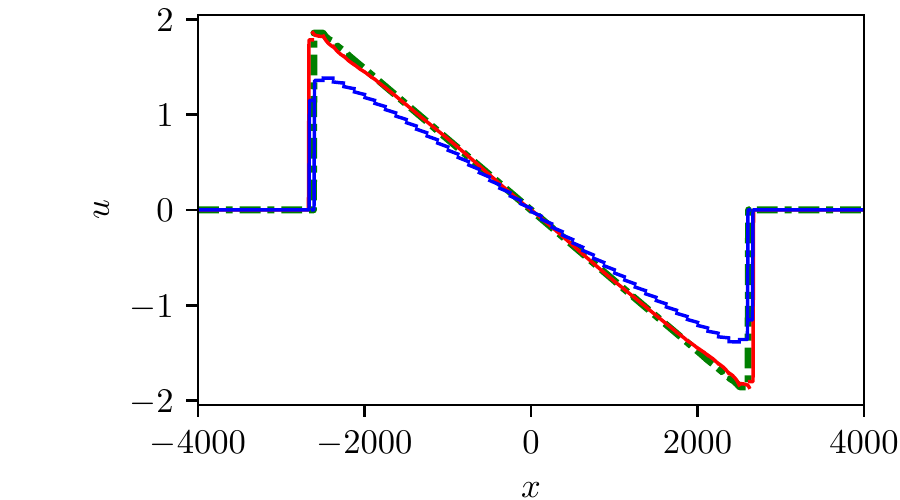}\\
  \includegraphics[width=0.33\textwidth]{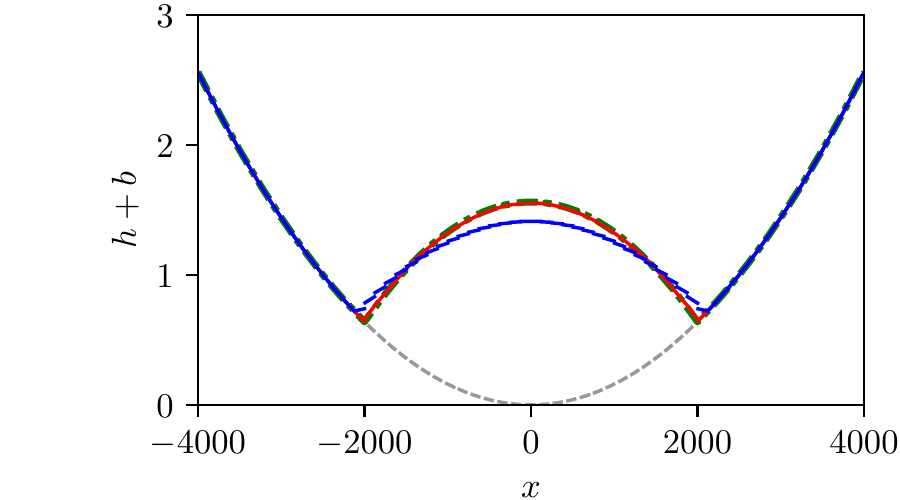}%
  \includegraphics[width=0.33\textwidth]{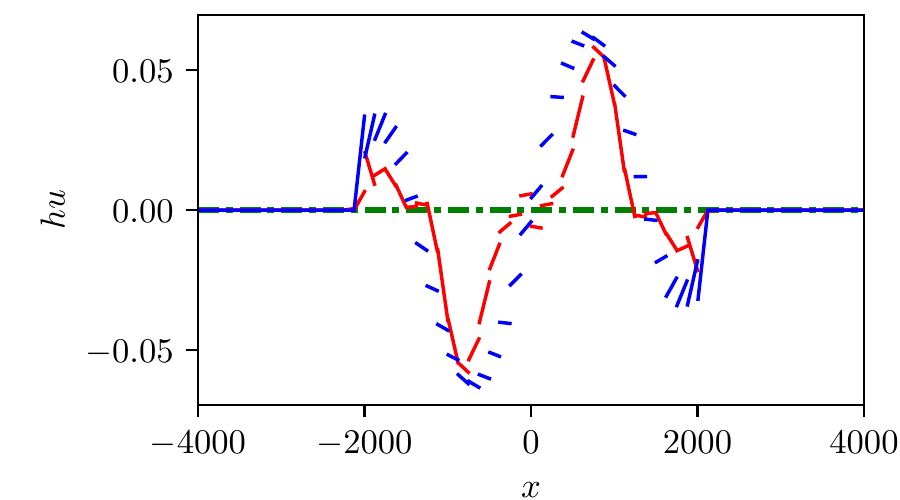}%
  \includegraphics[width=0.33\textwidth]{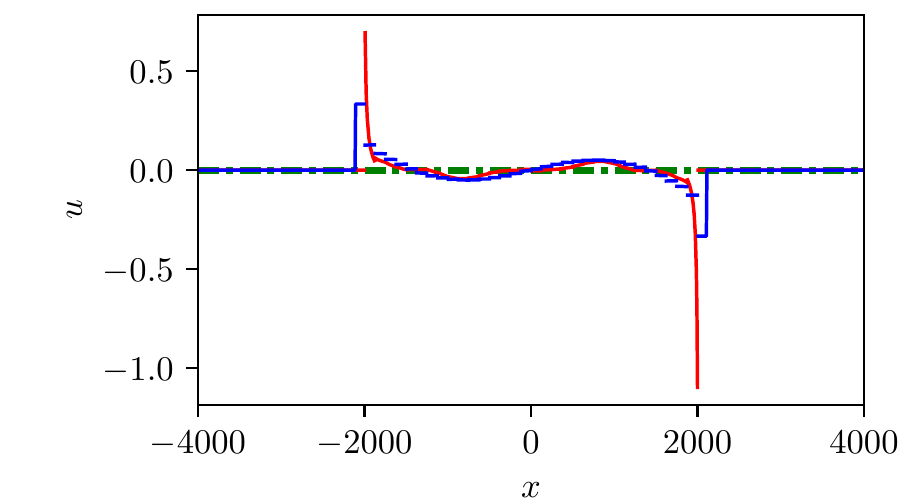}%
  \caption{Long wave resonance in a paraboloid basin: cross section of fluid depth, $x$-momentum and $x$-velocity (left to right) at times $t=1.5P$, $t=1.75P$ and $t=2P$ (top to bottom) at $y=0$. Exact solution (green dash-dotted), vertex-based limiter (red), edge-based limiter (blue), $\tolwet=10^{-2}$.\label{fig:thacker1_crossx_cmplim}}
\end{figure}

The results for fluid depth, $x$-momentum and $x$-velocity at times $t=1.5P,1.75P$ and $2P$ over a cross section $y=0$ are shown in figure \ref{fig:thacker1_crossx_cmplim}. Qualitatively, the vertex-based limiter (red line) shows slightly better results than the edge-based limiter (blue line). Note the small scale for the momentum at $t=1.5P$ and $t=2P$ in comparison with $t=1.75P$. The momentum plots in comparison with the velocity plots also show the action of the limiter: While the momentum, especially for the edge-based limiter, is non-monotone in some regions, the velocity is mostly monotone. Only near the wet/dry interface some spurious velocities are visible when the overall velocity is close to zero.

\begin{figure}
  \centering
  \includegraphics[width=0.33\textwidth]{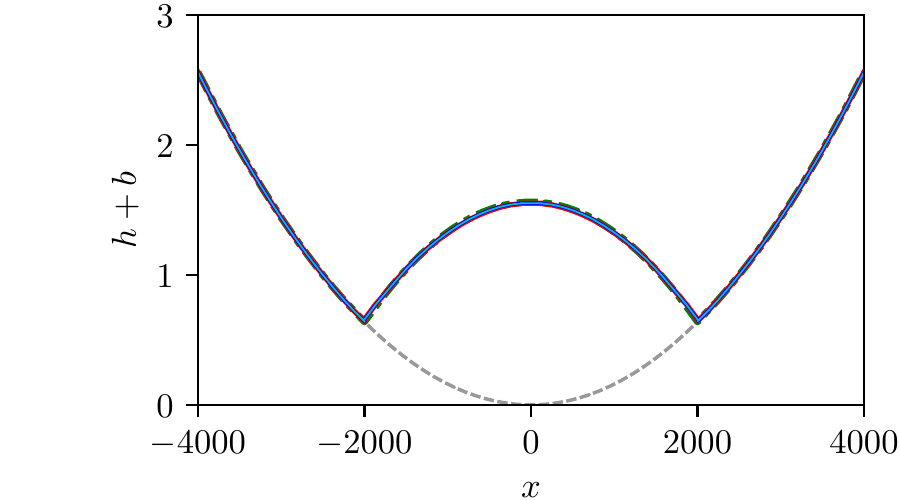}%
  \includegraphics[width=0.33\textwidth]{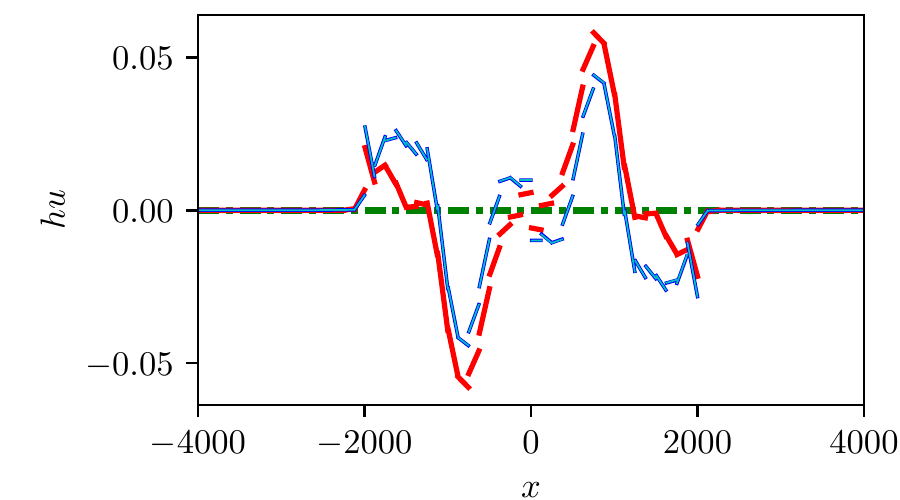}%
  \includegraphics[width=0.33\textwidth]{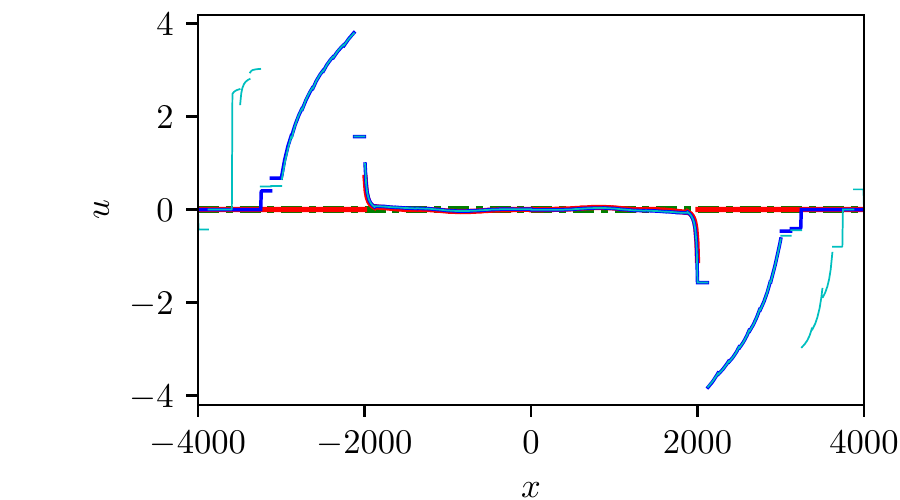}%
  \caption{Long wave resonance in a paraboloid basin: cross section of fluid depth, $x$-momentum and $x$-velocity (left to right) at time $t=2P$ at $y=0$ for the vertex-based limiter. Exact solution (green dash-dotted), solution with $\tolwet= 10^{-2}$ (red), with $\tolwet= 10^{-8}$ (blue) and with $\tolwet= 10^{-14}$ (cyan).\label{fig:thacker1_crossx_wdTOL}}
\end{figure}

\begin{figure}
  \centering
  \includegraphics[width=0.4\textwidth]{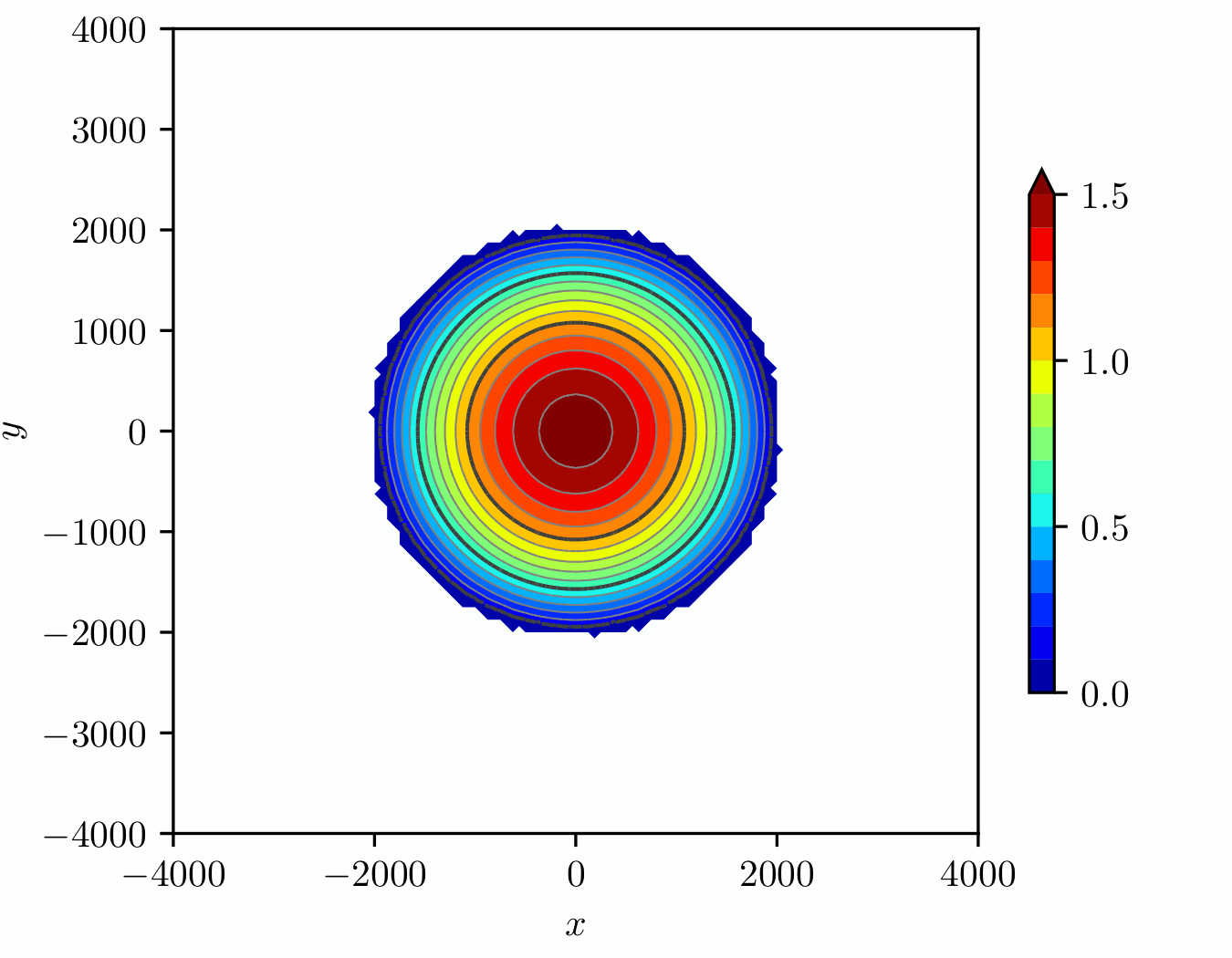}%
  \includegraphics[width=0.4\textwidth]{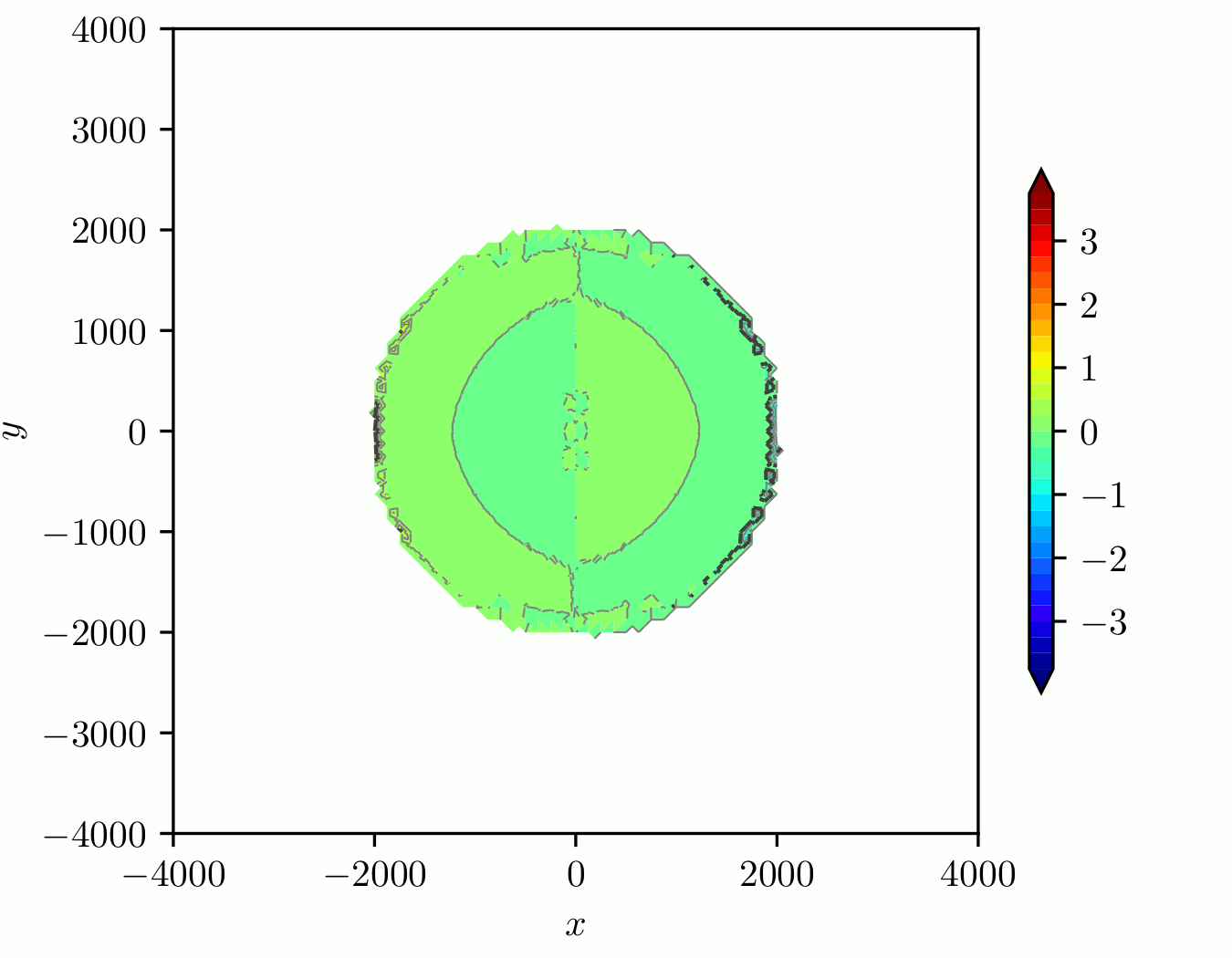}\\
  \includegraphics[width=0.4\textwidth]{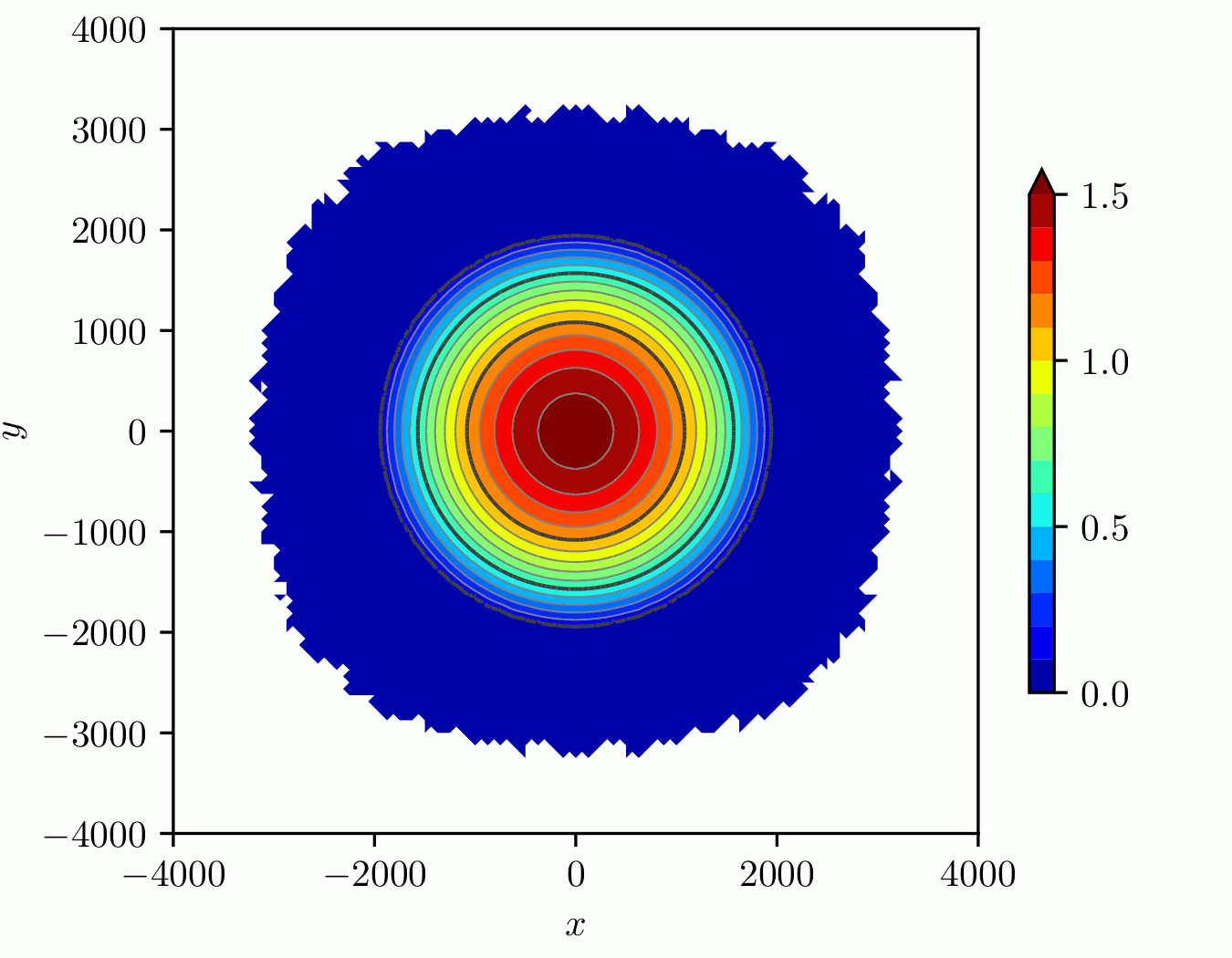}%
  \includegraphics[width=0.4\textwidth]{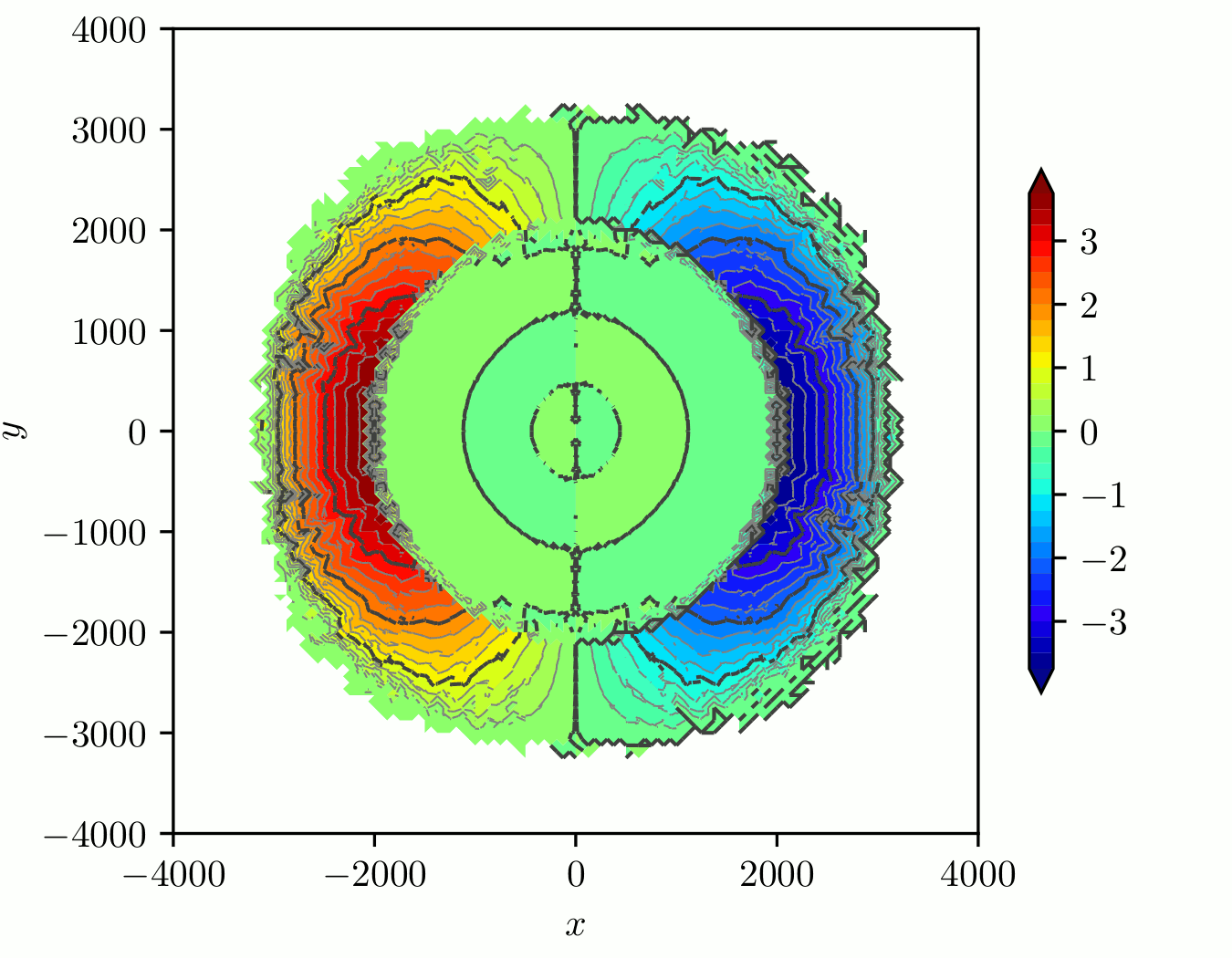}
  \caption{Long wave resonance in a paraboloid basin: 2d view of fluid depth (left) and $x$-velocity (right) at time $t=2P$ for the vertex-based limiter. Note, that only the area is colored, where the fluid depth is above the wet/dry tolerance. $\tolwet=10^{-2}$ (top) and $\tolwet=10^{-8}$ (bottom).\label{fig:thacker1_2d}}
\end{figure}

A comparison of simulation results at final time $t=2P$ with different wet/dry tolerances $\tolwet$ reveals that the results for the prognostic variables as well as the reconstructed velocities are largely insensitive to the chosen tolerance (see figure \ref{fig:thacker1_crossx_wdTOL}). However, as the wet/dry tolerance gets small a larger area is considered wet by the scheme. This is visible in the velocity plot, where spurious velocities start to appear in nearly dry regions.
We further illustrate the effect of the parameter $\tolwet$ by comparison of fully two-dimensional fields obtained with the vertex-based limiter (figure \ref{fig:thacker1_2d}). The top and bottom rows show the fluid depth and velocity with $\tolwet = 10^{-2}$ and $\tolwet = 10^{-8}$, respectively. The results for the fluid depth are largely identical with the exception that the area that the model recognizes as ``wet'' is much larger with a smaller tolerance. In the additional wet area obtained with a smaller tolerance, small values of momentum and fluid depth lead to spurious velocities. However, we note that in spite of the observed existence of spurious velocities, these are still bounded and their magnitudes are within the range of the exact solution to the problem.

\begin{figure}
 \centering
 \includegraphics[width=0.5\textwidth]{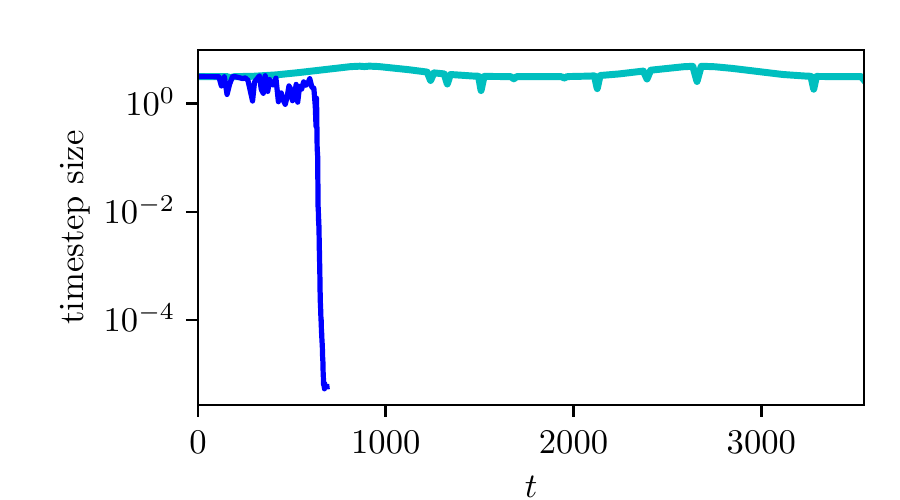}
 \caption{Long wave resonance in a paraboloid basin: Resulting time step $\Delta t$ over time by keeping $\mathrm{cfl} = 0.2$ fixed. Results with velocity-based limiting of the momentum (cyan) and with direct limiting in momentum (blue).\label{fig:thacker1_dt}}
\end{figure}

A major aspect of the velocity-based limiter becomes apparent when compared to the non-velocity-based version of that same limiter, i.e., limiting directly in the momentum variable. In figure \ref{fig:thacker1_dt} we show the maximal possible global time step $\Delta t$ for a fixed CFL number $\mathrm{cfl} = 0.2$. We allow the time step to vary over time based on the CFL number and compute it using the numerical velocity and fluid depth. Both simulation runs use the same version of the vertex-based limiter with respect to the fluid depth. We observe that the velocity-based limiter (blue line in figure \ref{fig:thacker1_dt}) allows for a reasonable time step that does not show much variation. In contrast, if we do not use the velocity-based limiter, the simulation result shows large accumulations of spurious velocities in the first drying phase, leading to an unreasonably small time step (cyan line in figure \ref{fig:thacker1_dt}) to the extent that we were not able to finish the simulation within reasonable time.

\subsection{Oscillatory flow in a parabolic bowl}

The second test case which goes back to \citet{Thacker1981} is also defined in a parabolic bowl, but describes a circular flow with a linear surface elevation in the wet part of the domain. It is the 2D analogue of the 1D test case described in \citet{Vater2015}. Here we follow the particular setup of \citet{Gallardo2007}. In a square domain $\Omega=[-2,2]^2$ with bottom topography $b = b(\Cx) = 0.1 \left(x^2 + y^2 \right)$, an analytical solution of the shallow water equations is given by
\begin{align*}
  h(\Cx,t) &= \max\left\{ 0, 0.1 \left(x \cos(\omega t) + y \sin(\omega t) + \tfrac{3}{4} \right) - b(\Cx) \right\} \\
  (u,v)(\Cx,t) &=
    \begin{cases}
      \frac{\omega}{2} \bigl( -\sin(\omega t), \cos(\omega t) \bigr) & \text{if } h(\Cx, t) > 0\\
      \fatvec{0} & \text{otherwise,}
    \end{cases}
\end{align*}
with $\omega = \sqrt{0.2 g}$.

\begin{figure}
  \centering
  \includegraphics[width=0.33\textwidth]{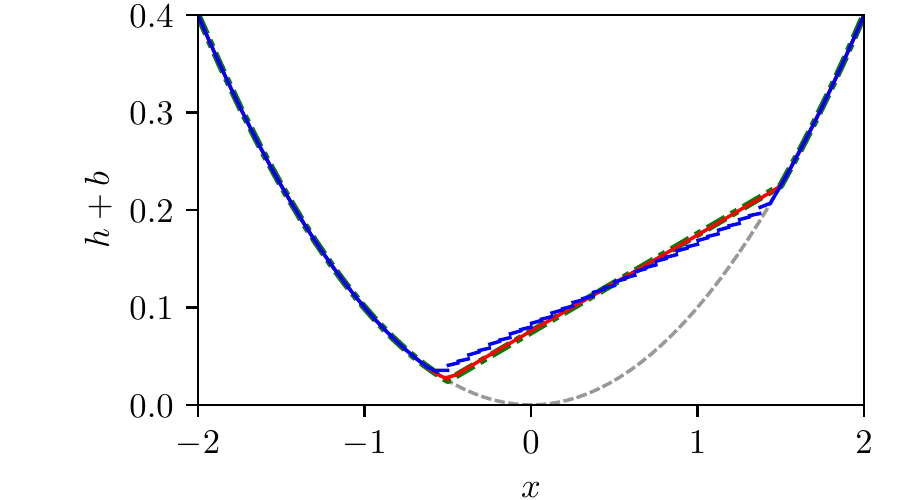}%
  \includegraphics[width=0.33\textwidth]{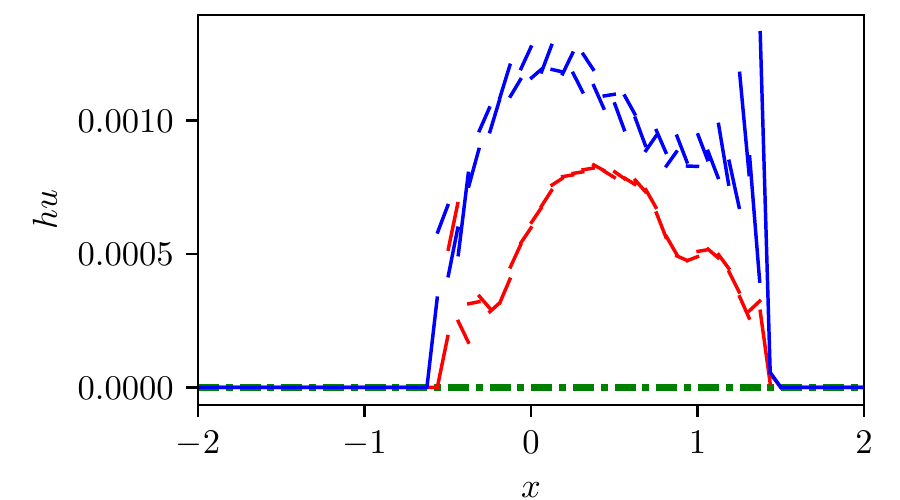}%
  \includegraphics[width=0.33\textwidth]{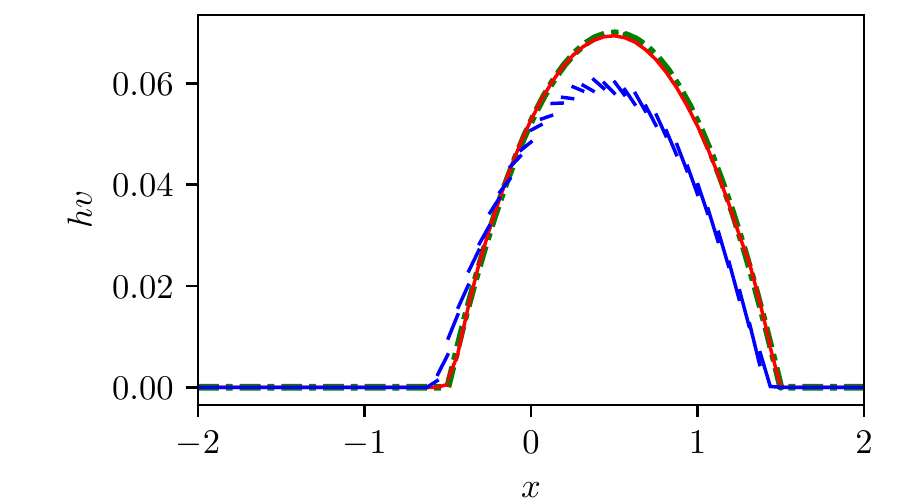}\\
  \includegraphics[width=0.33\textwidth]{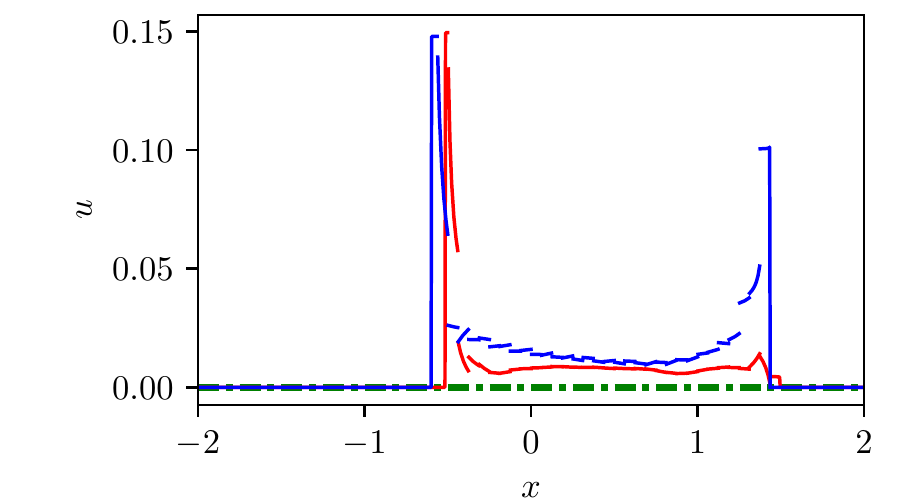}%
  \includegraphics[width=0.33\textwidth]{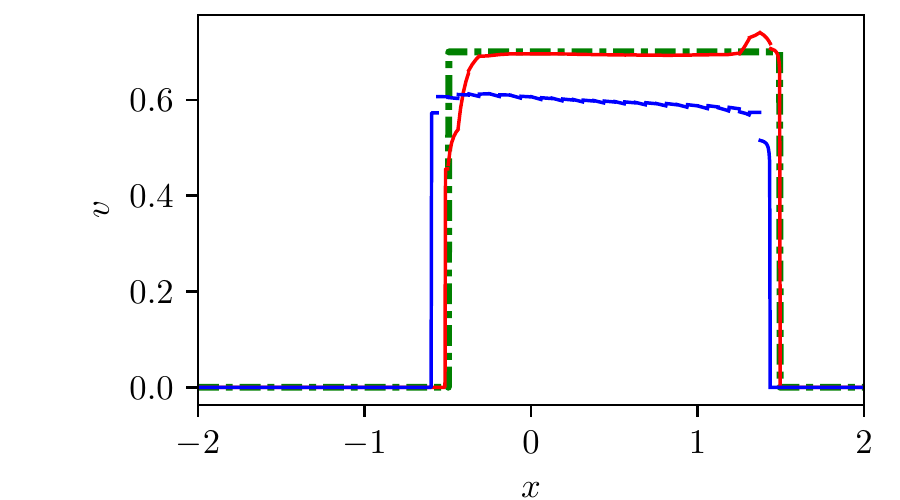}
  \caption{Oscillatory flow: cross section at $y=0$ for time $t=2P$. fluid depth, $x$-momentum, $y$-momentum, $x$-velocity, $y$-velocity (left to right, top to bottom). Exact solution (green dash-dotted), vertex-based limiter (red), edge-based limiter (blue), $\tolwet=10^{-3}$.\label{fig:thacker2_crossx}}
\end{figure}

\begin{figure}
  \centering
  \includegraphics[width=0.33\textwidth]{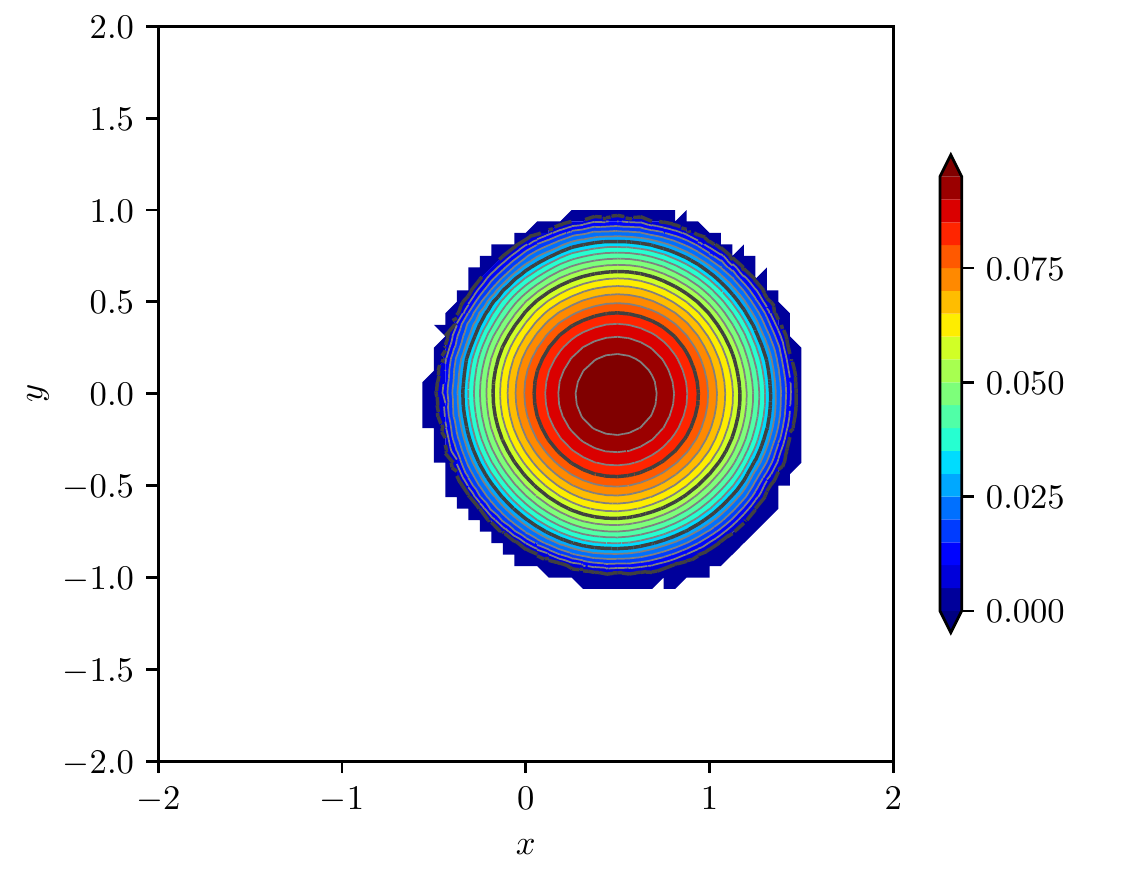}%
  \includegraphics[width=0.33\textwidth]{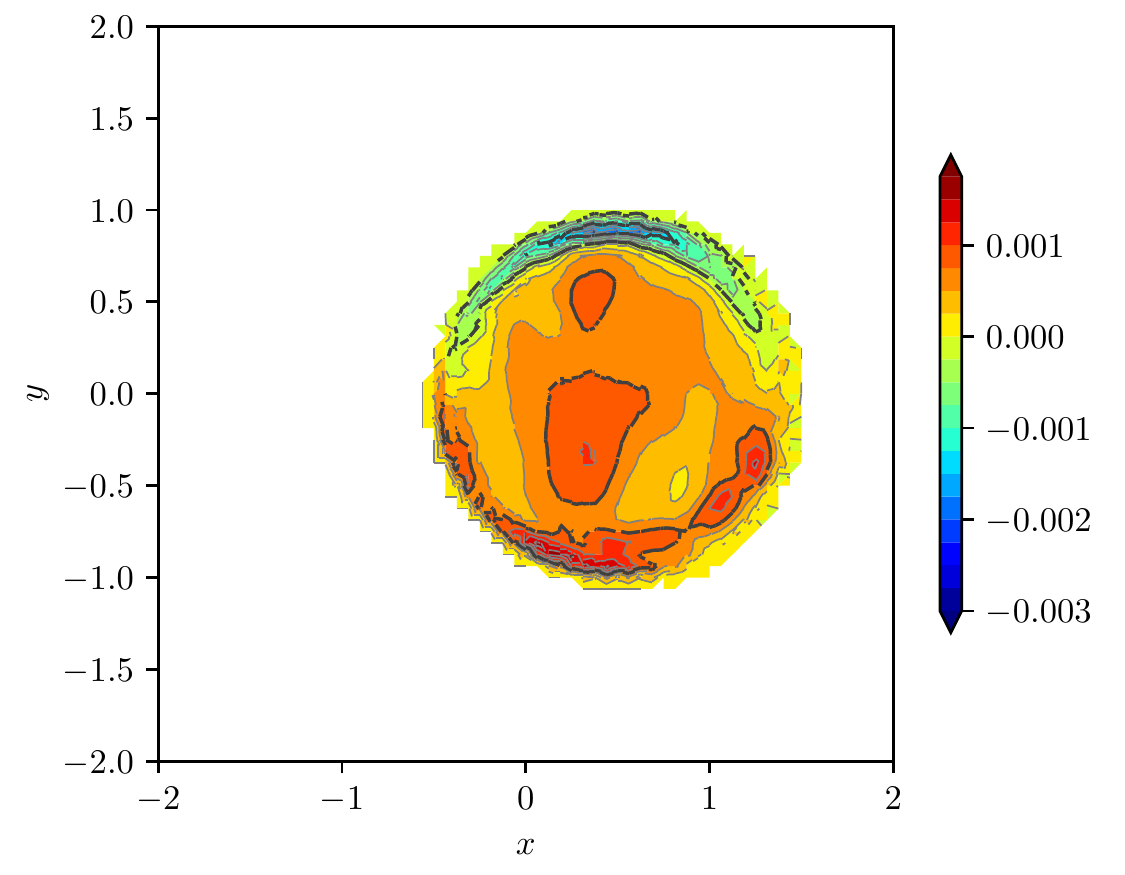}%
  \includegraphics[width=0.33\textwidth]{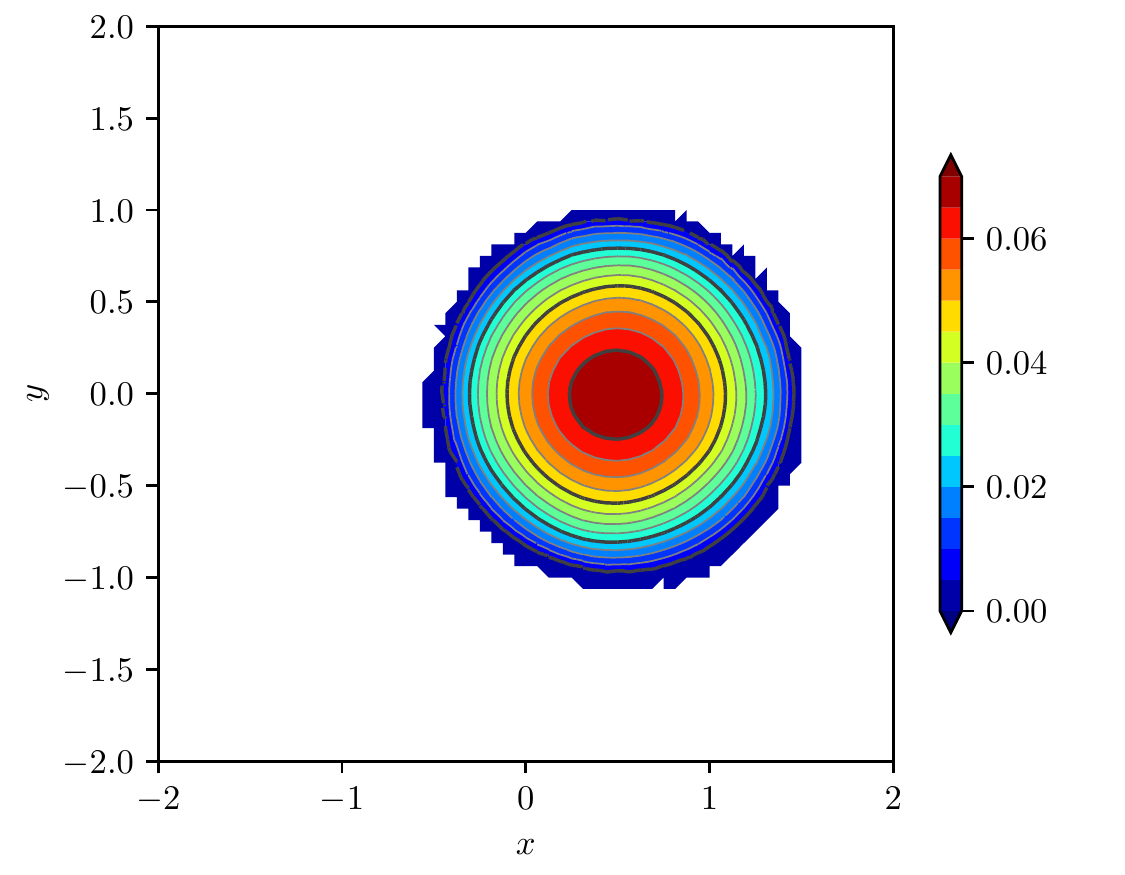}\\
  \includegraphics[width=0.33\textwidth]{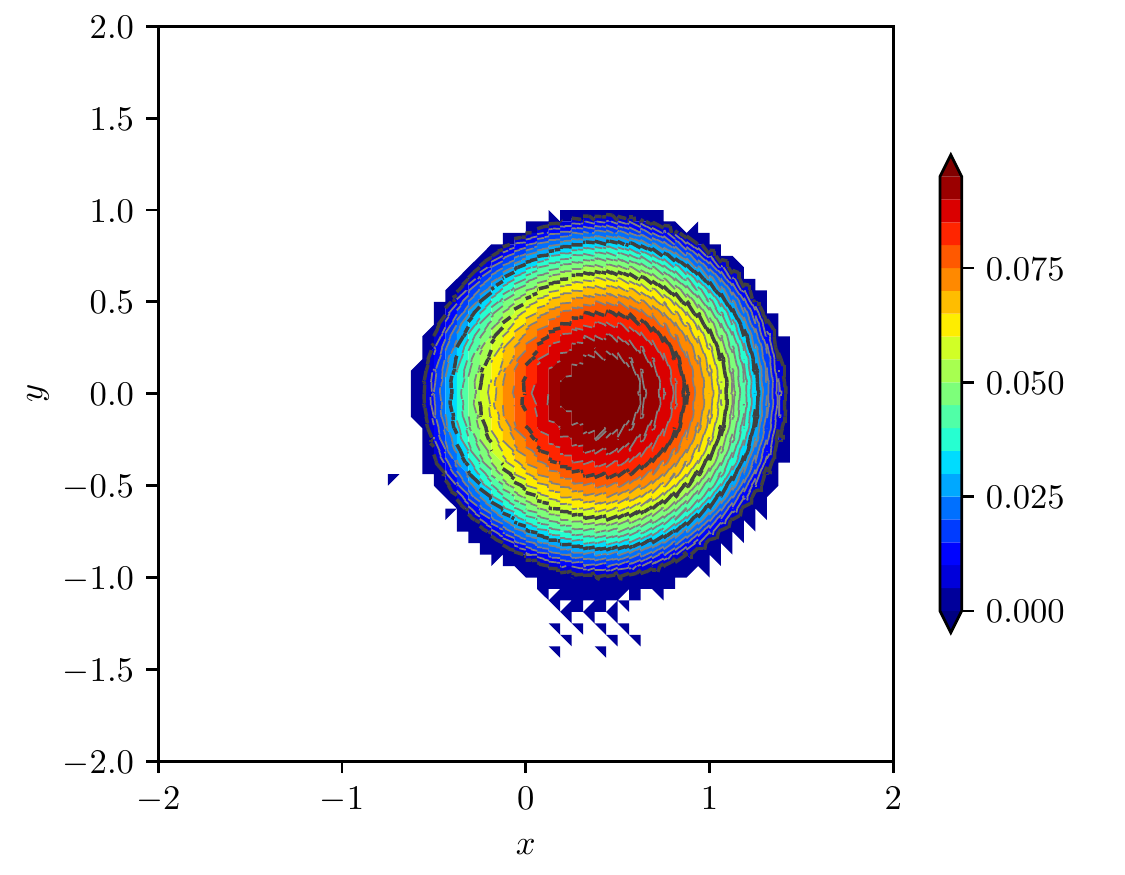}%
  \includegraphics[width=0.33\textwidth]{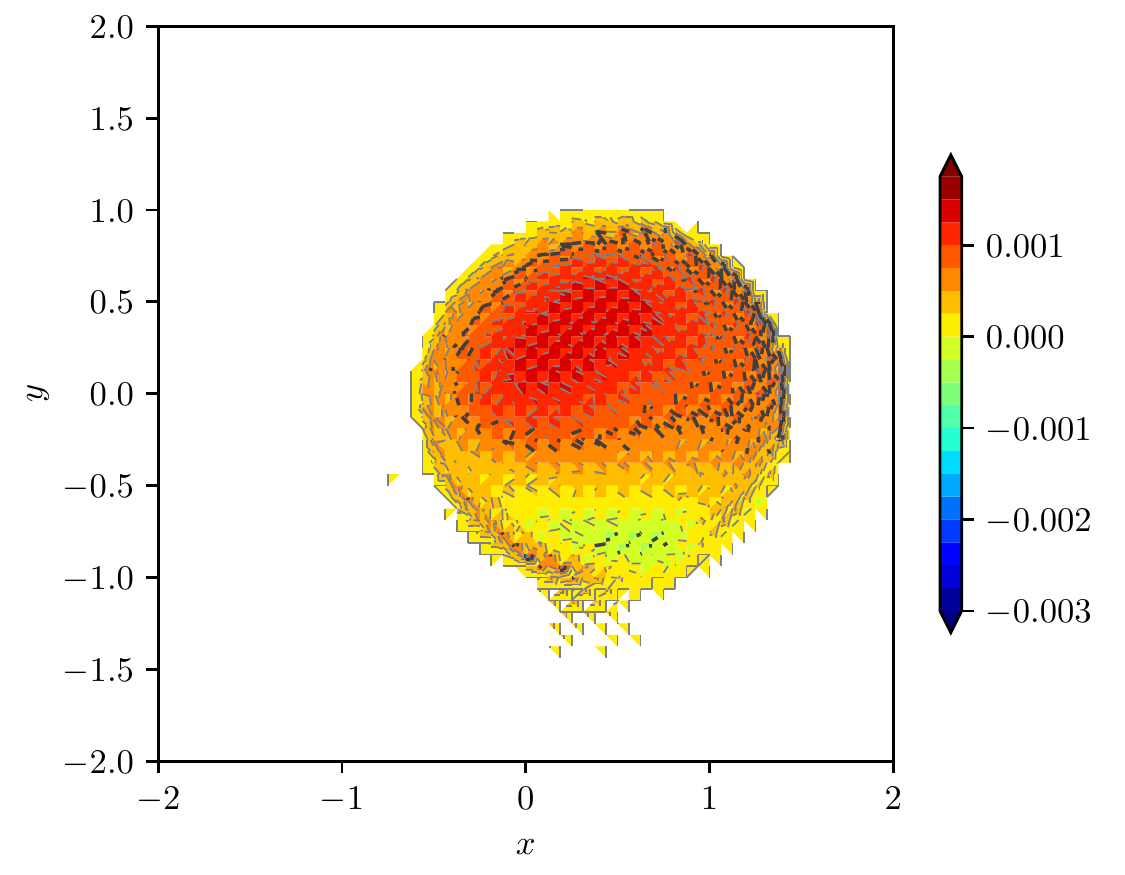}%
  \includegraphics[width=0.33\textwidth]{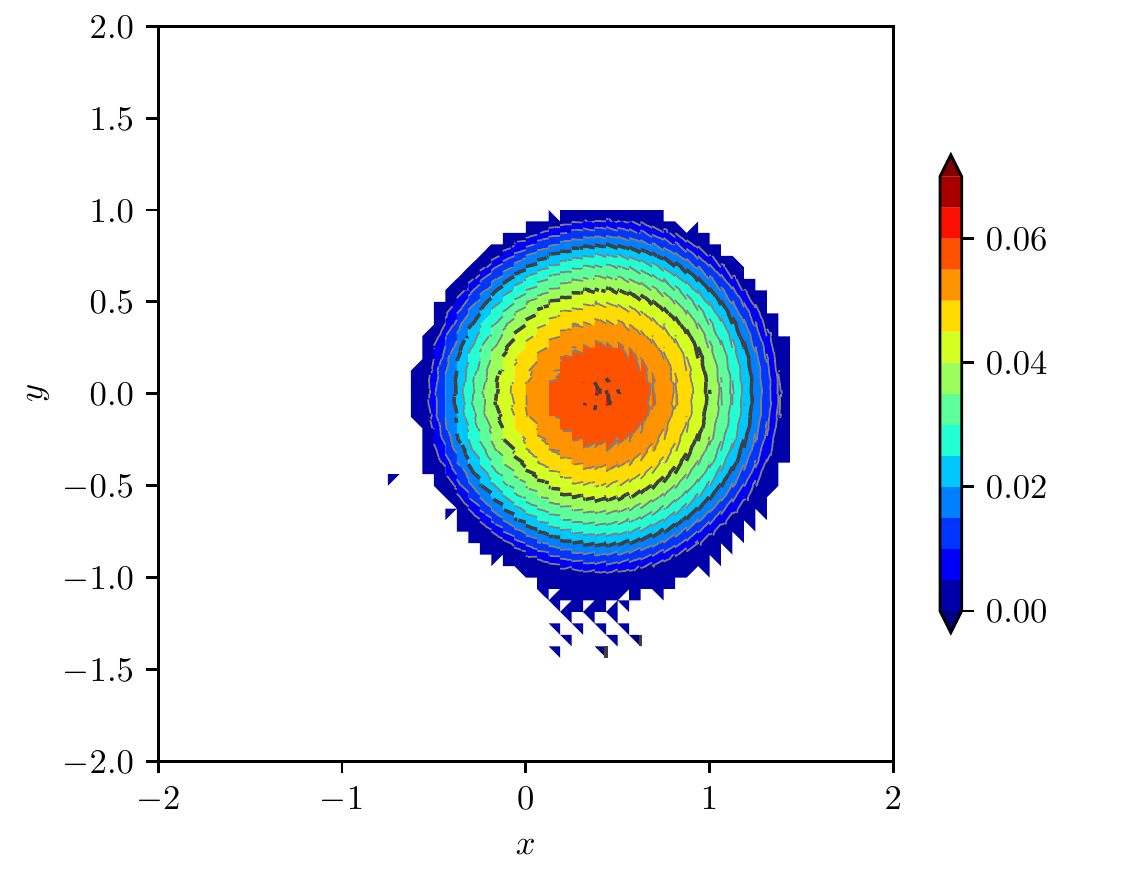}
  \caption{Oscillatory flow: 2d view for vertex-based limiter (top) and edge-based limiter (bottom). Fluid depth, $x$-momentum, $y$-momentum (left to right) at time $t=2P$. White areas denote where the fluid depth is below the wet/dry tolerance $\tolwet=10^{-3}$.\label{fig:thacker2_2d}}
\end{figure}

Starting with $t=0$ we ran simulations for two periods until $\Tend = 2P = 2\cdot(2\pi/\omega)$ of the oscillation with a time step of $\Delta t = P/1000 \approx 0.004487$. The spatial resolution is set to 8\,192 elements, which is a Cartesian grid with $64^2$ squares divided into two triangles of an edge length of $0.0625$ (leg of right angled triangle).
Figure \ref{fig:thacker2_crossx} shows cross sections over the line $y=0$ at time $t=2P$. The exact solution is plotted in green and the numerical approximation in red and blue for the vertex-based and edge-based limiter, respectively. The tolerance is chosen to $\tolwet=10^{-3}$. We observe good agreement of the numerical results with the analytical solution for fluid depth, momentum and velocity. Note, however, that the edge-based limiter tends to produce artificial discontinuities in the solution and to slightly under-predict the $y$-momentum due to a higher inherent diffusion, which results from the edge-based stencil. This also yields a too small velocity and is visible in the 2D plots in figure \ref{fig:thacker2_2d}. Moreover, the contour plot of the $x$-momentum (middle column) shows that the triggered discontinuities are clearly visible. The results with the vertex-based limiter (top) are smoother and show less diffusion.

\begin{figure}
  \begin{minipage}[c]{0.6\textwidth}
    \includegraphics[width=\textwidth]{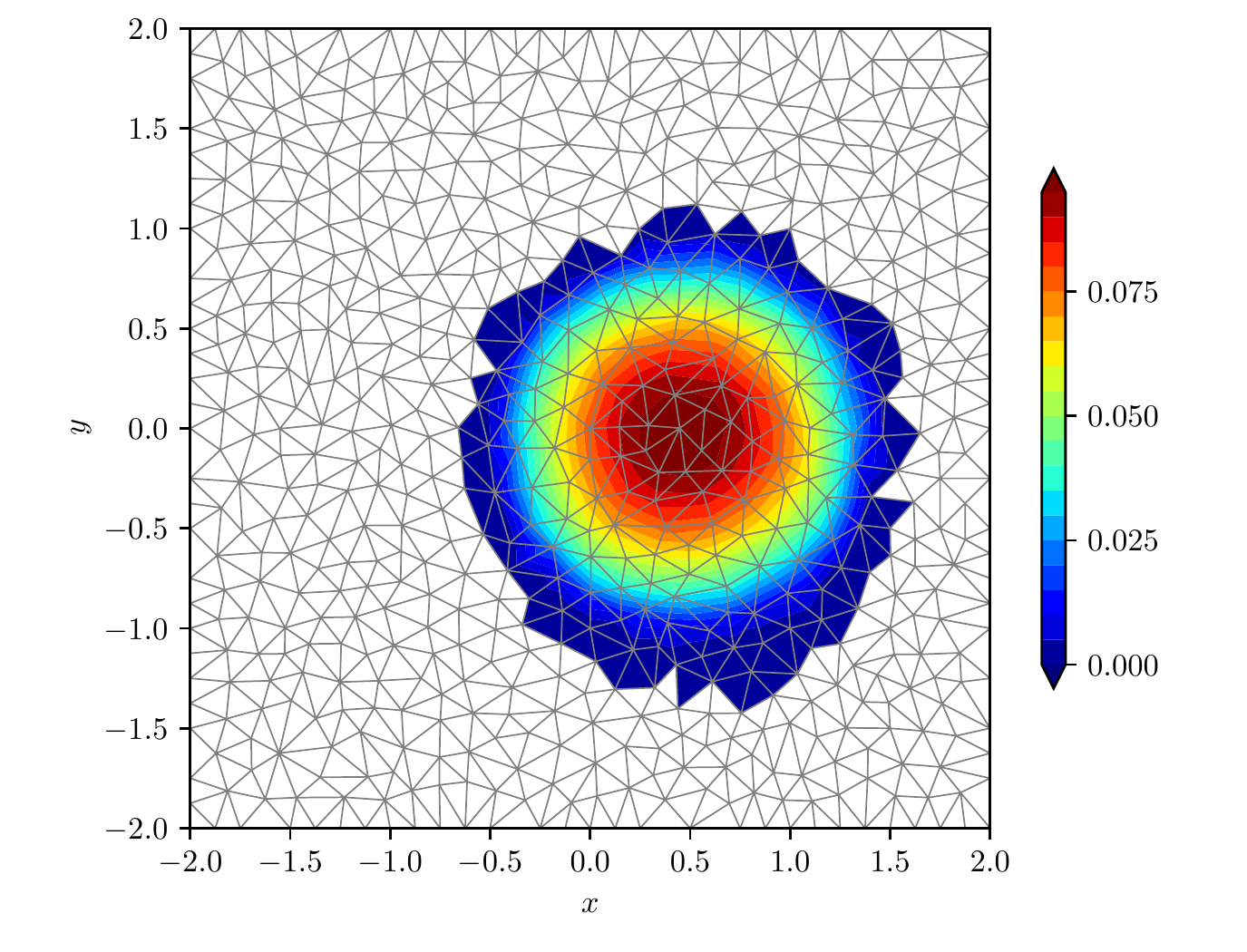}
  \end{minipage}
  \begin{minipage}[c]{0.4\textwidth}
    \includegraphics[width=\textwidth]{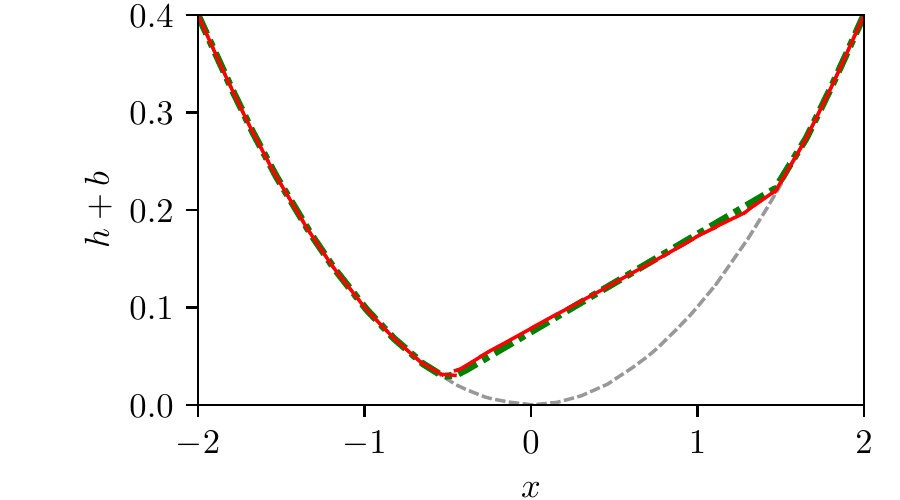}
    \includegraphics[width=\textwidth]{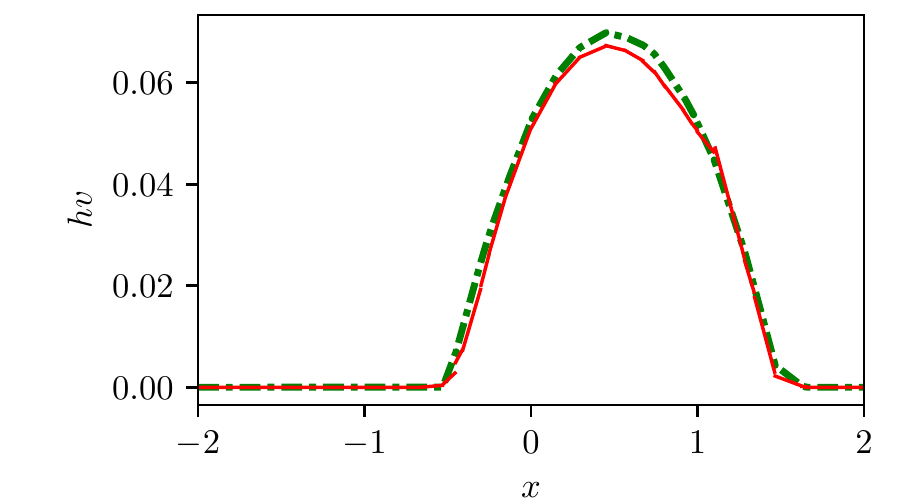}
  \end{minipage}
  \caption{Oscillatory flow on unstructured grid: 2d view of fluid depth with grid (left) and cross section at $y=0$ for fluid depth and $y$-momentum (right) at time $t=2P$. Exact solution (green dash-dotted), vertex-based limiter (red), $\tolwet=10^{-3}$.\label{fig:thacker2_2dunstruct}}
\end{figure}

To show that our limiting approach is applicable to arbitrary grids, figure \ref{fig:thacker2_2dunstruct} shows analogue simulation results on a highly unstructured Delaunay grid with 1233 elements. Note that this grid has a coarser resolution than the one used in figures \ref{fig:thacker2_crossx} and \ref{fig:thacker2_2d}. As can be seen from the cross section plots, the results are similar to those on the other grids.

\begin{figure}
  \centering
  \includegraphics[width=0.5\textwidth]{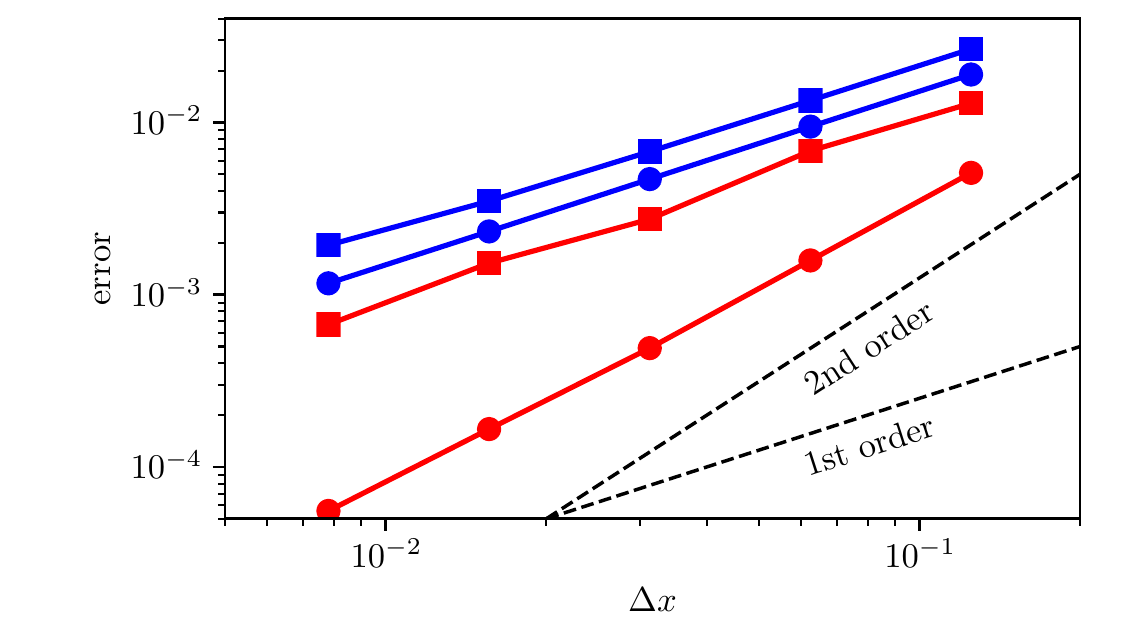}%
  \includegraphics[width=0.5\textwidth]{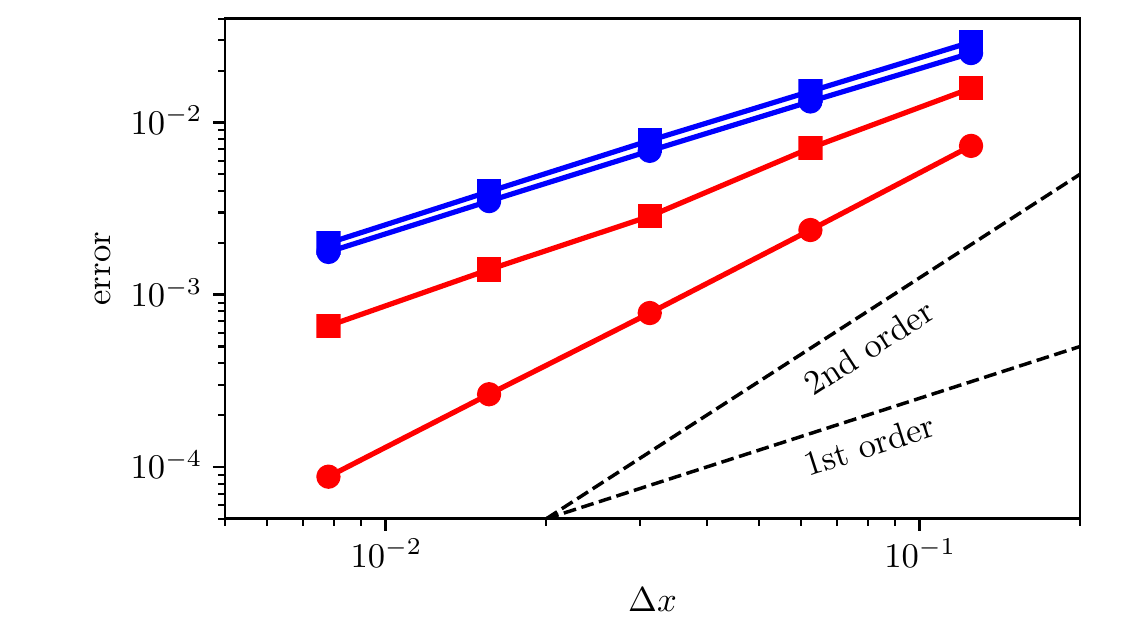}
  \caption{Oscillatory flow: errors in fluid depth (left) and momentum (right) measured in the $L^2$ norm (circles) and the $L^\infty$ norm (squares). Vertex-based limiter (red), edge-based limiter (blue).\label{fig:thacker2_convergence}}
\end{figure}

Besides the accuracy on fixed grids, also the convergence of the wetting and drying scheme is of interest. While we cannot expect second order convergence due to the non-smooth transition (kink) between wet and dry regions in the flow variables, the convergence rate should be at least approximately linear. For the convergence calculation we have computed the solution up to $t=2P$ on several grids with the number of cells ranging from $2\,048$ to $524\,288$ and fixed ratio $\Delta t / \Delta x$ and a wet/dry tolerance $\tolwet=10^{-8}$. The experimental convergence rate is then calculated by the formula
\begin{equation*}
  \gamma_c^f := \frac{\log(\|e_c\| / \|e_f\|)}{\log(\Delta x_c / \Delta x_f)} .
\end{equation*}
In this definition, $e_c$ and $e_f$ are the computed error functions of the solution on a coarse and a fine grid (denoted by the number of cells) and $\Delta x_c$ and $\Delta x_f$ are the corresponding grid resolutions. In figure \ref{fig:thacker2_convergence} and table \ref{table:thacker2_convergence} we show the results of this convergence analysis. The DG method converges with both limiters, however, the convergence rate that is achieved in the $L^2$ norm is higher with the vertex-based limiter ($\approx 1.6$) than with the edge-based limiter ($\approx 1$).

\begin{figure}
  \centering
  \includegraphics[width=0.5\textwidth]{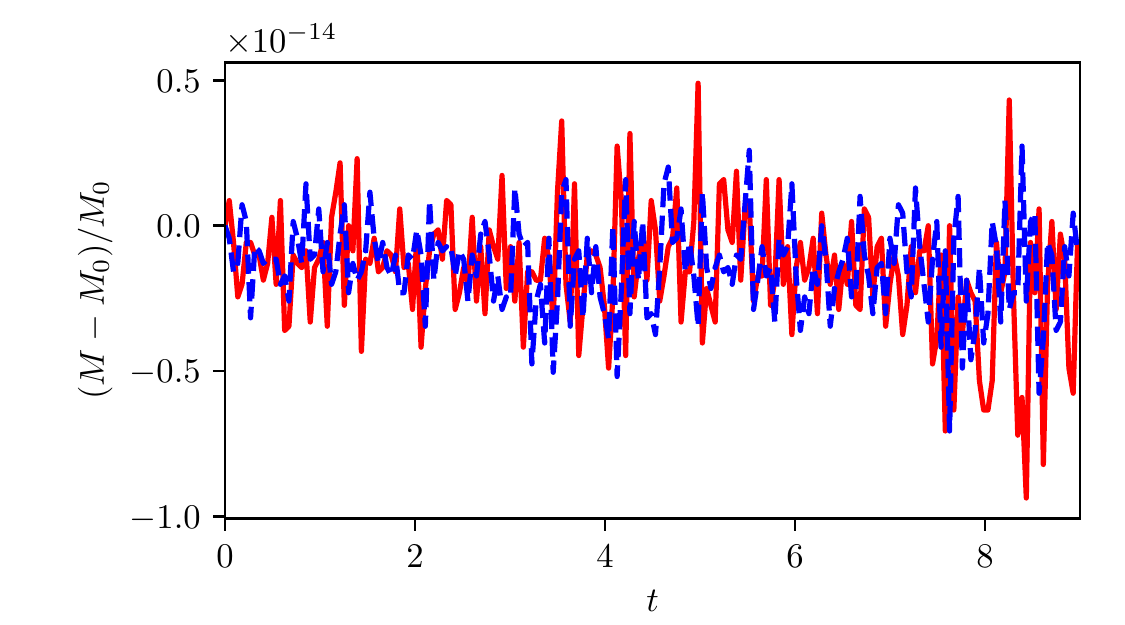}%
  \includegraphics[width=0.5\textwidth]{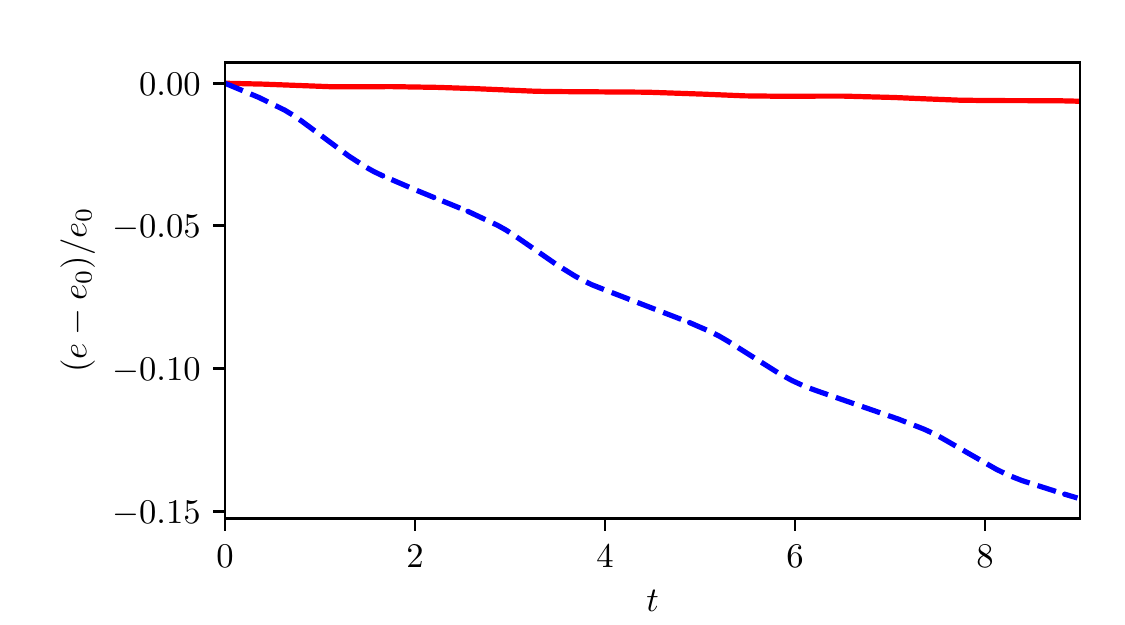}
  \caption{Oscillatory flow: time series of relative mass (left) and energy (right) changes for the vertex-based limiter (red) and edge-based limiter (blue).\label{fig:thacker2_mass}}
\end{figure}

The test case of an oscillatory flow in a parabolic bowl is also suitable to evaluate the conservation of mass and of total energy $E=\int_\Omega h(\vu\cdot\vu)/2+gh(h/2+b) \dint{x}$ for the numerical method, since there is no flow across the boundary of the domain. While mass conservation should hold up to machine accuracy, total energy can only hold within the approximation error. We can see that mass conservation is not affected by the slope limiters (left plot of figure \ref{fig:thacker2_mass}), while only the vertex-based limiter (right plot of same figure) nearly conserves energy. This indicates that the edge-based limiter exposes some numerical dissipation.

\begin{table}
  \renewcommand{\arraystretch}{1.2}
  \begin{tabular}{l|cccc}
                               & $L^2(h)$ & $L^2(m)$ & $L^\infty(h)$ & $L^\infty(m)$ \\\hline
    $\gamma_{2048}^{8192}$     & 1.6873   & 1.6230   & 0.9104        & 1.1587        \\
    $\gamma_{8192}^{32768}$    & 1.6903   & 1.5996   & 1.3190        & 1.3072        \\
    $\gamma_{32768}^{131072}$  & 1.5626   & 1.5671   & 0.8477        & 1.0294        \\
    $\gamma_{131072}^{524288}$ & 1.5779   & 1.5901   & 1.1845        & 1.0847        \\
    $\gamma_\mathrm{fitted}$   & 1.6289   & 1.5926   & 1.0690        & 1.1496
  \end{tabular}%
  \hspace{1ex}\hfill%
  \begin{tabular}{l|cccc}
                               & $L^2(h)$ & $L^2(m)$ & $L^\infty(h)$ & $L^\infty(m)$ \\\hline
    $\gamma_{2048}^{8192}$     & 1.0048   & 0.9332   & 0.9926        & 0.9494        \\
    $\gamma_{8192}^{32768}$    & 1.0125   & 0.9527   & 0.9860        & 0.9491        \\
    $\gamma_{32768}^{131072}$  & 1.0090   & 0.9694   & 0.9513        & 0.9834        \\
    $\gamma_{131072}^{524288}$ & 1.0012   & 0.9802   & 0.8538        & 0.9957        \\
    $\gamma_\mathrm{fitted}$   & 1.0077   & 0.9593   & 0.9505        & 0.9688
  \end{tabular}
  \caption{Oscillatory flow: convergence rates between different grid levels for fluid depth ($h$) and momentum ($m$) in the $L^2$ and $L^\infty$ norms. Vertex-based limiter (left) and edge-based limiter (right). Also displayed is the mean convergence rate $\gamma_\mathrm{fitted}$, which is obtained by a least squared fit.\label{table:thacker2_convergence}}
\end{table}

\begin{figure}
  \centering
  \includegraphics[width=0.5\textwidth]{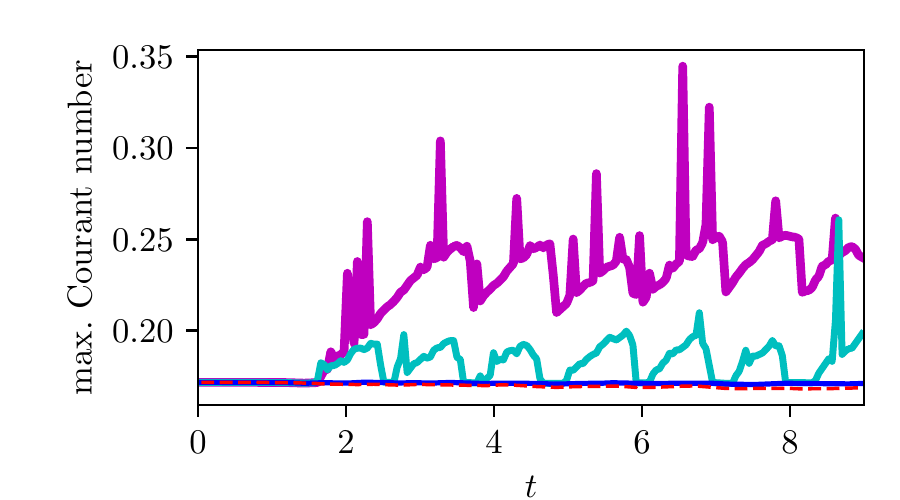}%
  \includegraphics[width=0.5\textwidth]{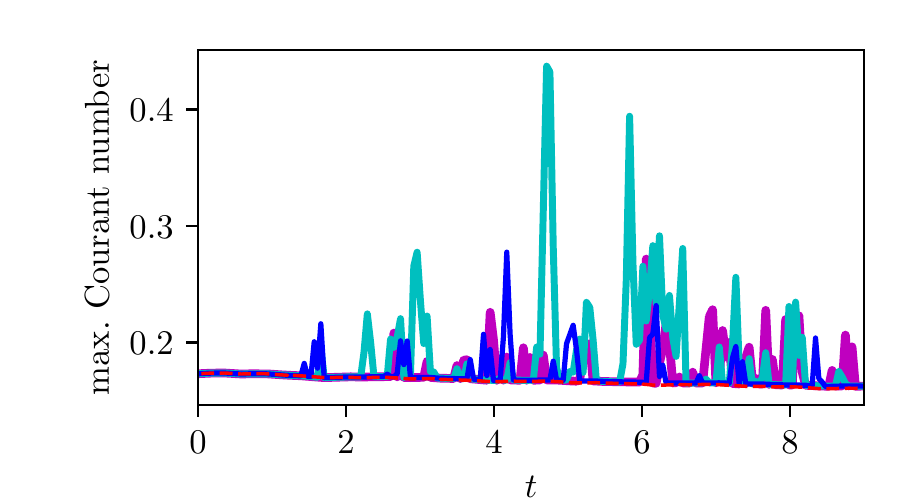}
  \caption{Oscillatory flow: time series of maximum Courant number with wet/dry tolerance $\tolwet=10^{-2}$ (red dashed), $\tolwet=10^{-4}$ (blue), $\tolwet=10^{-8}$ (cyan), $\tolwet=10^{-14}$ (magenta). Vertex-based limiter (left), edge-based limiter (right).\label{fig:thacker2_cfl}}
\end{figure}

Finally, we record the Courant number for this test case over time for different wet/dry tolerances. In figure \ref{fig:thacker2_cfl} we plot the Courant number for the vertex-based (left) and edge-based limiter (right) with wet/dry tolerances $\tolwet\in \{10^{-2},10^{-4},10^{-8},10^{-14} \}$. It can be observed that for all tolerances the Courant number stays bounded and mostly below the theoretical limit. However, when the wet/dry tolerance becomes smaller, spurious velocities start to arise and affect the Courant number. If the tolerance is set large enough ($\approx 10^{-4}$), we obtain a nearly constant Courant number over time, which is similar to the Courant number of the exact problem.

\subsection{Runup onto a complex three-dimensional beach}

\begin{figure}
  \includegraphics[width=0.45\textwidth]{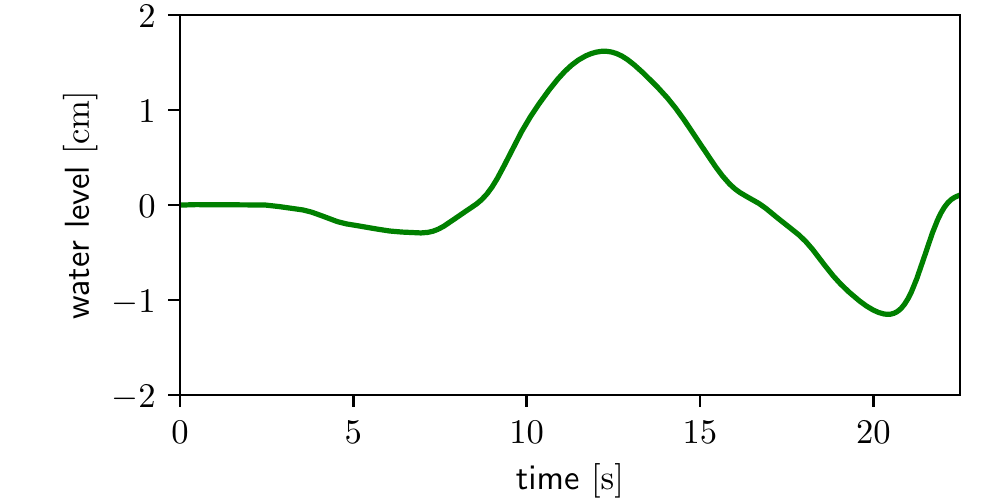}\hfill
  \includegraphics[width=0.45\textwidth]{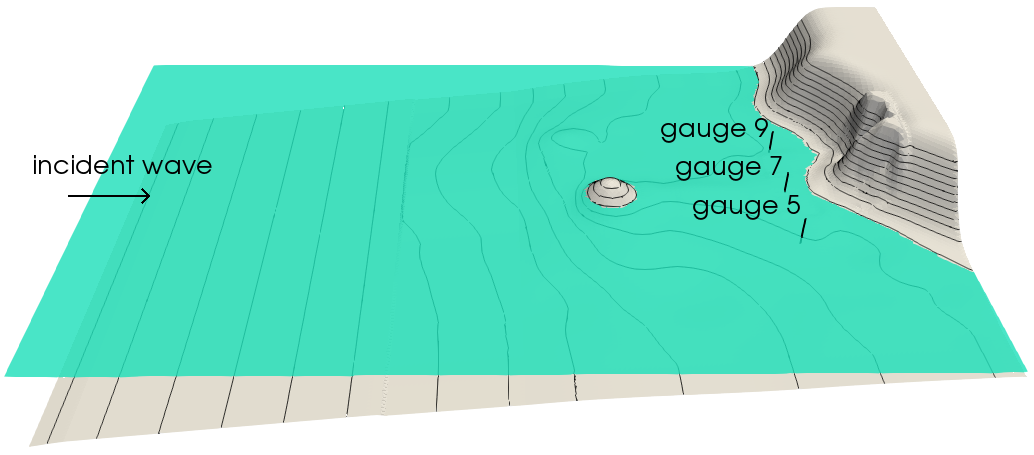}
  \caption{Okushiri: time series of incident wave which is used as boundary condition (left) and experimental setup (right).\label{fig:okushiri_setup}}
\end{figure}

The 1993 Okushiri tsunami caused many unexpected phenomena, such as an extreme runup height of 32m which was observed near the village of Monai on Okushiri island. The event was reconstructed in an 1/400 scale laboratory experiment, using a large-scale tank (205m long, 6m deep, 3.4m wide) at Central Research Institute for Electric Power Industry (CRIEPI) in Abiko, Japan \citep{Matsuyama2001}. For the test case the coastal topography in a domain of $5.448\mathrm{m} \times 3.402\mathrm{m}$ and the incident wave from offshore is provided. Beside the temporal and spatial variations of the shoreline location, the temporal evolution of the surface elevation at three specified gauge stations are of interest (figure \ref{fig:okushiri_setup}).

At the offshore boundary we set the incident wave as right-going simple wave. This means, given the fluid depth $h$ of the incident wave the $x$-velocity is computed by
\begin{equation}\label{eq:bcsimplewave}
  u = 2 (\sqrt{g h} - \sqrt{g h_0}) ,
\end{equation}
where $h_0=0.13535\mathrm{m}$ denotes the water depth at rest. At the other three boundaries there were reported to be walls. So we set reflective wall boundary conditions at these locations.

\begin{figure}
  \includegraphics[width=0.33\textwidth]{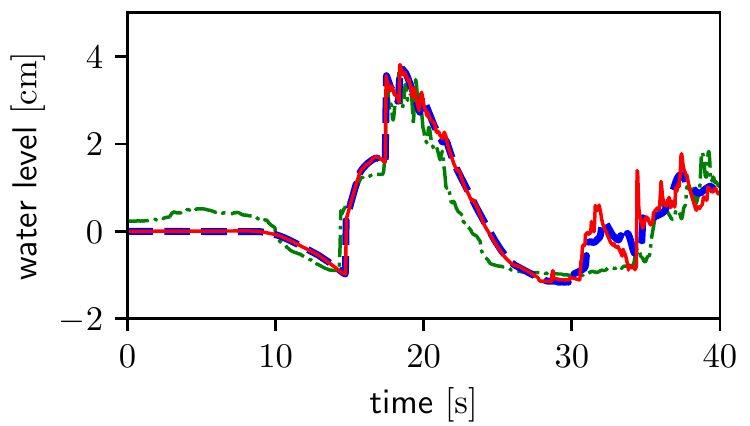}\hfill
  \includegraphics[width=0.33\textwidth]{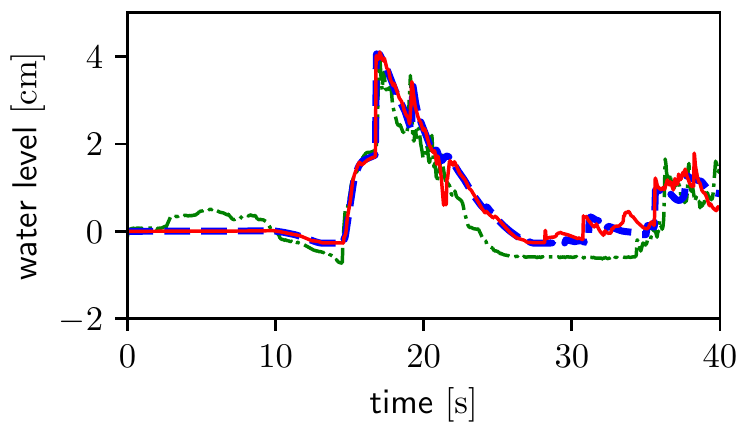}\hfill
  \includegraphics[width=0.33\textwidth]{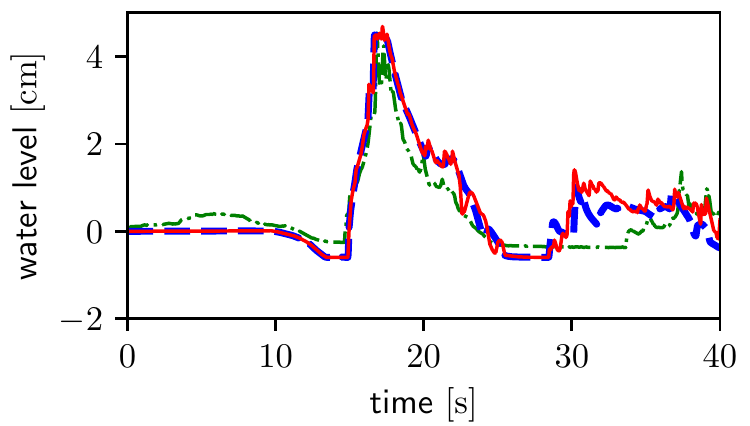}
  \caption{Okushiri: time series of gauge data: gauge 5 (left), gauge 7 (middle), gauge 9 (right), experimental data (green dash-dotted), vertex-based limiter (red), edge-based limiter (blue dashed). $\tolwet= 10^{-4}$.\label{fig:okushiri_timeseries}}
\end{figure}

\begin{figure}[t!]
  \centering
  \includegraphics[width=0.33\textwidth]{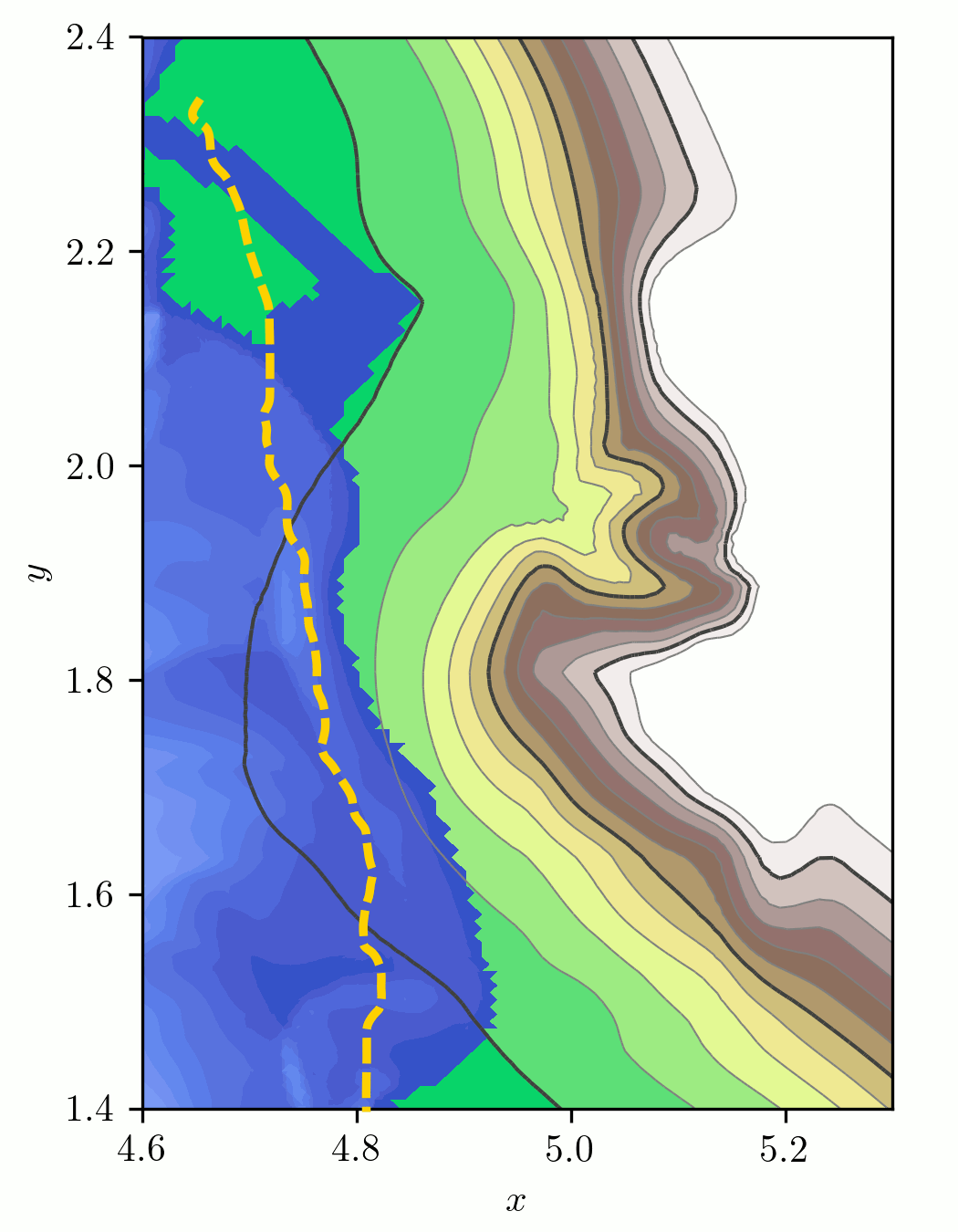}%
  \includegraphics[width=0.33\textwidth]{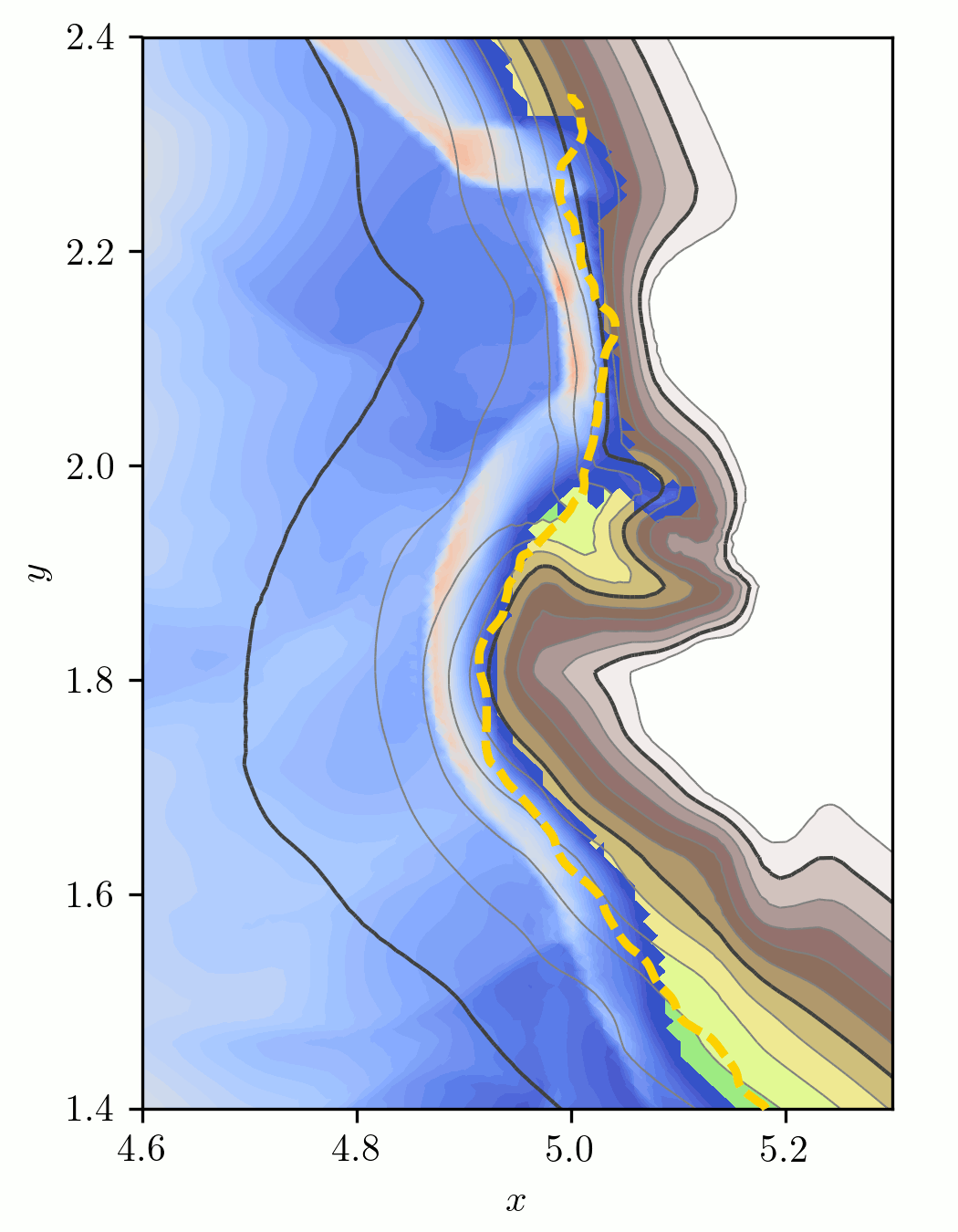}%
  \includegraphics[width=0.33\textwidth]{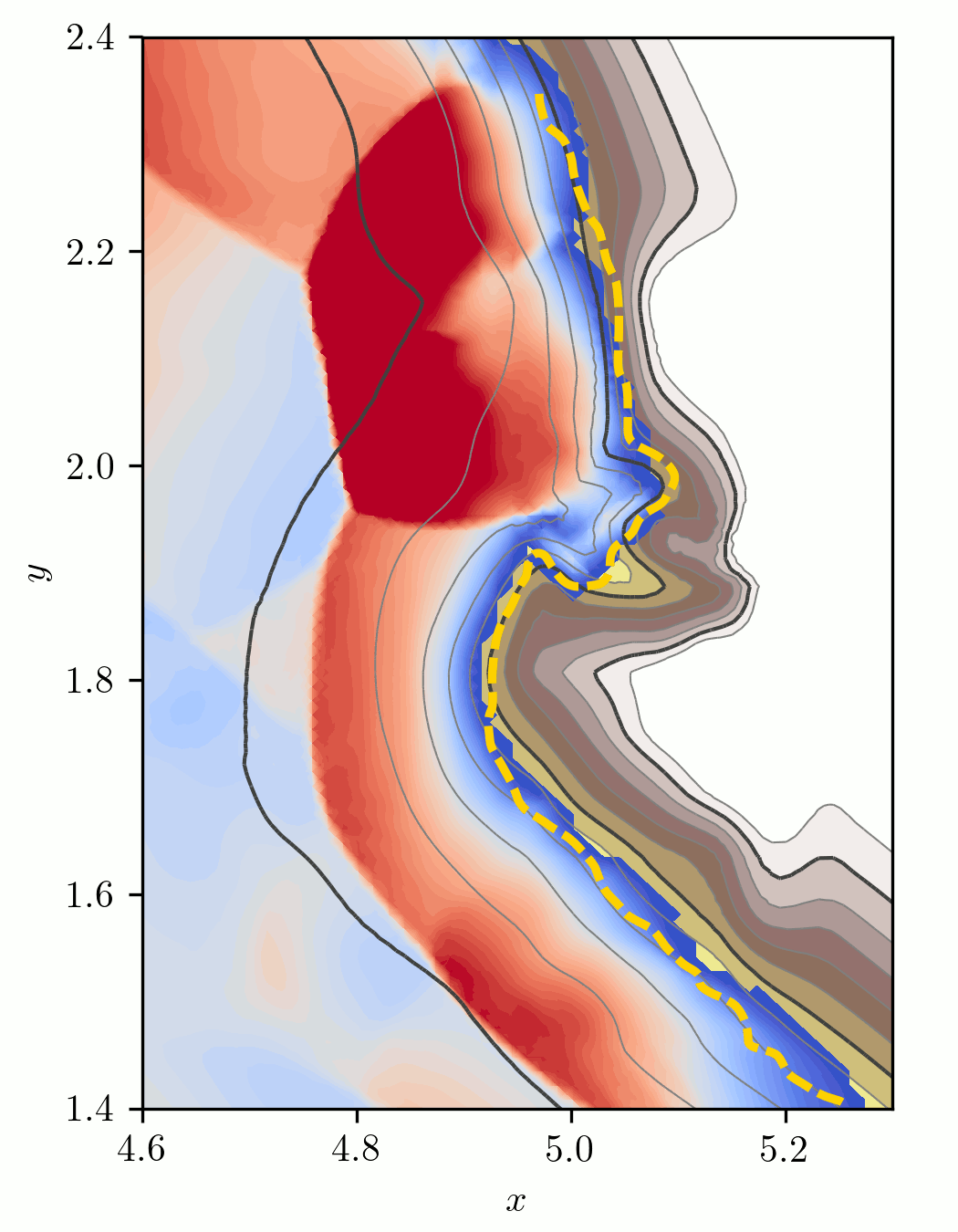}\\
  \includegraphics[width=0.33\textwidth]{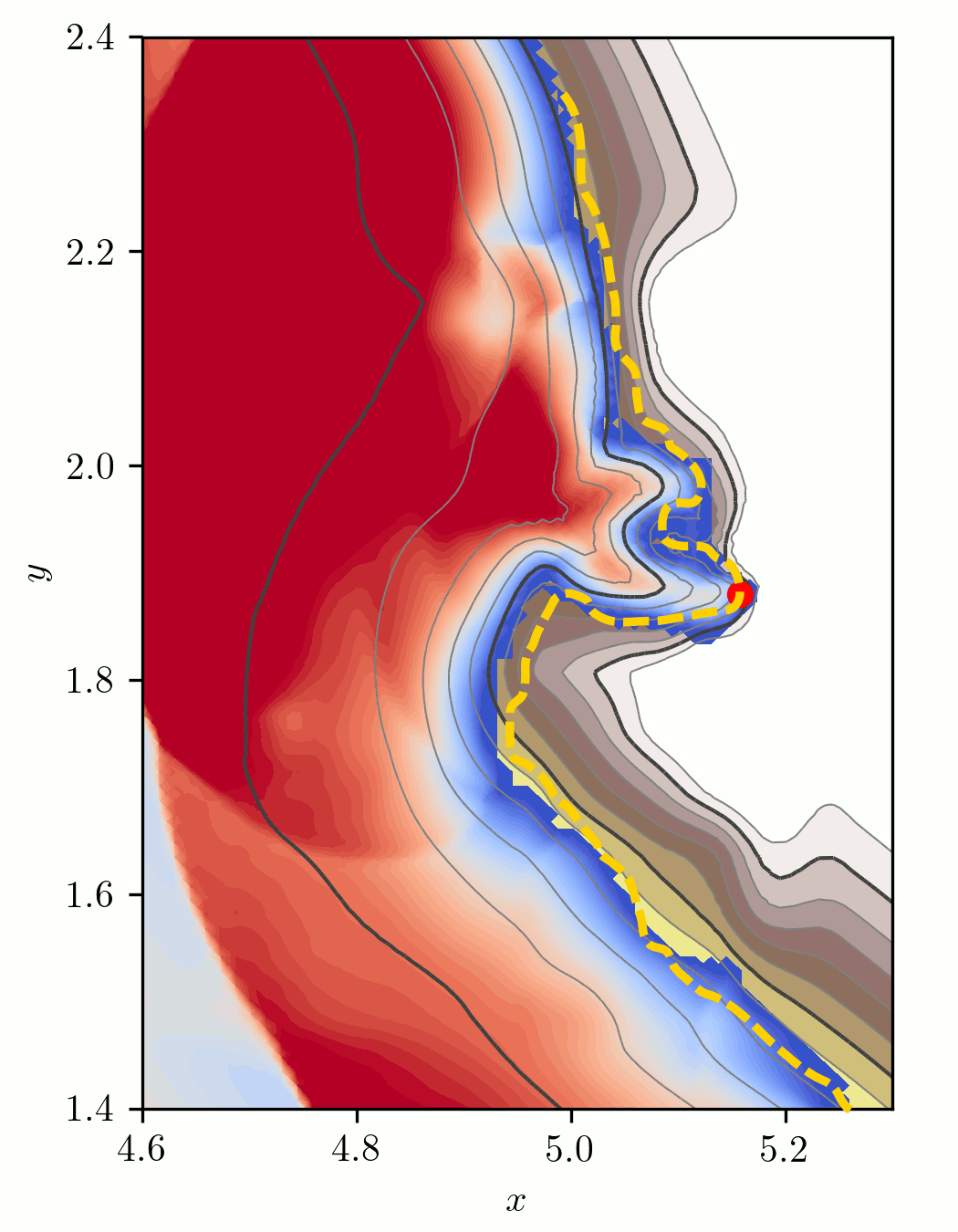}%
  \includegraphics[width=0.33\textwidth]{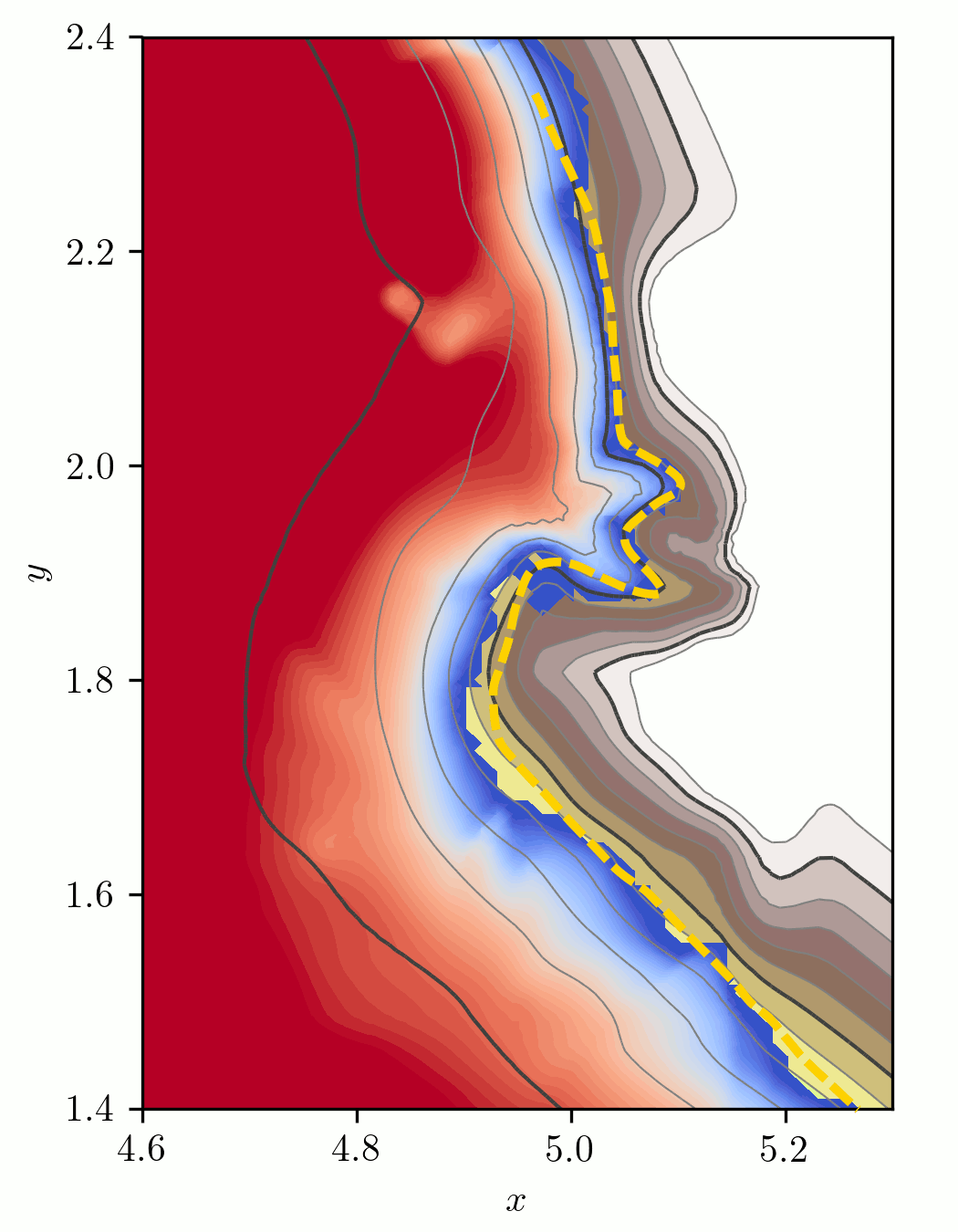}%
  \includegraphics[width=0.33\textwidth]{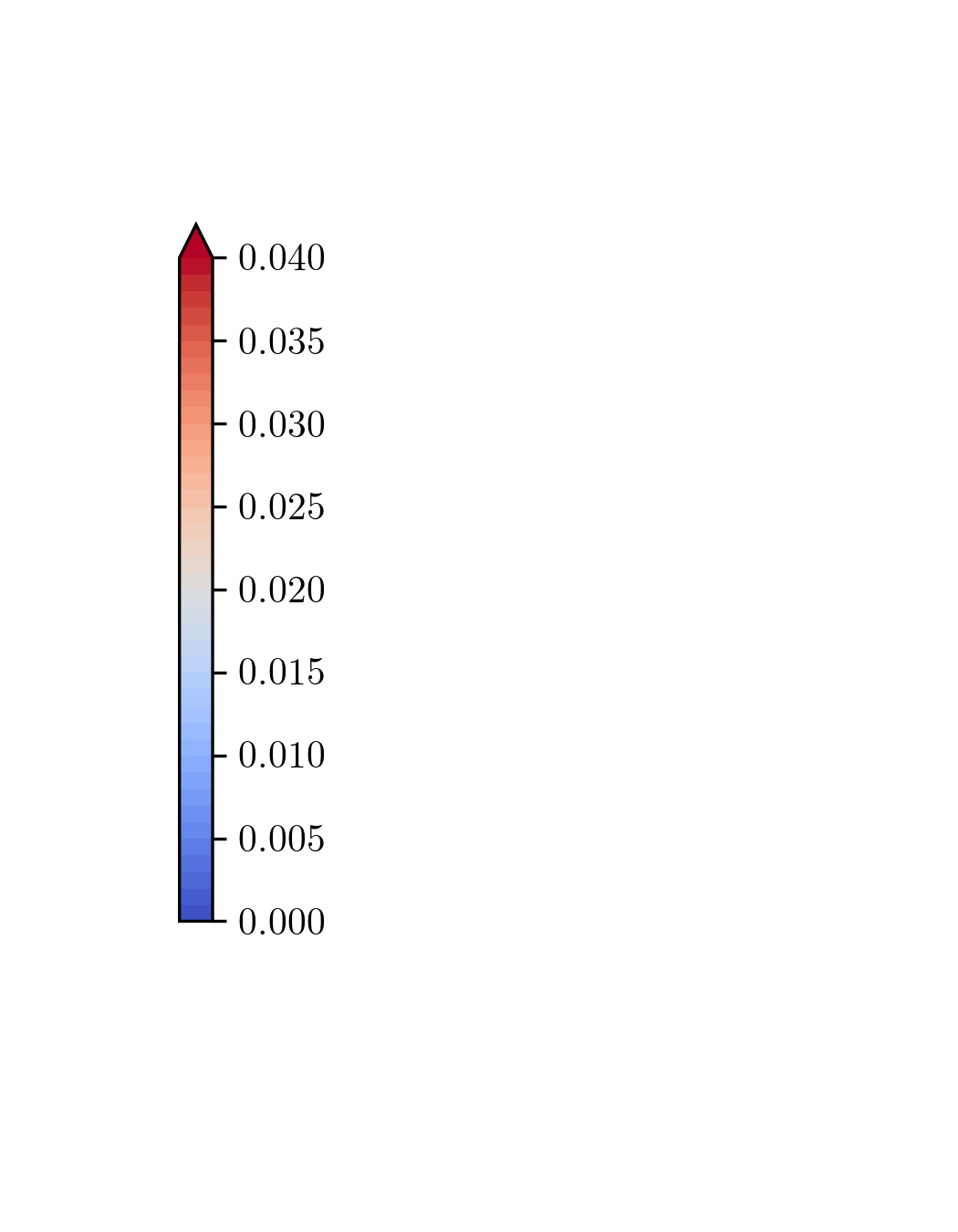}
  \caption{Okushiri: detailed contour plot of the coastal area at $t = 15.0, 15.5, 16.0, 16.5, 17.0$, (top left to right bottom), vertex-based limiter with $\tolwet= 10^{-4}$. Contour colors are given for the fluid depth. Contour lines represent the topography at 0.0, 0.01, \dots, 0.11. Also given is the approximate location of the shoreline from the experiment (yellow dashed) and the maximum observed runup (red circle), which happened around $t=16.5$.\label{fig:okushiri_contours}}
\end{figure}

We perform simulations with a time step of $\Delta t = 0.001$ until $40\,000$ steps ($\Tend = 40$) on a grid with 393\,216 elements ($384 \times 256$ rectangles divided into four triangles). The wet/dry tolerance is set to $\tolwet=10^{-4}$. The results are depicted in figures \ref{fig:okushiri_timeseries} and \ref{fig:okushiri_contours}. Figure \ref{fig:okushiri_timeseries} shows the comparison of the numerical results with experimental data at gauges 5, 7, and 9. Overall, we observe good agreement with both limiters (red and blue lines). Detailed contour plots of the coastal area together with the experimentally derived shoreline are shown in figure \ref{fig:okushiri_contours} for times $t = 15.0, 15.5, 16.0, 16.5, 17.0$. This shoreline is taken from \citep{NTHMP2011} and adjusted to the figures, which means it can only be used as a rough estimate. However, the flood line is represented well and we also demonstrate a good match of the maximum run-up (red dot) at $t=16.5$.

\subsection{Flow around a conical island}

\begin{figure}
  \centering
  \includegraphics[width=0.4\textwidth]{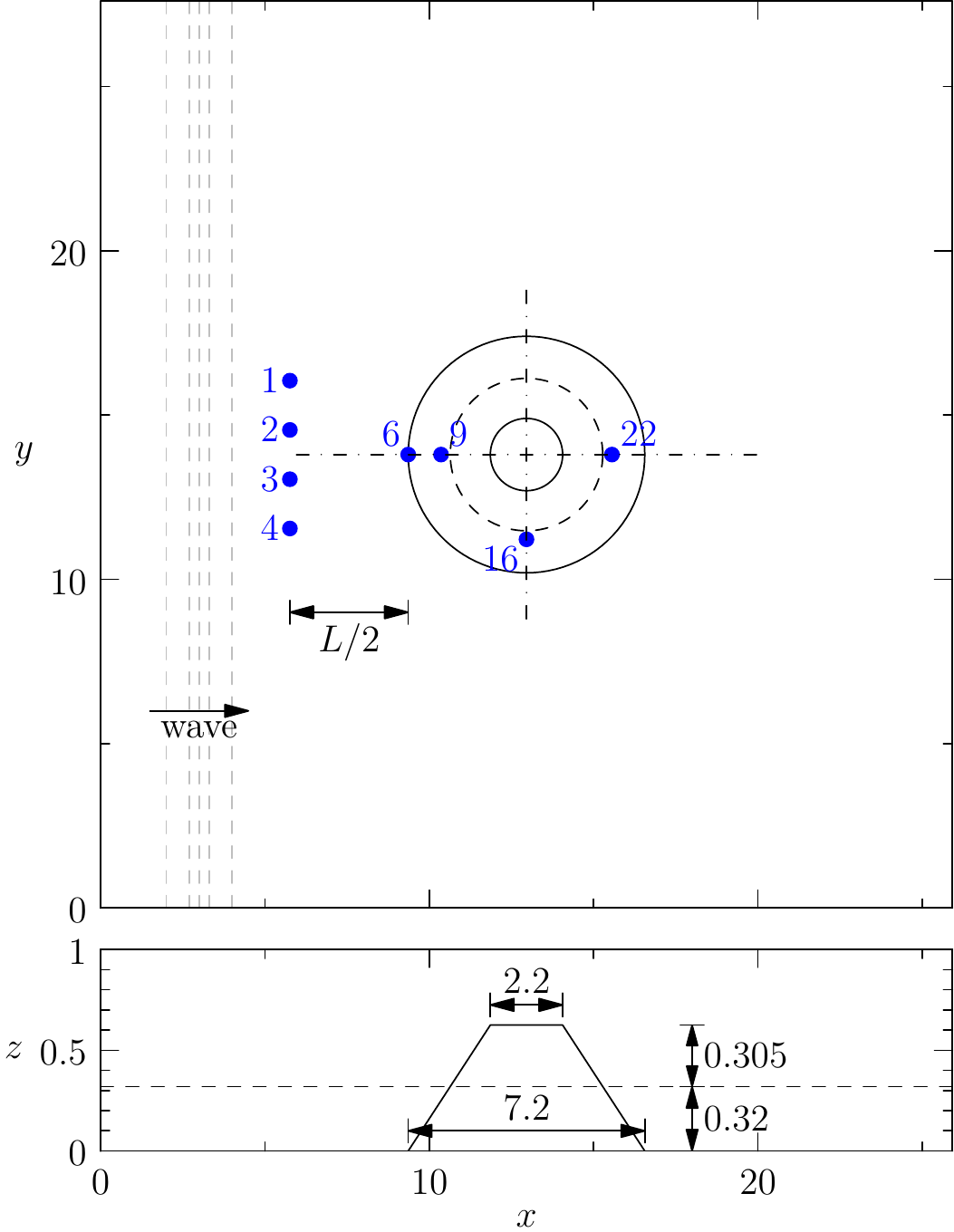}
  \caption{Conical Island: Top and side views of experimental setup with location of wave gauges (blue).\label{fig:conical_setup}}
\end{figure}

This test is part of a series of experiments carried out at the US Army Engineer Waterways Experiment Station in a $25\times28.2$m basin with a conical island situated at its center \citep{Briggs1995,Liu1995}. The experiment was motivated by the 1992 Flores Island tsunami run-up on Babi Island.
The conical island has its center at $\Cx_I = (12.96, 13.8)^\top$ and is defined by
\begin{align*}
b (\Cx) = \tilde{b}(r) = \begin{cases}
        0.625 & r \leq 1.1 \\
        (3.6-r)/4 & 1.1 \leq r \leq 3.6 \\
        0 & \text{otherwise,}
       \end{cases}
\end{align*}
with $r = \| \Cx - \Cx_I\|_2$ being the distance from the center (see figure \ref{fig:conical_setup}). The initial fluid depth and velocity are given by $h(\Cx,0) = \max\{0, h_0 - b(\Cx)\}$ and $\vu(\Cx,0) = \fatvec{0}$, where $h_0 = 0.32$. Three different solitary waves (denoted by case A, B and C) were generated by a wavemaker in the experiments at the left side of the domain, from which we only consider case A and C. Besides the trajectories of the wave paddle, time series of the surface elevation at 27 different gauge stations, 8 of which are freely available, were measured. The first four gauge stations were situated in a $L/2$ distance in $x$-direction from the toe of the beach, where $L/2$ is the distance at which the solitary wave height drops to $5\%$ of its maximum height, and $L$ defines the wave length. In the numerical simulations we describe the wave by an incoming analytical solitary wave through the boundary condition on the left side of the domain. In order to make the analytical wave compatible to measurements, it needs to be adjusted with the parameters given below. The wave is defined by
\begin{equation*}
  h_b(t) = h_0 + a \left( \frac{1}{\cosh(K (cT - ct - x_0))} \right)^{\! 2}
\end{equation*}
where $K = \sqrt{\tfrac{3 a}{4 h_0^3}}$ and $c = \sqrt{g h_0} \left(1+\tfrac{a}{2 h_0}\right)$. To obtain the other parameters, we fitted the solitary wave to the experimental data at the first four gauge stations. This resulted in an amplitude and time shift of $a = 0.014$ and $T = 8.85$ for case A. The parameter $x_0 = 5.76$ is the $x$-coordinate of the first four gauges. For case C the parameters are $a = 0.057$, $T = 7.77$ and $x_0 = 7.56$. Compared to the experiments this also includes a time shift of 20. As in the Okushiri test case the velocity at the boundary is defined to obtain a right running simple wave (cf.\ \eqref{eq:bcsimplewave}).

For the numerical simulations the domain is slightly adjusted to have dimensions $25.92 \times 27.60$, and the conical island is exactly centered. The domain is discretized into $1024 \times 1024$ uniform squares which are divided into two triangles (2\,097\,152 elements). The time step is $\Delta t = 0.0025$, and the computations are run until $\Tend = 20$. Results are computed using both limiters with a wet/dry tolerance of $\tolwet= 10^{-3}$. On the left side of the domain, we impose an inflow boundary condition to prescribe the solitary wave. Furthermore, a transparent boundary condition is set on the right side and wall boundary conditions are set at the top and the bottom of the domain.

\begin{figure}
  \includegraphics[width=0.5\textwidth]{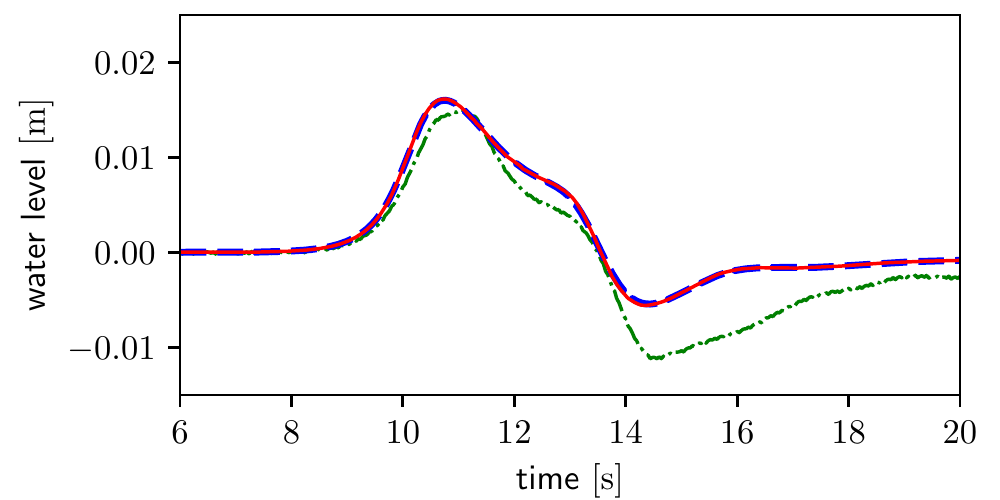}\hfill
  \includegraphics[width=0.5\textwidth]{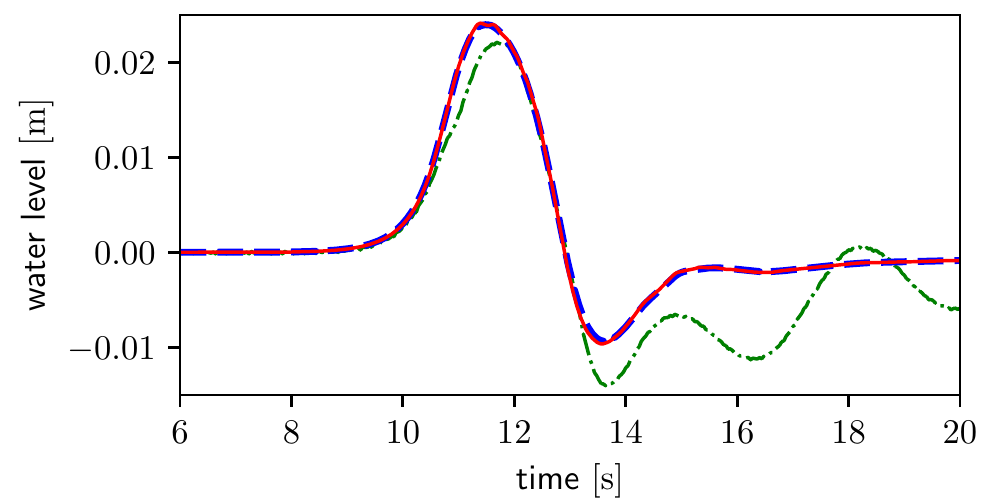}\\
  \includegraphics[width=0.5\textwidth]{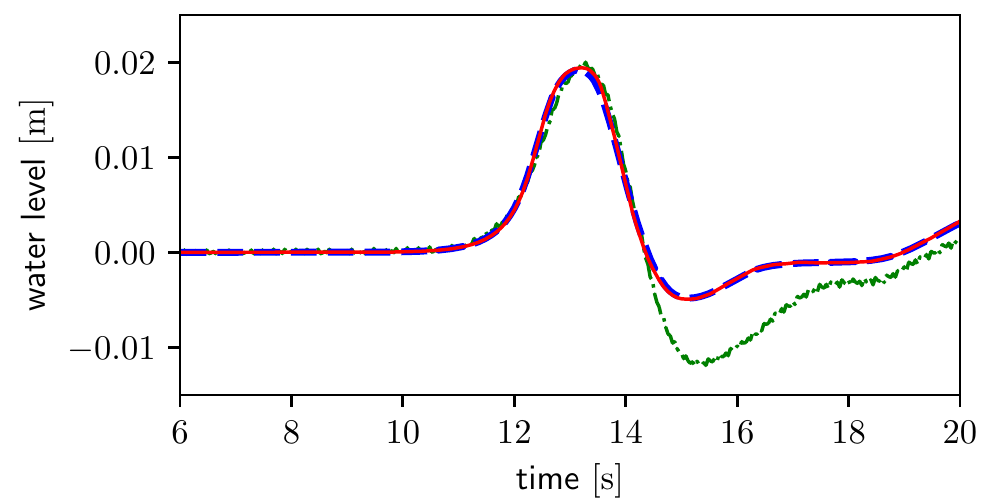}\hfill
  \includegraphics[width=0.5\textwidth]{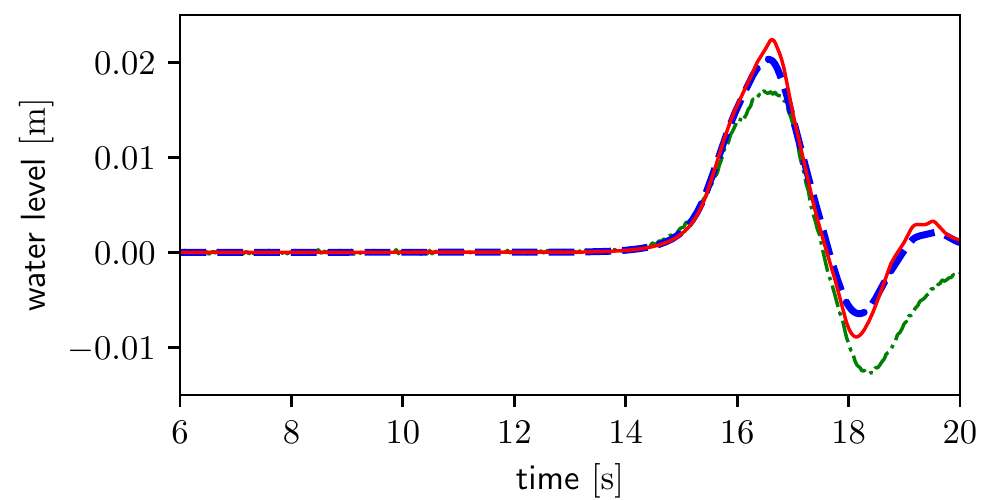}
  \caption{Conical Island: time series for case A of surface elevation at wave gauges 6, 9, 16, 22 (top left to bottom right) experimental data (green dash-dotted, shifted by $\Delta t = -20$), vertex-based limiter (red), edge-based limiter (blue dashed).\label{fig:conical_timeseriesA}}
\end{figure}

\begin{figure}
  \includegraphics[width=0.5\textwidth]{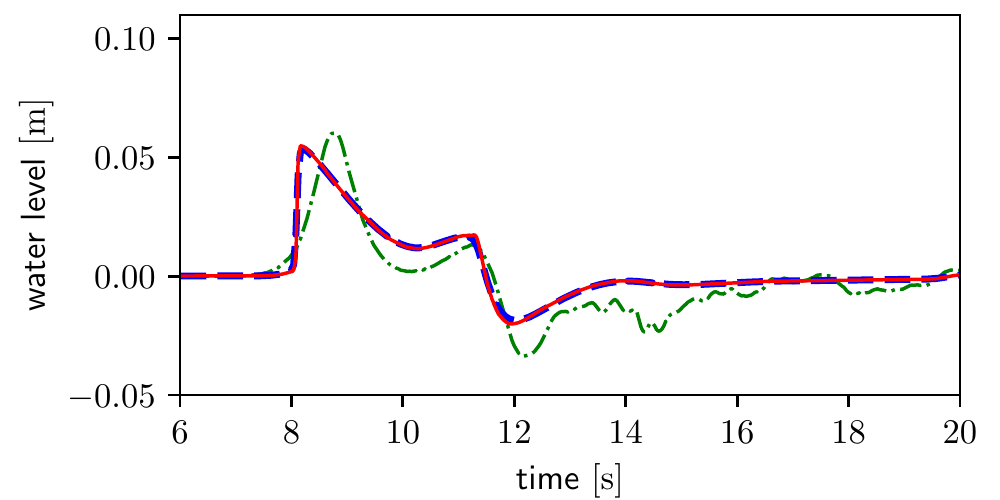}\hfill
  \includegraphics[width=0.5\textwidth]{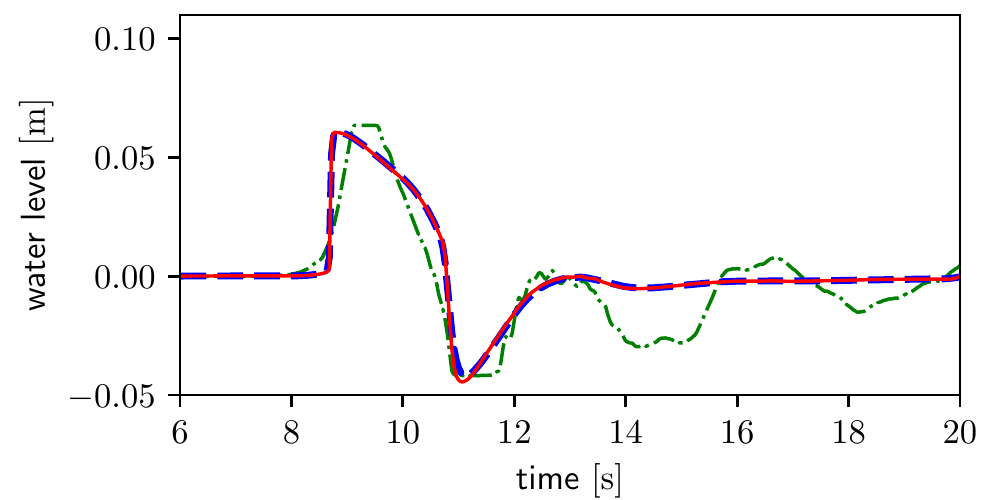}\\
  \includegraphics[width=0.5\textwidth]{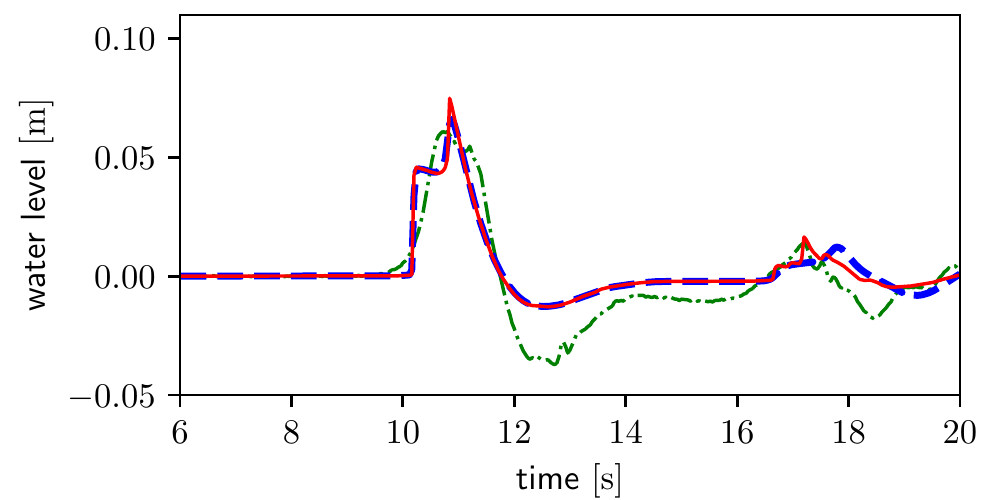}\hfill
  \includegraphics[width=0.5\textwidth]{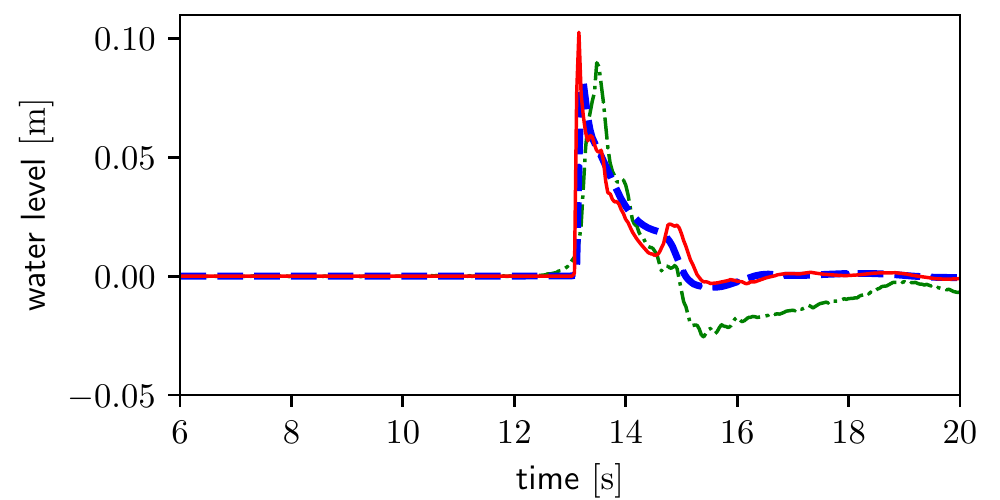}
  \caption{Conical Island: time series for case C of surface elevation at wave gauges 6, 9, 16, 22 (top left to bottom right) experimental data (green dash-dotted, shifted by $\Delta t = -20$), vertex-based limiter (red), edge-based limiter (blue dashed).\label{fig:conical_timeseriesC}}
\end{figure}

In figures \ref{fig:conical_timeseriesA} and \ref{fig:conical_timeseriesC} we compare the resulting time series of the surface elevation at gauge stations 6, 9, 16, and 22 for case A and C, respectively, with the experimental data. Note that some time series from the experiments were slightly shifted to have an initial zero water level. While gauge 6 and 9 are right in front of the island, gauge 16 is on the side and gauge 22 at the rear of it. It can be seen that for smaller wave amplitudes (case A) the experimental data can be reproduced well. On the other hand, for higher amplitudes in case C non-linear effects become dominant and are not balanced because of the lack of wave dispersion. The result is a steepening of the wave at the front and a flattening at the rear. Furthermore, we observe a general under-estimation of the trough after the first wave.

\begin{figure}[t!]
  \includegraphics[width=\textwidth]{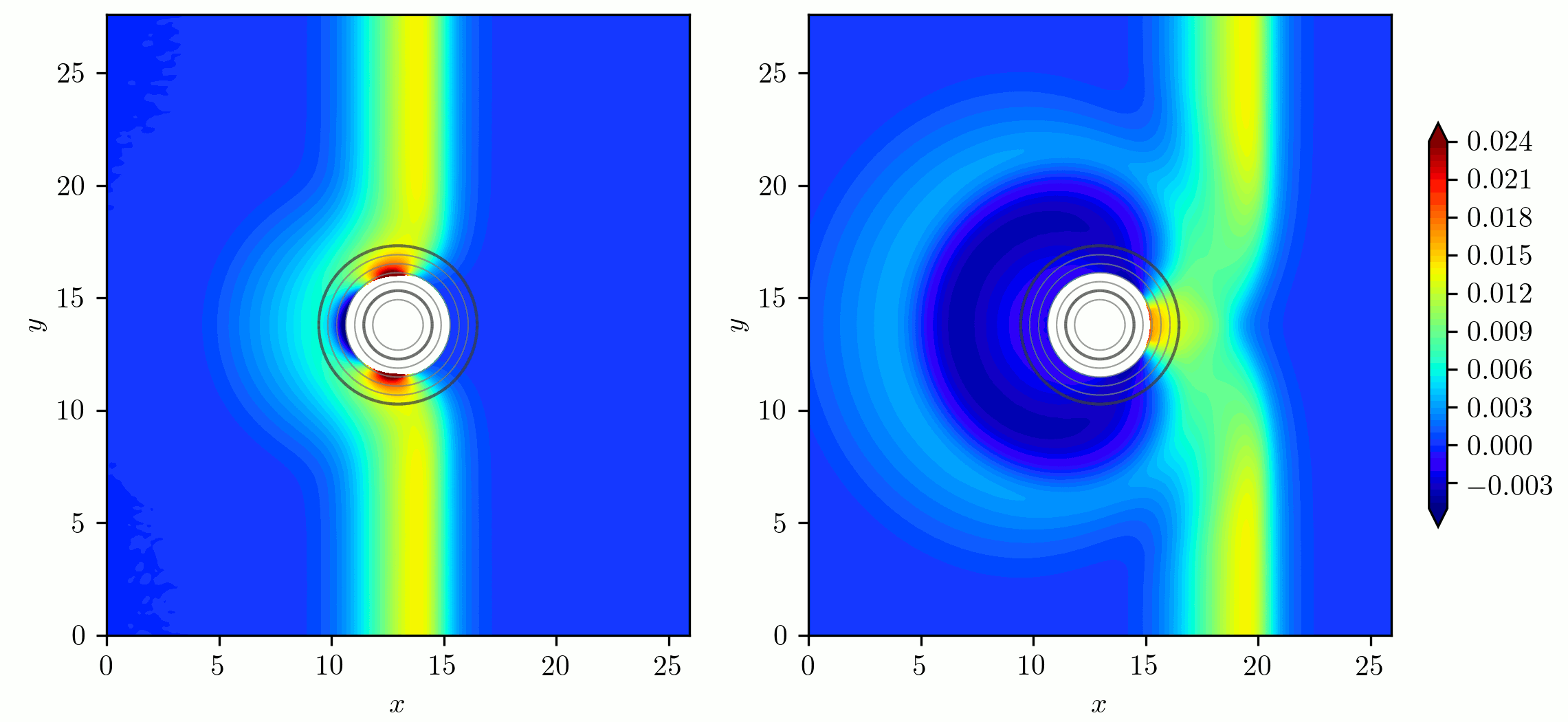}\\
  \includegraphics[width=\textwidth]{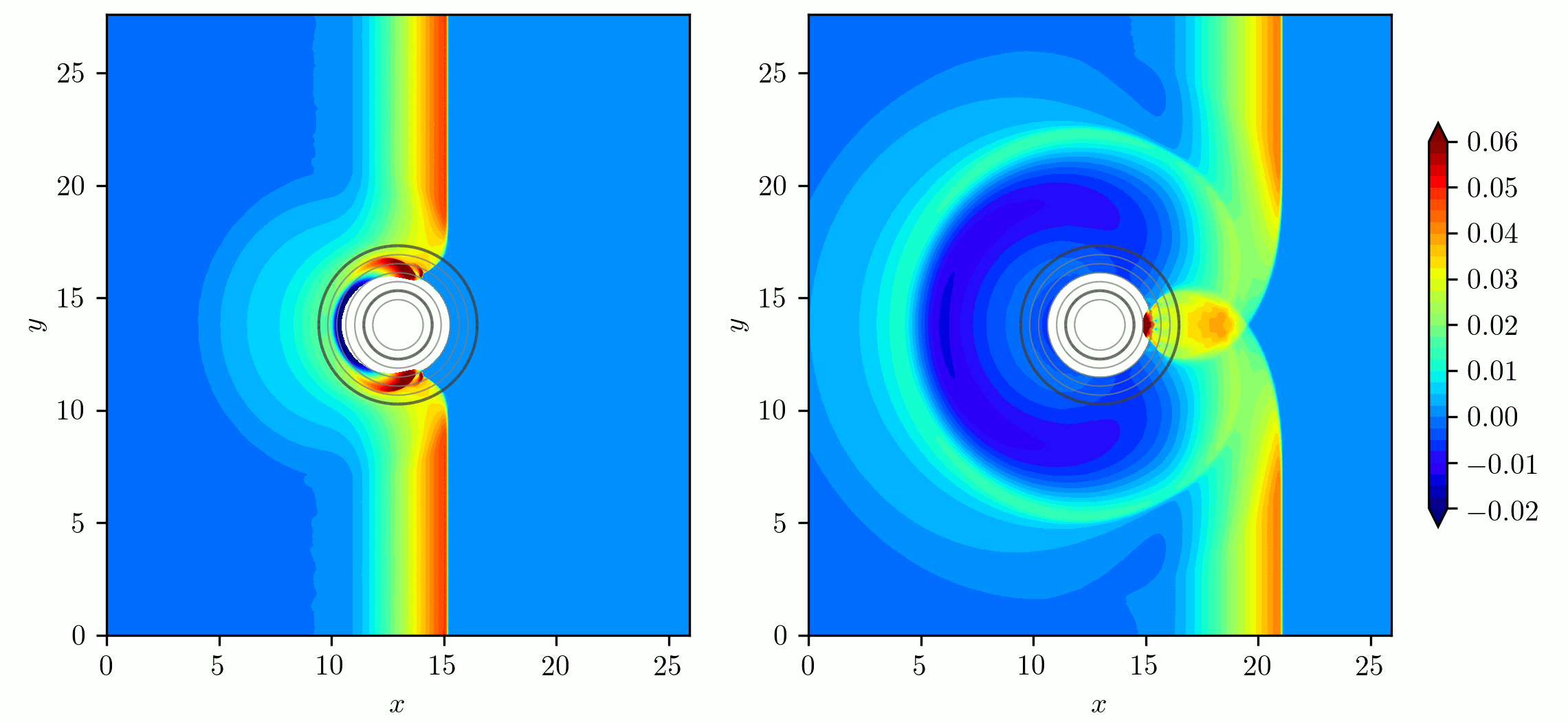}
  \caption{Conical Island: contour plot of surface elevation for case A (top) at $t=13$ (left) and $t=16$ (right) using the vertex-based limiter. Same for case C (bottom) at $t=11$ (left) and $t=14$ (right).\label{fig:conical_snapshots}}
\end{figure}

\begin{figure}
  \includegraphics[width=0.5\textwidth]{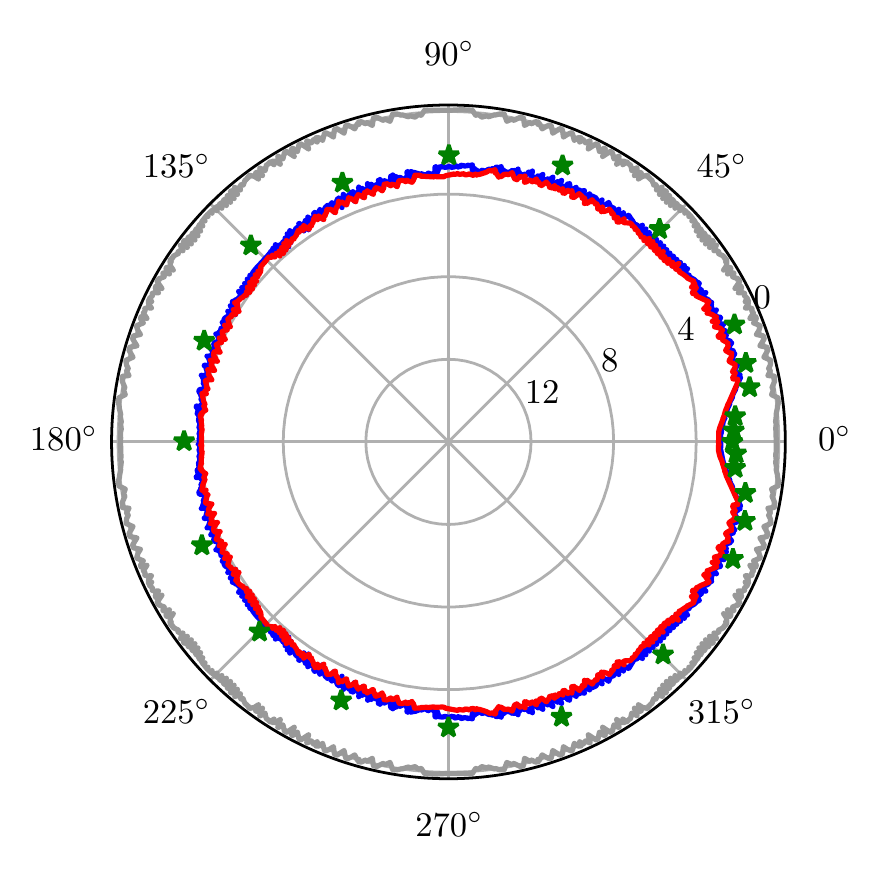}\hfill
  \includegraphics[width=0.5\textwidth]{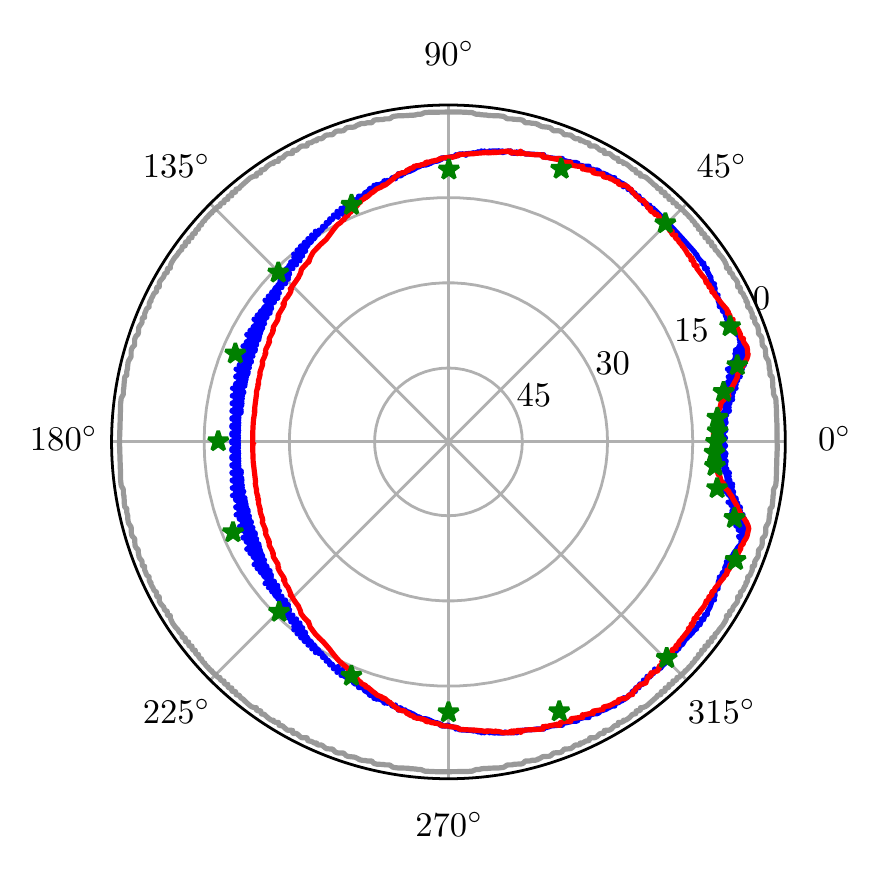}
  \caption{Conical Island: maximum vertical runup in $\mathrm{cm}$ for case A (left) and case C (right). Results using vertex-based limiter (red), edge-based limiter (blue), and from the experiments (green asterisks). Initial numerical shoreline is displayed in gray. Note the enlarged runup scale for case A.\label{fig:conical_maxrunups}}
\end{figure}

In figure \ref{fig:conical_snapshots} some snapshots of the simulations using the vertex-based limiter are displayed. They demonstrate that the initial wave correctly splits into two wave fronts after hitting the island. These wave fronts collide behind the island at a later time in the simulation. Finally, we compare the computed maximum run-up on the island with measurements from the experiments in figure \ref{fig:conical_maxrunups}. For both configurations of wave amplitude the simulations resemble measurements well and only slightly overestimate the runup. These deviations are larger in case C at the front of the island, where the wave first arrives. This behavior is probably due to the lack of wave dispersion and an imprecise representation of the wave generated by the wavemaker. Additional discrepancies might be related to the neglected bottom friction within the model. We attribute the better fit of the runup resulting from the simulation with the edge-based limiter to the additional diffusion introduced by this limiter, and not to a better physical modeling of the runup.

\section{Conclusions} \label{sec:Conclusions}

In this work a new wetting and drying treatment for RKDG2 methods applied to the shallow water equations is presented. The key  ingredients are a non-destructive limiting of the fluid depth combined with a velocity-based limiting of the momentum, the latter preconditioning the velocity computation near the shoreline, i.e., in areas of small fluid depth and momentum. This, in turn, guarantees a uniform time step with respect to the CFL stability constraint for explicit methods, which we explicitly report. The limiting strategy is complemented by a straightforward flux modification and a positivity-preserving limiter, which renders the scheme mass-conservative, well-balanced and stable for a wide range of flow regimes. It is a natural extension of a previously developed 1D scheme \citep{Vater2015} to the case of two-dimensional structured and unstructured triangular grids.

Originally designed to control linear stability, the chosen limiter for the fluid depth does not alter steady states at rest and small perturbations around them.
Two versions of the limiter are presented that differ in the selection of cells to be included into the limiting procedure. The ``edge-based'' version is based on the original Barth/Jespersen \citep{Barth1989} limiter. Due to its small stencil, it modifies states with constant gradients and therefore introduces additional diffusion into the method. The ``vertex-based'' version is an extension of the Barth/Jespersen limiter \citep{kuzmin_vertex-based_2010} especially designed for triangular grids and is non-destructive to linear states. It results in slightly more accurate computations in most situations, but has a larger stencil.

Only one single parameter $\tolwet$ enters the scheme, which controls the threshold in fluid depth considered to be dry. We show that the stability of the method is unaffected by this parameter. It solely determines the effective area which is considered wet by the discretization. A carefully chosen wet/dry tolerance thus leads to an accurate shoreline computation.

The presented test cases range from simple configurations where the analytical solution is known to the reproduction of laboratory experiments. They illustrate the method's applicability to a variety of flow regimes and verify its numerical properties: well-balancing in the case of a lake at rest, accurate representation of the shoreline, even in case of fast transitions, and convergence to the exact solution. Comparison with laboratory experiments shows good agreement. Some of the test cases are benchmark problems for the evaluation of operational tsunami models \citep{Synolakis2007}. With the successful simulation of these problems, we could show that the presented model satisfies the requirements for its application to realistic geophysical problems.

Future research will concentrate on the extension of the current scheme to adaptive grids and its application to tsunami and storm surge simulations. In this respect, additional source terms, such as the parametrization of sub-grid roughness by bottom friction and wind drag must be incorporated into the model. Furthermore, possibilities to extend the proposed concept to higher than second-order RKDG methods are investigated.

\section*{Acknowledgments}
This work benefited greatly from free software products. Without these tools -- such as python, the gnu FORTRAN compiler and the Linux operating system -- a lot of tasks would not have been so easy to realize. It is our pleasure to thank all developers for their excellent products.

\bibliographystyle{unsrtnat}
\bibliography{InundationDG2d}

%

\end{document}